# RANDOM WALKS ON SUPERCRITICAL PERCOLATION CLUSTERS[1]

By Martin T. Barlow

*University of British Columbia*

We obtain Gaussian upper and lower bounds on the transition density $q_t(x, y)$ of the continuous time simple random walk on a supercritical percolation cluster $\mathcal{C}_\infty$ in the Euclidean lattice. The bounds, analogous to Aronsen's bounds for uniformly elliptic divergence form diffusions, hold with constants $c_i$ depending only on $p$ (the percolation probability) and $d$. The irregular nature of the medium means that the bound for $q_t(x, \cdot)$ holds only for $t \geq S_x(\omega)$, where the constant $S_x(\omega)$ depends on the percolation configuration $\omega$.

**0. Introduction.** In this paper we study the simple random walk on the infinite component of supercritical bond percolation in the lattice $\mathbb{Z}^d$. We recall the definition of percolation [see Grimmett (1999)]: For edges $e = \{x, y\} \in \mathbb{E}_d = \{\{x, y\} : |x - y| = 1\}$, we have i.i.d. Bernoulli r.v. $\eta_e$, with $\mathbb{P}_p(\eta_e = 1) = p \in [0, 1]$, defined on a probability space $(\Omega, \mathcal{F}, \mathbb{P}_p)$. Edges $e$ with $\eta_e = 1$ are called *open* and the open cluster $\mathcal{C}(x)$ that contains $x$ is the set of $y$ such that $x$ and $y$ are connected by an open path. It is well known that there exists $p_c \in (0, 1)$ such that when $p > p_c$ there is a unique infinite open cluster, which we denote $\mathcal{C}_\infty = \mathcal{C}_\infty(\omega)$.

For each $\omega$ let $Y = (Y_t, t \geq 0, P_\omega^x, x \in \mathcal{C}_\infty)$ be the continuous time simple random walk (CTSRW) on $\mathcal{C}_\infty$; $Y$ is the process that waits an exponential mean 1 time at each vertex $x$ and then jumps along one of the open edges $e$ that contains $x$, with each edge chosen with equal probability. If we write $\nu_{xy}(\omega) = 1$ if $\{x, y\}$ is an open edge and 0 otherwise, and set $\mu(x) = \sum_y \nu_{xy}$, then $Y$ is the Markov process with generator

$$\mathcal{L}_\omega f(x) = \mu(x)^{-1} \sum_y \nu_{xy}(f(y) - f(x)), \qquad x \in \mathcal{C}_\infty. \tag{0.1}$$

Received January 2003; revised January 2004.
[1]Supported in part by a NSERC (Canada) grant, by CNRS (France) and by the Centre Borel (Paris).

A number of papers have studied this process or the closely related discrete time random walk $X = (X_n, n \geq 0, \widehat{P}_\omega^x, x \in \mathcal{C}_\infty)$. De Gennes (1976) discussed the link between the behavior of $X$ and resistance properties of $\mathcal{C}_\infty$, and coined the term "the ant in the labyrinth" to describe its motion. It is believed that, in all dimensions, there is no infinite cluster if $p = p_c$ (this is known to be true if $d = 2$ or $d \geq 19$). Kesten (1986a) defined the *incipient infinite cluster* when $d = 2$, and proved [Kesten (1986b)] that the random walk on this set has subdiffusive behavior.

For the supercritical case $p > p_c$, $Y$ is expected to have long time behavior similar to the random walk on $\mathbb{Z}^d$, and this has been confirmed in several ways. De Masi, Ferrari, Goldstein and Wick (1989) proved an invariance principle $(d = 2)$, while Grimmett, Kesten and Zhang (1993) proved that $Y$ is recurrent if $d = 2$ and transient when $d \geq 3$.

More recent papers have studied the transition density of $Y$ with respect to $\mu$:

$$(0.2) \qquad q_t(x,y) = q_t^\omega(x,y) = P_\omega^x(Y_t = y)\mu(y)^{-1}.$$

Theorem 1.2 of Mathieu and Remy (2004) gives the (quenched) bound

$$q_t^\omega(0,y) \leq c_1(p,d)t^{-d/2}, \qquad t \geq T_0(\omega), y \in \mathcal{C}_\infty,$$

$\mathbb{P}_p$-a.s. on the set $\{\omega : 0 \in \mathcal{C}_\infty\}$, while a similar (annealed) estimate was given by Heicklen and Hoffman (2000), but with an extra logarithmic factor.

The main result of this paper is the following two-sided bound on $q_t$.

THEOREM 1. *Let $p > p_c$. There exists $\Omega_1 \subseteq \Omega$ with $\mathbb{P}_p(\Omega_1) = 1$ and r.v. $S_x, x \in \mathbb{Z}^d$, such that $S_x(\omega) < \infty$ for each $\omega \in \Omega_1$, $x \in \mathcal{C}_\infty(\omega)$. There exist constants $c_i = c_i(d,p)$ such that for $x, y \in \mathcal{C}_\infty(\omega)$, $t \geq 1$ with*

$$(0.3) \qquad S_x(\omega) \vee |x - y|_1 \leq t,$$

*the transition density $q_t^\omega(x,y)$ of $Y$ satisfies*

$$(0.4) \quad c_1 t^{-d/2} \exp(-c_2|x-y|_1^2/t) \leq q_t^\omega(x,y) \leq c_3 t^{-d/2} \exp(-c_4|x-y|_1^2/t).$$

REMARKS. 1. The usual graph distance on $\mathbb{Z}^d$ is denoted $|x - y|_1 = \sum_{i=1}^d |x_i - y_i|$.

2. The CTSRW on $\mathbb{Z}^d$ (i.e., $p = 1$) satisfies these bounds with $S_x \equiv 1$. For $|x - y|_1 > t$, we have bounds which, up to constants, depend only on the tail of the Poisson distribution and not on the geometry of $\mathcal{C}_\infty$; see Lemma 1.1.

3. If $G$ is any finite graph which can be embedded in $\mathbb{Z}^d$, then $\mathcal{C}_\infty$ contains infinitely many copies of $G$ (attached at one point to the rest of $\mathcal{C}_\infty$). Since (0.4) does not hold uniformly for all such graphs, it is clear that we cannot expect (0.4) for all $x, y, t$ with $|x - y|_1 \leq t$. This irregularity of $\mathcal{C}_\infty$ is taken care of by the random variable $S_x$; after an initial period of possible bad behavior, $q_t(x, \cdot)$ settles down to a distribution with Gaussian tails.

RANDOM WALKS AND PERCOLATION 3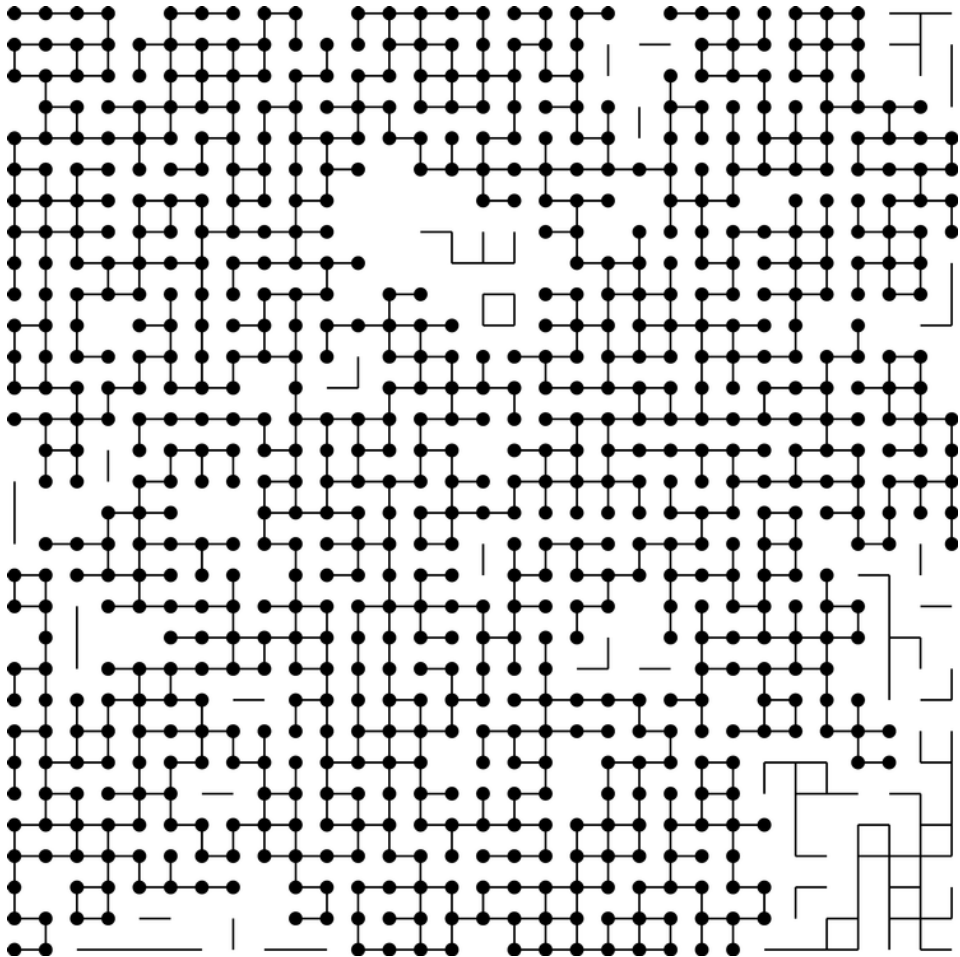

Fig. 1. *Bond percolation with $p = 0.53$. The vertices in the largest open cluster are marked with black circles.*

4. In (0.4) we can replace $|x - y|_1$ by the graph distance in $\mathcal{C}_\infty$ and, in fact, this is the result that we prove first (Proposition 6.1).

5. The constants $c_i$ are nonrandom, and depend only on $p$ and $d$. For $p$ sufficiently close to 1 it would, in principle, be possible to estimate them by careful tracking of the various constants in this paper. However, for general $p \in (p_c, 1)$ the constants arise in a noneffective fashion—if $p > p_c$, then we know that certain kinds of good behavior occur in cubes of side $k \geq k_0(p, d)$, but have no control on $k_0$. The constants $c_i$ then depend on $k_0$.

6. The tail of the random variable $S_x$ satisfies

(0.5) $$\mathbb{P}_p(x \in \mathcal{C}_\infty, S_x \geq n) \leq c \exp(-c' n^{\varepsilon_d});$$

see Lemma 2.24 and Section 6.

7. Similar bounds hold for the discrete time r.w. $X$. The proofs (which are not given here) run along the same lines, but with some extra (mainly minor) difficulties due to the discrete time. Whereas $\mathbb{Z}^d$ and $\mathcal{C}_\infty$ are bipartite we need to replace $q_t(x,y)$ in (0.4) with $p_{n-1}(x,y) + p_n(x,y)$.

8. It seems very likely that this theorem holds for other lattices in $\mathbb{Z}^d$ and that the "random walk" part of the proofs below transfers easily to other situations, but many of the percolation estimates that we use have been proved only for the Euclidean lattice.

Similar estimates hold in the annealed case.

THEOREM 2. *Let $p > p_c$. There exist constants $c_i = c_i(d,p)$ such that for $x, y \in \mathbb{Z}^d$, $t \geq 1$ with $|x - y|_1 \leq t$,*

$$
(0.6) \quad c_1 t^{-d/2} \exp(-c_2 |x-y|_1^2/t) \leq \mathbb{E}_p(q_t^\omega(x,y) | x, y \in \mathcal{C}_\infty) \\
\leq c_3 t^{-d/2} \exp(-c_4 |x-y|_1^2/t).
$$

An immediate consequence of Theorem 1 is that $Y$ is transient if and only if $d \geq 3$, but of course this was already known from Grimmett, Kesten and Zhang (1993). As an example of a new application, the off-diagonal bounds in Theorem 1 enable us to control harmonic functions on $\mathcal{C}_\infty$. Write $d_\omega(x,y)$ for the graph distance on $\mathcal{C}_\infty$ and let $B_\omega(x,R) = \{y : d_\omega(x,y) < R\}$ for balls. A function $h : B_\omega(x_0, R+1) \to \mathbb{R}$ is *harmonic* on $B_\omega(x,R)$ if $\mathcal{L}h(x') = 0$, $x' \in B_\omega(x,R)$. We have the following Harnack inequality.

THEOREM 3. *Let $p > p_c$. There exists $c_1 = c_1(p,d)$, $\Omega_1 \subseteq \Omega$ with $\mathbb{P}_p(\Omega_1) = 1$ and $R_0(x,\omega)$ such that $R_0(x,\omega) < \infty$ for each $\omega \in \Omega_1$, $x \in \mathcal{C}_\infty$. If $R \geq R_0(x,\omega)$ and $h : B_\omega(x, 2R+1) \to (0,\infty)$ is a positive harmonic function on $B_\omega(x, 2R)$, then writing $B = B_\omega(x,R)$,*

$$(0.7) \quad \sup_B h \leq c_1 \inf_B h.$$

This leads immediately to the Liouville property for positive harmonic functions.

THEOREM 4. (a) *Let $h : \mathcal{C}_\infty \to \mathbb{R}$ be positive and harmonic on $\mathcal{C}_\infty$. Then $h$ is constant.*

(b) *Let $\mathcal{T}$ denote the tail $\sigma$-field of $Y$. There exists $\Omega_2$ with $\mathbb{P}_p(\Omega_2) = 1$ such that for each $\omega \in \Omega_2$ and for each $F \in \mathcal{T}$, either $P_\omega^x(F) = 0$ for all $x \in \mathcal{C}_\infty(\omega)$ or else $P_\omega^x(F) = 1$ for all $x \in \mathcal{C}_\infty(\omega)$.*

REMARK. The Liouville property for bounded harmonic functions on $\mathcal{C}_\infty$ is already known; see Kaimanovitch (1990) and Lemma 4.6 of Benjamini, Lyons and Schramm (1999).

We can also use Theorem 1 to estimate $E_\omega^x |Y_t - x|^2$.

THEOREM 5. *Let $p > p_c$. There exists $\Omega_1 \subseteq \Omega$ with $\mathbb{P}_p(\Omega_1) = 1$ and r.v. $S_x', x \in \mathbb{Z}^d$, such that $S_x'(\omega) < \infty$ for each $\omega \in \Omega_1$, $x \in \mathcal{C}_\infty(\omega)$. There exist constants $0 < c_5 \leq c_6$ such that for $x \in \mathcal{C}_\infty(\omega)$, $t \geq S_x'$,*

$$(0.8) \qquad c_5 t \leq E_\omega^x |Y_t - x|^2 \leq c_6 t.$$

This result implies that the (annealed) invariance principle, proved in De Masi, Ferrari, Goldstein and Wick (1989) for $d = 2$, can be extended to $d \geq 3$. See Sidoravicius and Sznitman (2003) for a discussion of this and for a much more delicate quenched invariance principle when $d \geq 4$.

The proof of Theorem 1 breaks into two fairly distinct parts. First (Section 2) we prove suitable geometric and analytic properties of $\mathcal{C}_\infty$. Then (Sections 3–5) we use "heat kernel" techniques to obtain the estimate (0.4). These techniques originate in the work of De Giorgi, Moser and Nash on divergence form elliptic equations; more recently they have been employed to study random walks on graphs. While they have been very successful in a wide variety of algebraic and geometric contexts, this has almost always been in circumstances in which the same regularity condition holds for all balls of a given size $r$.

A guide to the kind of properties we need is given by the following theorem from Delmotte, which is a translation to graphs of results from Grigor'yan (1992) and Saloff-Coste (1992) on manifolds. (The version given here has been adapted to the CTSRW on a positive density subgraph $\mathcal{G}$ of $\mathbb{Z}^d$.)

THEOREM A [see Delmotte (1999)]. *Let $\mathcal{G} = (G, E)$ be a subgraph of $\mathbb{Z}^d$ with graph distance $d$. Let $c_0$, $c_1$ and $c_2$ be positive constants. Suppose that $\mathcal{G}$ satisfies the following two conditions.*

(a) *For all balls $B(x, R)$ in $G$,*

$$(V_d) \qquad c_0 R^d \leq \mu(B(x, R)) \leq c_1 R^d.$$

(b) *For any ball $B = B(x, R)$ and function $f : B \to \mathbb{R}$, writing $\bar{f}_B = \sum_{x \in B} f(x) \mu(x) / \mu(B)$,*

$$(PI) \qquad \sum_{x \in B} (f(x) - \bar{f}_B)^2 \mu(x) \leq c_2 R^2 \sum_{x \in B} \sum_{\substack{y \in B \\ y \sim x}} (f(x) - f(y))^2.$$

*Then the transition density $q_t(x, y)$ of $Y$ satisfies, for $t \geq d(x, y) \vee 1$,*

$$(0.9) \quad c_3 t^{-d/2} \exp(-c_4 d(x,y)^2 / t) \leq q_t(x, y) \leq c_5 t^{-d/2} \exp(-c_6 d(x,y)^2 / t).$$

The first condition is of regular volume growth of balls in $\mathcal{G}$, and can be replaced by a more general "volume doubling" condition. The second is a family of Poincaré or spectral gap inequalities for $\mathcal{G}$.

This theorem suggests that to prove Theorem 1, we should obtain volume growth and Poincaré inequalities for $\mathcal{C}_\infty$. Some results of this kind are in the literature; see Pisztora (1996) and Deuschel and Pisztora (1996) for volume growth estimates, and Benjamini and Mossel (2003) for a Poincaré inequality. These results show that for fixed $x \in \mathcal{C}_\infty$ the probability that ($V_d$) and (PI) hold for a ball $B_\omega(x, R)$ increases to 1 as $R \to \infty$.

However, on its own this is not enough to give Theorem 1. There are now several proofs of Theorem A [Delmotte (1999) used Moser iteration to prove a parabolic Harnack inequality], but all involve iterative methods or differential inequalities which use ($V_d$) and (PI) for many balls of different sizes. The exact definition of "good" and "very good" balls is given in Section 1.7, but roughly speaking we say a ball $B_\omega(y, r)$ is good if ($V_d$) and (PI) hold, and a ball $B_\omega(x, R)$ is very good if all subballs $B_\omega(y, r) \subseteq B_\omega(x, R)$ are good for $r \geq R^{1/(d+2)}$. We need to prove that all sufficiently large balls $B_\omega(x_0, R)$ (centered at a fixed $x_0$) are very good, and to do this we have to extend some of the estimates in the literature. This is done in Section 2. The estimates in Pisztora (1996) and Deuschel and Pisztora (1996) are enough for the volume growth bounds, but more work is needed for the Poincaré inequality. As in Benjamini and Mossel (2003), we prove this from an isoperimetric inequality, which was obtained by Benjamini and Mossel (2003) and Mathieu and Remy (2004). We use the methods of those papers, but the need for better control of the probabilities means that we have to rework some of these arguments to identify more precisely the set of percolation configurations $\omega$ for which a ball is good or very good. If $p$ is sufficiently close to 1, then a fairly direct counting argument [see Benjamini and Mossel (2003) and Mathieu and Remy (2004)] is all that is needed, but for general $p > p_c$, we have to use a renormalization argument, as in Antal and Pisztora (1996). In this paper we follow quite closely the approach of Mathieu and Remy (2004); there is a gap in the renormalization argument of Benjamini and Mossel (2003).

While percolation arguments generally use cubes in $\mathbb{Z}^d$, the heat kernel estimates work most naturally if we use the "chemical" or graph distance $d_\omega(x, y)$ on $\mathcal{C}_\infty$. We can compare these two metrics using the main theorem of Antal and Pisztora (1996).

Many of the methods used to derive (0.9) from ($V_d$) and (PI) are not suitable for the percolation context. For example, Saloff-Coste (1992) proved in that (PI) and ($V_d$) imply a Nash estimate: for $f : G \to \mathbb{R}$,

$$\text{(N)} \qquad \int |\nabla f|^2 \geq c \|f\|_2^{2+4/d} \|f\|_1^{-4/d},$$

and Carlen, Kusuoka and Stroock (1987) proved that (N) is equivalent to

(0.10) $\qquad q_t(x,y) \leq c' t^{-d/2} \qquad$ for all $x, y \in G, t \geq 1$.

However, since $\mathcal{C}_\infty$ contains copies of $\{0, \ldots, n\}$ for all $n$, (0.10) is clearly false for $\mathcal{C}_\infty$. More generally, we cannot use any method which relies on global Sobolev inequalities; it is necessary to use local methods. [One approach here might be that of rooted or anchored isoperimetric inequalities as in Thomassen (1992) or Benjamini, Lyons and Schramm (1999). However, the correct extension to Nash or Sobolev inequalities has not yet been made.]

In Section 3 we obtain an initial global upper bound on $q_t$ using the Poincaré inequality directly, following the approach of Kusuoka and Zhou (1992). We can then obtain the off-diagonal upper bound in (0.4) using a method of Nash (1958), Bass (2002) and Barlow and Bass (1989).

We use the method of Fabes and Stroock (1986), also based on ideas of Nash, to obtain a local lower bound; that is, for $q_t(x,y)$ if $d(x,y) \leq t^{1/2}$. However, a difficulty arises in extending this to prove (0.4) for points $x, y$ with $d_\omega(x,y) \approx t^{1-\varepsilon}$. The standard technique is chaining: using a sequence of small balls $B(z_i, r)$ that connect $x$ and $y$, and the Chapman–Kolmogorov equations. It turns out that we need to take $r \approx t/d_\omega(x,y)$, so we may need balls so small that we cannot be sure that they are very good. This problem is resolved by an additional percolation argument: for some fixed $r_1 = r_1(p,d) \gg 1$, we can show that the collection of good balls of size $r \geq r_1$ is large enough so that a suitable chain $[B(z_i, r), 0 \leq i \leq m]$ of very good balls exists with $z_0$ close to $x$ and $z_m$ close to $y$. (see Theorems 2.23 and 5.4). This argument needs renormalization techniques, even if $p$ is close to 1.

Section 1 contains a brief account of various known facts on random walks on graphs which are used in the rest of the paper. The percolation arguments are given in Section 2. Sections 3–5 were written for a general graph $\mathcal{G}$ that satisfies appropriate volume growth and Poincaré inequalities, and can be read independently of Section 2. Upper bounds on $q_t$ are obtained in Section 3. Section 4 proves a weighted Poincaré inequality from the unweighted (weak) Poincaré inequalities derived in Section 2, using methods of Saloff-Coste and Stroock (1991) and Jerison (1986). This is then used in Section 5 to obtain lower bounds on $q_t$. Section 6 then ties these results together and gives the proofs of Theorems 1–4.

We use $c_i$ to denote constants whose values are fixed within each argument and use $C.$ to denote constants fixed within each section; $c_{2.3.1}$ denotes the constant $c_1$ of Lemma 2.3, and $c$ and $c'$ are constants whose values may change on each appearance. The constants $C_i$, $c_i$, $c$ and $c'$ are always strictly positive. The notation $k = k(p,d)$ means that the constant $k$ depends only on $p$ and $d$.

**1. Graphs and random walks.** In this section we review some well-known facts concerning graphs, random walks, and isoperimetric Cheeger and Poincaré inequalities. Let $\mathcal{G} = (G, E)$ be an infinite, locally finite, connected graph. We define weights $\nu_{xy}$ by

$$\nu_{xy} = \begin{cases} 1, & \text{if } \{x, y\} \in E, \\ 0, & \text{otherwise,} \end{cases}$$

and set $\mu(x) = \sum_y \nu_{xy}$. We extend $\mu$ to a measure on $G$ and extend $\nu$ to a measure on $E$. Given a function $f : G \to \mathbb{R}$, we define $|\nabla f| : E \to \mathbb{R}$ by $|\nabla f|(e) = |f(x) - f(y)|$ if $e = \{x, y\}$. We write

$$\int f = \int f \, d\mu = \sum_{x \in V} f(x) \mu(x),$$

$$\int |\nabla f|^p = \int |\nabla f|^p \, d\nu = \sum_{e \in E} (|\nabla f|(e))^p \nu_e.$$

Let $d(x, y)$ be the graph distance on $G$ and define $B(x, r) = \{y : d(x, y) < r\}$. We assume we have a global upper bound on the size of balls: for any $x \in G$, $r \geq 1$,

(1.1) $$\mu(B(x, r)) \leq C_0 r^d.$$

Note that this implies that for each $x \in G$,

(1.2) $$1 \leq \mu(x) \leq C_0.$$

Let $Y = (Y_t, t \geq 0, P^x, x \in G)$ be the continuous time random walk on $G$; this is the Markov process with generator

(1.3) $$\mathcal{L}f(x) = \mu(x)^{-1} \sum_y \nu_{xy}(f(y) - f(x)).$$

Thus $Y$ waits at $x$ for an exponential mean 1 random time and then moves to a neighbor of $x$ at random. We define the transition density of $Y$ with respect to $\mu$ (or heat kernel density) by

$$q_t(x, y) = \mu(y)^{-1} P^x(Y_t = y).$$

Note that by (1.2) we have $q_t(x, y) \leq 1$ for all $x, y, t$. Given any points $x_0, \ldots, x_k \in G$ and times $t_1, \ldots, t_k$, then by the Markov property, if $t = \sum_k t_k$,

(1.4) $$q_t(x_0, x_k) \geq \mu(x_k)^{-1} \prod_{i=1}^k P^{x_{i-1}}(Y_{t_i} = x_i)$$
$$= \mu(x_k)^{-1} \prod_{i=1}^k q_{t_i}(x_{i-1}, x_i) \mu(x_i) \geq \prod_{i=1}^k q_{t_i}(x_{i-1}, x_i).$$

We begin by recalling some general bounds on $q_t$. These are not given in full generality, but just as they apply to the situation here.

RANDOM WALKS AND PERCOLATION 9

LEMMA 1.1. *Let $\mathcal{G}$ satisfy* (1.1).

(a) *There exist constants $c_i = c_i(d, C_0)$ such that, writing $D = d(x,y)$,*

$$(1.5) \qquad q_t(x,y) \leq c_1 \exp\left(\frac{-c_2 D^2}{t}\right), \qquad D \leq t,$$

$$(1.6) \qquad c\exp\left(-c_3 D\left(1 + \log \frac{D}{t}\right)\right) \leq q_t(x,y)$$
$$\leq c' \exp\left(-c_4 D\left(1 + \log \frac{D}{t}\right)\right), \qquad D \geq t \geq 1$$

*and*

$$(1.7) \qquad \frac{c_3}{(t \log t)^{d/2}} \leq q_t(x,x) \leq \frac{c_4}{t^{1/2}}, \qquad t \geq 1.$$

(b) *If $d(x,y) \geq R \geq 2$ and $t \leq c_5 R^2 / \log R$, then*

$$(1.8) \qquad q_t(x,y) \leq c_6 t^{-d}.$$

PROOF. (a) See Corollaries 11 and 12 of Davies (1993) for (1.5) and (1.6). For the discrete time random walk $X$ on $\mathcal{G}$, the lower bound in (1.7) is immediate from Theorem 2.2 of Coulhon and Grigor'yan (2003) and (1.1), while the upper bound follows from Theorem 2.3 of Coulhon and Grigor'yan (2003) and the fact that $\mu(B(x,r)) \geq r$ for all $r \geq 1$. These bounds then transfer to $q_t$ by integration.

(b) If $t < c_5 R^2 / \log R$, then $t \log t \leq 2c_5 R^2$ provided $R \geq c_7$. Hence $t^{-d} \geq \exp(-2dc_5 R^2/t)$ and, taking $c_5$ sufficiently small, the bound (1.8) is an easy consequence of (1.5). If $R < c_8$, then $t \leq c_9$ and (1.8) still holds on adjusting the constant $c_6$. □

We now review some geometric and analytic inequalities on $\mathcal{G}$. Let $H \subseteq G$ be finite, write

$$E(H) = \{e = \{x,y\} : x, y \in H\}$$

and call $\mathcal{H} = (H, E(H))$ the *induced subgraph* on $H$. We define the measures $\mu_0$ on $H$ and $\nu_0$ on $E(H)$ in the same way as $\mu$ and $\nu$ are defined for $\mathcal{G}$. Note that while $\nu_0$ and $\nu$ agree on $E(H)$, we have only

$$(1.9) \qquad \mu_0(A) \leq \mu(A) \leq C_0 \mu_0(A), \qquad A \subseteq H.$$

We now assume that $\mathcal{H}$ is connected. For $A_1, A_2 \subseteq H$ let

$$(1.10) \qquad \partial_E(A_1, A_2) = \{e = \{x,y\} : x \in A_1, y \in A_2\}.$$

Let

$$i(A) = \frac{\nu(\partial_E(A, H - A))}{\mu_0(A)}, \qquad A \subseteq H,$$

and define the isoperimetric constant

$$I_H = \min_{0 < \mu_0(A) \le \mu_0/2(H)} i(A).$$

Closely related is the Cheeger constant: Let

$$\chi(A) = \frac{\mu_0(H)\nu(\partial_E(A, H - A))}{\mu_0(A)\mu_0(H - A)}$$

and define

(1.11) $$J_H = \min_{A \ne \varnothing, H} \chi(A).$$

LEMMA 1.2. *The minimum in* (1.11) *is attained by a set $A$ such that $A$ and $H - A$ are connected.*

This is quite well known. For a proof, see, for example, Section 3.1 of Mathieu and Remy (2004).

LEMMA 1.3 [see Mathieu and Remy (2004), Section 3.1, and Benjamini and Mossel (2003)]. *Let $H$ be finite and connected.*

(a) *There exists $I_H \ge 2/\mu_0(H)$.*
(b) *If $I_H^* = \min\{i(A) : 0 < \mu(A) \le \tfrac{1}{2}\mu(H), A \text{ and } H - A \text{ are connected}\}$, then*

(1.12) $$I_H \le I_H^* \le 2I_H.$$

PROOF. (a) Let $0 < \mu_0(A) \le \tfrac{1}{2}\mu_0(H)$. Then since $H$ is connected, $\partial_E(A, H - A)$ is nonempty, so $i(A) \ge 1/\mu_0(A) \ge 2/\mu_0(H)$.

(b) The left-hand bound in (1.12) is obvious. If $0 < \mu_0(A) \le \tfrac{1}{2}\mu_0(H)$, then since $1 \le \mu_0(H)/\mu_0(H - A) \le 2$ we have $i(A) \le \chi(A) \le 2i(A)$. This immediately implies that $I_H \le J_H \le 2I_H$. Let $A$ be a minimal set for $J_H$. By Lemma 1.2 we can assume that $A$ and $H - A$ are connected. We can also take $\mu_0(A) \le \tfrac{1}{2}\mu_0(H)$. Then $I_H^* \le i(A) \le \chi(A) = J_H \le 2I_H$, proving the right-hand bound in (1.12). □

PROPOSITION 1.4. *Let $H \subseteq G$ be finite. Suppose that $I_H \ge \alpha^{-1}$.*

(a) *If $f : H \to \mathbb{R}$, then*

$$\min_a \int_H |f - a|^2 \, d\mu_0 \le c_1 \alpha^2 \int_{E(H)} |\nabla f|^2 \, d\nu.$$

(b) *If $f : H \to \mathbb{R}$, then*

(1.13) $$\min_a \int_H |f - a|^2 \, d\mu \le c_1 C_0 \alpha^2 \int_{E(H)} |\nabla f|^2 \, d\nu.$$

PROOF. (a) This result is well known. For a recent proof, see Lemma 3.3.7 of Saloff-Coste (1997).

(b) Since (1.9) can be used to replace $\mu_0$ with $\mu$, this follows immediately from (a). □

Inequality (1.13) is a *Poincaré* or *spectral gap* inequality for $H$. The minimum in (1.13) is of course attained by the value $a = \bar{f}_H = \int_H f\, d\mu/\mu(H)$. The following result is immediate from Lemma 1.3 and Proposition 1.4.

COROLLARY 1.5. *Let $H \subseteq G$ be finite and connected, and let $P_H$ be the best constant in the Poincaré inequality* (1.13).

(a) *There exists $P_H \leq c\mu(H)^2$.*
(b) *If $i(A) \geq \alpha^{-1}$ for all $A \subseteq H$ such that $A$ and $H - A$ are connected, then $P_H \leq c\alpha^2$.*

We note the following discrete version of the Gauss–Green lemma.

LEMMA 1.6. *Let $f, g \in L^2(G, \mu)$. Then*

$$(1.14) \quad \sum_{x \in G} g(x)\mathcal{L}f(x)\mu(x) = \tfrac{1}{2} \sum_x \sum_y (f(x) - f(y))(g(x) - g(y))\nu_{xy}.$$

In the sequel we need the following definitions.

DEFINITION 1.7. Let $C_V$, $C_P$ and $C_W \geq 1$ be fixed constants. We say $B(x, r)$ is $(C_V, C_P, C_W)$-*good* if

$$(1.15) \quad C_V r^d \leq \mu(B(x, r))$$

and the weak Poincaré inequality

$$(1.16) \quad \int_{B(x,r)} (f - \bar{f}_{B(x,r)})^2\, d\mu \leq C_P r^2 \int_{E(B(C_W x, r))} |\nabla f|^2\, d\nu$$

holds for every $f : B(x, C_W r) \to \mathbb{R}$.

We say $B(x, R)$ is $(C_V, C_P, C_W)$-*very good* if there exists $N_B = N_{B(x,R)} \leq R^{1/(d+2)}$ such that $B(y, r)$ is good whenever $B(y, r) \subseteq B(x, R)$, and $N_B \leq r \leq R$. We can always assume that $N_B \geq 1$. Usually the values of $C_V, C_P$ and $C_W$ are clear from the context and we just use the terms "good" and "very good."

**2. Percolation estimates.** We work with both bond and site percolation in $\mathbb{Z}^d$. We regard $\mathbb{Z}^d$ as a graph, with edge set $\mathbb{E}_d = \{\{x,y\} : |x-y| = 1\}$ and write $x \sim y$ to mean $\{x,y\} \in \mathbb{E}_d$. Given $Q \subseteq \mathbb{Z}^d$, define the internal and external boundaries of $A \subseteq Q$ by

$$\partial_i(A|Q) = \{y \in A : y \sim x \text{ for some } x \in A^c \cap Q\},$$
$$\partial_e(A|Q) = \{y \in A^c \cap Q : y \sim x \text{ for some } x \in A\} = \partial_i(Q - A|Q).$$

We begin with the notation for site percolation. Let $q \in (0,1)$ and $\Omega_s = \{0,1\}^{\mathbb{Z}^d}$, and define the coordinate maps $\zeta_x(\omega) = \omega(x)$. Let $\mathbb{Q}_q$ be the probability measure on $\Omega_s$ which makes the $\zeta_x$ i.i.d. Bernoulli r.v. with $\mathbb{Q}_q(\zeta_x = 1) = q$. We call those $x$ such that $\zeta_x = 1$ the *open sites* and write $\mathcal{O} = \mathcal{O}(\omega) = \{x : \zeta_x = 1\}$.

For $A \subseteq \mathbb{Z}^d$ we define the graph distance $d_A(x,y)$ to be the smallest $k$ such that there exists a path $\gamma = \{x_0, x_1, \ldots, x_k\} \subseteq A$ with $x = x_0$, $x_k = y$ and $x_{i-1} \sim x_i$, $1 \leq i \leq k$. If there is no such $k$, then $d_A(x,y) = \infty$. [We have $d_A(x,x) = \infty$ if $x \notin A$.] We write $d_\omega(x,y) = d_{\mathcal{O}(\omega)}(x,y)$ and refer to a path $x = x_0, x_1, \ldots, x_k = y$ such that each $x_i$ is open as an *open path*. We say $A$ is *connected* if $d_A(x,y) < \infty$ for all $x, y \in A$.

Now write

$$\mathcal{C}(x) = \{y : d_\omega(x,y) < \infty\}$$

for the connected open cluster that contains $x$. Write also, given $Q \subseteq \mathbb{Z}^d$,

$$\mathcal{C}_Q(x) = \{y : d_{Q \cap \mathcal{O}(\omega)}(x,y) < \infty\}.$$

This is the set of points connected to $x$ by an open path within $Q$. We call sets of the form $\mathcal{C}(x)$ *open clusters* and call the sets $\mathcal{C}_Q(x)$ *open $Q$ clusters*. It is well known that there exists $q_c = q_c(d) \in (0,1)$ such that if $q > q_c$, then $\mathbb{Q}_q$-a.s. there is a unique infinite open cluster $\mathcal{C}_\infty$. However, for site percolation we are interested only in the case when $q$ is either close to 1 or close to 0. Given a cube $Q \subseteq \mathbb{Z}^d$, the set $Q \cap \mathcal{C}_\infty$ in general is not connected. We write

$$\mathcal{C}^\vee(Q) \text{ for the largest open } Q \text{ cluster.}$$

(We adopt some procedures for breaking ties.)

Write $|x - y|_\infty = \max_i |x_i - y_i|$, let $\mathbb{E}_d^* = \{\{x,y\} : |x-y|_\infty = 1\}$ and write $x \overset{*}{\sim} y$ if $\{x,y\} \in \mathbb{E}_d^*$. We also need to consider site percolation in the graph $(\mathbb{Z}^d, \mathbb{E}_d^*)$. We say that $A \subseteq \mathbb{Z}^d$ is *-connected if $A$ is connected in the graph $(\mathbb{Z}^d, \mathbb{E}_d^*)$ and we define the clusters $\mathcal{C}^*(x)$ analogously.

DEFINITION 2.1. 1. Let $Q$ be a cube of side $n$ in $\mathbb{Z}^d$. We write $s(Q) = n$ for the side length of $Q$. Let $Q^+ = A_1 \cap \mathbb{Z}^d$ and $Q^\oplus = A_2 \cap \mathbb{Z}^d$, where $A_1$ and $A_2$ are the cubes in $\mathbb{R}^d$ with the same center as $Q$ and with sides $\frac{3}{2}n$ and $\frac{6}{5}n$, respectively.

2. A cluster $\mathcal{C}$ in a cube $Q$ is *crossing* for a cube $Q' \subseteq Q$ if for all $d$ directions there exists an open path in $\mathcal{C} \cap Q'$ that connects the two opposing faces of $Q'$.

3. The *diameter* of a set $A$ is defined by $\text{diam}(A) = \max\{|x-y|_\infty : x, y \in A\}$.

4. Given a set $A$, we write $|A|$ for the number of elements in $A$.

REMARK. In the arguments in this section, we frequently need to assume that a cube $Q$ is sufficiently large. More precisely, we need that $s(Q) \geq n_1$, where $n_1 = n_1(d, p)$ is a constant that depends only on $d$ and $p$. We make this assumption in our proofs whenever necessary without stating it explicitly each time. Unless otherwise indicated, the statements of the results are true for all $n$; this can be ensured by adjusting the constants so that the result is automatic for small cubes.

Let $Q$ be a cube in $\mathbb{Z}^d$. Define the event

$$(2.1) \quad K(Q, \lambda) = \left\{\omega : \mathcal{C}^\vee(Q) \text{ is crossing for } Q \text{ and } \frac{|\mathcal{C}^\vee(Q)|}{|Q|} > \lambda\right\}.$$

The following estimate was proved in Theorem 1.1 of Deuschel and Pisztora (1996).

LEMMA 2.2. *Let $Q$ be a cube of side $n$ and $\lambda < 1$. Then there exists $q_0 = q_0(d, \lambda) < 1$ and $c_i = c_i(\lambda, d)$ such that if $q \in [q_0, 1)$, then*

$$\mathbb{Q}_q(K(Q, \lambda)^c) \leq c_1 \exp(-c_2 n^{d-1}).$$

Let $\Sigma = \{x \in \mathbb{Z}^d : x \sim 0\} \cup \{0\}$. For $\sigma \in \Sigma$ define the shifted set of open sites $\mathcal{O}_\sigma = \{x - \sigma : x \in \mathcal{O}\}$. Let

$$(2.2) \quad \beta = 1 - 2(1+d)^{-1} < (d-1)/d.$$

Set for $r \in \mathbb{N}$ and $\varepsilon > 0$,

$$F(Q, r, \sigma, \varepsilon) = \{\text{any } *\text{-connected set } A \subseteq Q \text{ with } |A| = r$$
$$\text{satisfies } |A \cap \mathcal{O}_\sigma| \geq (1-\varepsilon)|A|\},$$
$$F(Q, \varepsilon) = \bigcap_{r \geq s(Q)^\beta} \bigcap_{\sigma \in \Sigma} F(Q, r, \sigma, \varepsilon).$$

Note that the event $F(Q, r, \sigma, \varepsilon)$ is increasing and that although $K(Q, \lambda)$ is not in general increasing, it is if $\lambda > \frac{1}{2}$.

Recall the following bounds on the tail of the binomial.

LEMMA 2.3. *Let $X \sim \text{Binomial}(n,p)$ and $\lambda \in (0,1)$. Then*
$$\mathbb{P}(X < \lambda n) \leq e^{-nb(\lambda,p)},$$
*where $b(\lambda, p) \to \infty$ as $p \to 1$ for each fixed $\lambda \in (0,1)$.*

LEMMA 2.4 [see Grimmett (1999), Section 4.2]. *The number of $*$-connected sets $A$ with $|A| = r$ containing a fixed point $x_0 \in \mathbb{Z}^d$ is bounded by $\exp(c_1 r)$.*

LEMMA 2.5. *Let $\varepsilon \in (0,1)$ and let $Q$ be a cube of side $n$. There exists $q_1 = q_1(\varepsilon, d) > q_c$ such that if $q \geq q_1$, then*
$$\mathbb{Q}_q(F(Q,\varepsilon)^c) \leq c\exp(-n^\beta).$$

PROOF. If $A$ is a fixed connected set in $Q$ with $|A| = r$, then by Lemma 2.3,
$$\mathbb{Q}_q(|A \cap \mathcal{O}_\sigma| < (1-\varepsilon)|A|) \leq e^{-rb(1-\varepsilon,q)}.$$
Given $\varepsilon$, choose $q_1$ large enough so that $b(\lambda, p) > 2 + c_{2.4.1}$ for $\lambda \geq 1 - \varepsilon$, $q \geq q_1$. Then since there are at most $(n+1)^d \exp(c_{2.4.1} r)$ $*$-connected sets in $Q$ of size $r$,
$$\mathbb{Q}_q(F(Q,r,\sigma,\varepsilon)^c) \leq (n+1)^d \exp(c_{2.4.1} r - rb(1-\varepsilon, q)) \leq n^d e^{-2r}$$
and, as $|\Sigma| = 2d + 1$,
$$\mathbb{Q}_q(F(Q,\varepsilon)^c) \leq \sum_{r=n_0}^\infty (2d+1)(n+1)^d e^{-2r}$$
$$\leq cn^d \exp(-2n^\beta) \leq c'\exp(-n^\beta). \quad \square$$

We collect from Deuschel and Pisztora (1996) the following results on the boundaries of discrete sets contained in cubes.

LEMMA 2.6. *Let $Q$ be a cube in $\mathbb{Z}^d$.*

(a) *Let $A \subsetneq Q$ be $*$-connected. Let $\Lambda_j$, $1 \leq j \leq k$, be the connected components of $Q - A$. Then $\partial_i(\Lambda_j|Q)$ and $\partial_e(\Lambda_j|Q)$, $1 \leq j \leq k$, are $*$-connected.*
(b) *Let $A \subseteq Q$ with $|A| \leq (15/16)|Q|$. Then*

(2.3) $\quad |A| \leq c_1 |\partial_i(A|Q)|^{d/(d-1)} \quad$ *and* $\quad |A| \leq c_1 |\partial_e(A|Q)|^{d/(d-1)}.$

PROOF. Part (a) is Lemma 2.1(ii) of Deuschel and Pisztora (1996), while the discrete isoperimetric inequality (2.3) is assertion (A.3) on page 480 of Deuschel and Pisztora (1996). $\square$

The following result is based on ideas in Mathieu and Remy (2004). Recall from (1.10) the definition of $\partial_E(A_1, A_2)$.

PROPOSITION 2.7. *Let $\varepsilon < 1/(4d+2)$ and $\lambda \geq \frac{7}{8}$. Suppose that both $F(Q,\varepsilon)$ and $K(Q,\lambda)$ occur for a cube $Q$ with side $n$. Let $A \subseteq Q$ be connected.*

(a) *If $A \cap \mathcal{C}^\vee(Q) = \varnothing$, then $|A| < cn^{\beta d/(d-1)} < n$.*
(b) *If $|A| \leq \frac{3}{4}|Q|$ and $A \cap \mathcal{C}^\vee(Q) \neq \varnothing$, then there exists $c_1$ such that*

$$(2.4) \qquad |\partial_E(A \cap \mathcal{C}^\vee(Q), (Q-A) \cap \mathcal{C}^\vee(Q))| \geq c_1 n^{-1}|A|.$$

PROOF. We write $\mathcal{C}^\vee = \mathcal{C}^\vee(Q)$.

(a) Let $A_0$ be the connected component of $Q - \mathcal{C}^\vee$ that contains $A$. By Lemma 2.6, $\partial_i(A_0|Q)$ is $*$-connected. The definition of $A_0$ implies that $\partial_i(A_0|Q) \cap \mathcal{O} = \varnothing$. Since $F(Q,\varepsilon)$ occurs, this implies that $|\partial_i(A_0|Q)| < n^\beta$ [since otherwise $|\partial_i(A_0|Q) \cap \mathcal{O}| > 0$]. Hence by the discrete isoperimetric inequality (2.3),

$$|A| \leq |A_0| \leq c|\partial_i(A_0|Q)| \leq cn^{\beta d/(d-1)} < n.$$

(b) Now let $A \cap \mathcal{C}^\vee \neq \varnothing$ and $|A| \leq n$. So there exists $x \in A \cap \mathcal{C}^\vee$. Since $|\mathcal{C}^\vee| \geq \frac{7}{8}n^d$, there exists $y \in (Q-A) \cap \mathcal{C}^\vee$ and whereas $\mathcal{C}^\vee$ is connected, there therefore exists $\{x',y'\} \in \partial_E(A \cap \mathcal{C}^\vee(Q), (Q-A) \cap \mathcal{C}^\vee(Q))$. So

$$|\partial_E(A \cap \mathcal{C}^\vee(Q), (Q-A) \cap \mathcal{C}^\vee(Q))| \geq 1 \geq n^{-1}|A|.$$

It remains to consider the case $A \cap \mathcal{C}^\vee \neq \varnothing$, $|A| > n$. Let

$$\Lambda_1 = \bigcup\{\mathcal{C}_Q(y): y \in \partial_e(A|Q) \cap \mathcal{O} - \mathcal{C}^\vee\}, \qquad A_1 = A \cup \Lambda_1.$$

Let $C_i$, $1 \leq i \leq k$, be the connected components of $Q - A_1$ and let

$$\Lambda_2 = \bigcup\{C_i: C_i \cap \mathcal{C}^\vee = \varnothing\}, \qquad A_2 = A_1 \cup \Lambda_2.$$

It is clear from the construction of $A_1$ and $A_2$ as a union of connected sets each connected to $A$ or $A_1$ that $A_2$ is connected. Note that since $\Lambda_1 \cup \Lambda_2 \subseteq Q - \mathcal{C}^\vee$ and $K(Q,\lambda)$ holds, we have $|A_2| \leq |A| + |Q - \mathcal{C}^\vee| \leq \frac{7}{8}|Q|$.

We now show that

$$(2.5) \qquad \partial_e(A|Q) \cap \mathcal{C}^\vee = \partial_e(A_2|Q) \cap \mathcal{O}.$$

First let $y \in \partial_e(A|Q) \cap \mathcal{C}^\vee$, so that $y \sim x$ with $x \in A$. Whereas $y \in \mathcal{C}^\vee$, we cannot have $y \in \Lambda_1 \cup \Lambda_2$, so $y \notin A_2$ and thus $y \in \partial_e(A_2|Q) \cap \mathcal{O}$.

Now let $y \in \partial_e(A_2|Q) \cap \mathcal{O}$, so that $y \sim z$ for some $z \in A_2$. If $z \in \Lambda_2$, then $z \in C_i$ for some $i$. Hewever, because $y \notin A_1$ and $y \sim z$, we have $y \in C_i$, a contradiction, so $z \in A_1$. If $z \in \Lambda_1$, then $z \in \mathcal{C}(x)$ for some $x \in \partial_e(A|Q) \cap \mathcal{O} - \mathcal{C}^\vee$. Again we have $y \in \mathcal{C}(x)$, a contradiction, so we must have $z \in A$ and, therefore, $y \in \partial_e(A|Q)$. If $y \notin \mathcal{C}^\vee$, then because $y \in \mathcal{O}$, $\mathcal{C}(y)$ is included in $\Lambda_1$ and $y$ is in $A_1$. Hence we have $y \in \mathcal{C}^\vee$. This completes the proof of (2.5).

Let $\Gamma_1, \ldots, \Gamma_l$ be the connected components of $Q - A_2$, arranged in decreasing order of size.

CASE 1. Suppose that $|\Gamma_1| > \frac{1}{2}|Q|$. Let $A_3 = A_2 \cup (\bigcup_{i=2}^l \Gamma_j)$. By Lemma 2.6, $\partial_e(A_3|Q) = \partial_i(\Gamma_1|Q)$ is $*$-connected. We also have $\partial_e(A_3|Q) \subseteq \partial_e(A_2|Q)$. Note that since $n < |A| \leq |A_3| \leq \frac{1}{2}|Q|$, by the isoperimetric inequality,

$$|\partial_e(A_3|Q)| \geq c_2|A_3|^{1-1/d} \geq c_2 n^{1-1/d} > n^\beta.$$

Since $F(Q,\varepsilon)$ holds, writing $\mathcal{O}' = \bigcap_{\sigma \in \Sigma} \mathcal{O}_\sigma$,

$$|\partial_e(A_3|Q) \cap (\mathcal{O}')^c| \leq \sum_\sigma |\partial_e(A_3|Q) \cap (\mathcal{O}_\sigma)^c| \leq (2d+1)\varepsilon|\partial_e(A_3|Q)|.$$

We deduce

$$|\partial_e(A_3|Q) \cap \mathcal{O}'| \geq (1 - (2d+1)\varepsilon)|\partial_e(A_3|Q)| \geq \tfrac{1}{2}|\partial_e(A_3|Q)|.$$

If $y \in \partial_e(A_3|Q) \cap \mathcal{O}'$, then $y \in \partial_e(A_2|Q) \cap \mathcal{O}$ and therefore, by (2.5), $y \in \partial_e(A|Q) \cap \mathcal{C}^\vee$. So there exists $x \in A$ with $y \sim x$. Whereas $y \in \mathcal{O}'$, we have $x \in \mathcal{O}$ and thus $x \in \mathcal{C}^\vee$. Hence $\{x,y\} \in \partial_E(A \cap \mathcal{C}^\vee, (Q-A) \cap \mathcal{C}^\vee)$. Thus

$$\partial_E(A \cap \mathcal{C}^\vee, (Q-A) \cap \mathcal{C}^\vee) \geq |\partial_e(A_3|Q) \cap \mathcal{O}'|$$
$$\geq \tfrac{1}{2}|\partial_e(A_3|Q)|$$
$$\geq \tfrac{1}{2}c_2|A_3|^{1-1/d} \geq \tfrac{1}{2}c_2 n^{-1}|A_3| \geq \tfrac{1}{2}c_2 n^{-1}|A|,$$

proving (2.4) in this case.

CASE 2. Suppose $|\Gamma_1| \leq \frac{1}{2}|Q|$. Note that

$$|\partial_E(A \cap \mathcal{C}^\vee, (Q-A) \cap \mathcal{C}^\vee)| = \sum_j |\partial_E(\Gamma_j \cap \mathcal{C}^\vee, (Q - \Gamma_j) \cap \mathcal{C}^\vee)|.$$

We have

(2.6) $$|\partial_i(\Gamma_j|Q) \cap \mathcal{O}'| \leq |\partial_E(\Gamma_j \cap \mathcal{C}^\vee, (Q - \Gamma_j) \cap \mathcal{C}^\vee)|.$$

For if $y \in \partial_i(\Gamma_j|Q) \cap \mathcal{O}'$, then $y \in \Gamma_j$ and there exists $x \in A_2$ with $x \sim y$. So $y \in \partial_e(A_2|Q) \cap \mathcal{O}$ and hence $y \in \partial_e(A|Q) \cap \mathcal{C}^\vee$ by (2.5). Whereas $y \in \mathcal{O}'$, then $x \in \mathcal{O}$, so that $x \in \mathcal{C}^\vee$ and hence $x \in A$ since $x$ cannot be in $\Lambda_1 \cup \Lambda_2$. Therefore $\{y,x\} \in \partial_E(\Gamma_j \cap \mathcal{C}^\vee, (Q - \Gamma_j) \cap \mathcal{C}^\vee)$.

Next, each $\Gamma_j$ intersects $\mathcal{C}^\vee$ by the construction of $A_2$. If $|\partial_i(\Gamma_j|Q)| \geq n^\beta$, then

$$|\partial_i(\Gamma_j|Q) \cap \mathcal{O}'| \geq \tfrac{1}{2}|\partial_i(\Gamma_j|Q)| \geq c|\Gamma_j|^{1-1/d} \geq c'n^{-1}|\Gamma_j|.$$

If $|\partial_i(\Gamma_j|Q)| < n^\beta$, then $|\Gamma_j| \leq cn^{\beta d/(d-1)} < n$ and

$$|\partial_E(\Gamma_j \cap \mathcal{C}^\vee, (Q - \Gamma_j) \cap \mathcal{C}^\vee)| \geq 1 \geq n^{-1}|\Gamma_j|.$$

Combining these estimates we obtain

$$|\partial_E(A \cap \mathcal{C}^\vee, (Q - A) \cap \mathcal{C}^\vee)| = \sum_j |\partial_E(\Gamma_j \cap \mathcal{C}^\vee, (Q - \Gamma_j) \cap \mathcal{C}^\vee)|$$

$$\geq cn^{-1} \sum_j |\Gamma_j| \geq cn^{-1}(|Q - A_2|).$$

Whereas $|Q - A_2| \geq \frac{1}{8}n^d \geq c|A|$, this proves (2.4). $\square$

We now follow the renormalization argument of Mathieu and Remy (2004), which uses techniques introduced by Antal and Pisztora (1996). We introduce a second percolation process, which is bond percolation on $(\mathbb{Z}^d, \mathbb{E}_d)$. Let $p \in (0, 1)$: We set $\Omega_b = \{0, 1\}^{\mathbb{E}_d}$, let $\eta_e$, $e \in \mathbb{E}_d$, be the coordinate maps and let $\mathbb{P}_p$ on $\Omega_b$ be the probability measure which makes $\eta_e$ i.i.d. Bernoulli r.v. with $\mathbb{P}_p(\eta_e = 1) = p$. We call the edges $e$ such that $\eta_e = 1$ *open edges* and we write $\mathcal{O}_E = \{e : \eta_e = 1\}$. An *open path* is any path $\gamma = \{x_0, \ldots, x_k\}$ such that each $\{x_{i-1}, x_i\} \in \mathcal{O}_E$. As for site percolation, we write $d_\omega(x, y)$ for the length of the shortest open path connecting $x$ and $y$, and set $d_\omega(x, y) = \infty$ if there is no such open path. For a set $A$ we write $d_A(x, y)$ for the length of the shortest open path contained in $A$ connecting $x$ and $y$. Set $\mathcal{C}(x) = \{y : d_\omega(x, y) < \infty\}$. Let $\theta(p) = \mathbb{P}_p(|\mathcal{C}(0)| = \infty)$ and $p_c = \inf\{p : \theta(p) > 0\}$. Then if $p > p_c$, there exists $\mathbb{P}_p$-a.s. a unique infinite cluster. We always assume that $p > p_c$. We define $\mathcal{C}_\infty$, $\mathcal{C}_Q(x)$ and $\mathcal{C}^\vee(Q)$ in the same way as for site percolation.

Let $n \geq 16$ be fixed and let $Q$ be a cube of side $n$. Recall the definition of $Q^+$ and crossing clusters from Definition 2.1. We define the event $R_0(Q)$ in a similar fashion to the event $R_{\mathbf{i}}^{(n)}$ in Antal and Pisztora (1996) and we set

$R_0(Q) = \{$there exists a unique crossing cluster $\mathcal{C}$ in $Q^+$ for $Q^+$, all open paths

contained in $Q^+$ of diameter greater than $\frac{1}{8}n$ are connected to $\mathcal{C}$

in $Q^+$ and $\mathcal{C}$ is crossing for each cube $Q' \subseteq Q$ with $s(Q') \geq n/8\}$,

$R(Q) = R_0(Q) \cap \{\mathcal{C}^\vee(Q) \text{ is crossing for } Q\} \cap \{\mathcal{C}^\vee(Q^+) \text{ is crossing for } Q^+\}.$

Note that if $\omega \in R(Q)$, then this forces $\mathcal{C}^\vee(Q) \subset \mathcal{C}^\vee(Q^+)$ and that $\mathcal{C}^\vee(Q^+)$ is the unique crossing cluster given by the event $R_0(Q)$.

Now let $k \geq 17$ and consider a tiling of $\mathbb{Z}^d$ by disjoint cubes

(2.7) $$T(x) = \{y \in \mathbb{Z}^d : x_i \leq y_i < x_i + k, 1 \leq i \leq d\}$$

with side $k - 1$. Let

(2.8) $$\varphi(x) = \mathbb{1}_{R(T(x))}.$$

LEMMA 2.8. (a) *Let $Q$ be a cube of side $k-1$ and let $p > p_c$. There exists $c_1 = c_1(p,d)$ such that*

(2.9) $$\mathbb{P}_p(R(Q)^c) \leq c\exp(-c_1 k).$$

(b) *The process $(\varphi_x, x \in \mathbb{Z}^d)$ dominates Bernoulli site percolation with parameter $q^*(k)$, where $q^*(k) \to 1$ as $k \to \infty$.*

PROOF. (a) The bound $\mathbb{P}_p(R_0(Q)^c) \leq c\exp(-c'n)$ follows from Theorem 3.1 of Pisztora (1996) ($d \geq 3$) and Theorem 5 of Penrose and Pisztora (1996) for $d = 2$. The estimate

$$\mathbb{P}_p(\mathcal{C}^\vee(Q) \text{ is not crossing for } Q) \leq c\exp(-c'n^{d-1})$$

follows from Theorem 1.2 of Pisztora (1996) ($d \geq 3$) and Theorem 1 of Couronné and Messikh (2003) ($d = 2$).

(b) By (a) we have $\mathbb{P}_p(\varphi(x) = 1) \to 1$ as $k \to \infty$. The r.v. $\varphi(x)$ and $\varphi(y)$ are independent if $|x-y|_\infty \geq 3$, so the result is immediate from Theorem 0.0 of Liggett, Schonmann and Stacey (1997). [We remark that using this theorem means that the events $S_{\mathbf{i}}^{(N)}$ defined in (2.5) of Antal and Pisztora (1996) are no longer needed.] □

Recall the definition of $q_1(\varepsilon, d)$ from Lemma 2.5 and choose $k_0 = k_0(p,d)$ large enough so that $q^*(k) \geq q_1(\varepsilon, d)$ for all $k \geq k_0$. We fix $k = k_0$ and refer to the process $\varphi$ as the *macroscopic percolation* process. We write $\widetilde{\mathcal{O}}, \widetilde{\mathcal{C}}(x)$, to denote the open sites, open clusters for the macroscopic process and $\widetilde{F}$ for associated events.

Let $\widetilde{Q}$ be a macroscopic cube of side $m$ and associate with $\widetilde{Q}$ the microscopic cube $Q = \bigcup\{T^+(x), x \in \widetilde{Q}\}$. Any (microscopic) cube which can be obtained in this way is called a *special cube*. Define $T'(x)$ to be $T(x)$ if $x$ is in the interior of $\widetilde{Q}$ [so that $T^+(x)$ is in the interior of $Q$]. Otherwise, if $x \in \partial_i(\widetilde{Q}|\mathbb{Z}^d)$, let $T'(x)$ be $T(x)$ together with all points in $T^+(x)$ which are closer to $T(x)$ than any $T(y)$, $y \neq x$, $y \in \widetilde{Q}$. Thus $T'(x)$ is a $d$-dimensional rectangle and $Q$ is the disjoint union of the $T'(x)$, $x \in \widetilde{Q}$. Note that each side of $T'(x)$ ($x \in \widetilde{Q}$) is less than $\frac{3}{2}(k_0 - 1)$.

Fix $\varepsilon_0 = (4d+2)^{-1}$.

LEMMA 2.9. *Let $Q$ be a special cube with $s(Q) = n$ and let $\widetilde{Q}$ be the associated macroscopic cube. Suppose that $m = s(\widetilde{Q}) \geq m_0(k_0)$ and that the event $\widetilde{K}(\widetilde{Q}, \frac{7}{8}) \cap \widetilde{F}(\widetilde{Q}, \varepsilon_0)$ holds for $(\varphi_x, x \in \mathbb{Z}^d)$. Then the cluster $\mathcal{C}^\vee(Q)$ satisfies*

$$|\mathcal{C}^\vee(Q)| \geq c_1 n^d, \qquad \text{diam}(\mathcal{C}^\vee(Q)) = n.$$

PROOF. Write $\mathcal{C}_x = \mathcal{C}^\vee(T^+(x))$. If $x, y \in \widetilde{\mathcal{O}} \cap \widetilde{Q}$ and $x \sim y$, then the events $R(T(x))$ and $R(T(y))$ force the clusters $\mathcal{C}_x$ and $\mathcal{C}_y$ to be connected. Thus there exists a $Q$-cluster $\mathcal{C}'$ with $\bigcup \{\mathcal{C}_x, x \in \widetilde{\mathcal{C}}^\vee(\widetilde{Q})\} \subseteq \mathcal{C}'$. It follows immediately that

$$|\mathcal{C}'| \geq |\widetilde{\mathcal{C}}^\vee(\widetilde{Q})| \geq \tfrac{7}{8} m^d \geq c_1 n^d.$$

Also, since $\widetilde{\mathcal{C}}^\vee$ is crossing for $\widetilde{Q}$, we deduce that $\mathcal{C}'$ is crossing for $Q$ and $\operatorname{diam}(\mathcal{C}') = n$.

It remains to prove that $\mathcal{C}' = \mathcal{C}^\vee(Q)$. Suppose not. Choose $m_0$ so that $(\tfrac{3}{2} k_0)^d < \tfrac{7}{8} m^d$ and therefore the cluster $\mathcal{C}^\vee(Q)$ is not contained in any one cube $T^+(x)$. Let $x \in \widetilde{\mathcal{C}}^\vee(\widetilde{Q})$. Then $\mathcal{C}^\vee(Q) \cap \mathcal{C}_x = \varnothing$ and so $\mathcal{C}^\vee(Q) \cap T^+(x)$ consists of clusters which have diameter less than $k_0/8$. Since $\mathcal{C}^\vee(Q)$ contains points outside $T^+(x)$, it follows that $\mathcal{C}^\vee(Q) \cap T(x) = \varnothing$.

So, if $\Gamma_1, \ldots, \Gamma_k$ are the connected components of $\widetilde{Q} - \widetilde{\mathcal{C}}^\vee(\widetilde{Q})$, we deduce that for some $j$, $\mathcal{C}^\vee(Q) \subseteq \bigcup \{T'(x), x \in \Gamma_j\}$. By Proposition 2.7(a) we deduce that $|\Gamma_j| \leq c_2 m^{d\beta/(d-1)}$ and so if $m_0$ is large enough,

$$|\mathcal{C}^\vee| \leq c_2 (\tfrac{3}{2} k_0)^d m^{d\beta/(d-1)} < \tfrac{7}{8} m^d,$$

giving a contradiction. Thus $\mathcal{C}' = \mathcal{C}^\vee(Q)$. □

Let $Q$ be a special cube with $s(Q) = n$ and let $\widetilde{Q}$ be the associated macroscopic cube. For sets $A_1, A_2 \subseteq Q$, let

$$\partial_E(A_1, A_2 | \mathcal{O}_E) = \{\{x, y\} : \{x, y\} \in \mathcal{O}_E, x \in A_1, y \in A_2\}.$$

Let $A \subseteq \mathcal{C}^\vee = \mathcal{C}^\vee(Q)$ and let $\Gamma = \mathcal{C}^\vee - A$ both be connected. Set

$$\widetilde{A} = \{x \in \widetilde{Q} : A \cap T'(x) \neq \varnothing\},$$
$$\widetilde{\Gamma} = \{x \in \widetilde{Q} : \Gamma \cap T'(x) \neq \varnothing\}.$$

Note that

(2.10) $$|A| \geq |\widetilde{A}| \geq (\tfrac{3}{2} k)^{-d} |A|.$$

Whereas $A$ and $\Gamma$ are connected, it is clear that $\widetilde{A}$ and $\widetilde{\Gamma}$ are connected. Let $\widetilde{\mathcal{C}}^\vee = \widetilde{\mathcal{C}}^\vee(\widetilde{Q})$ be the largest open (macroscopic) cluster in $\widetilde{Q}$.

LEMMA 2.10. *Let $Q$ be a special cube with $s(Q) = n$ and let $\widetilde{Q}$ be the associated macroscopic cube. Suppose that $m = s(\widetilde{Q}) \geq m_0(k_0)$ and that the event $\widetilde{K}(\widetilde{Q}, \tfrac{7}{8}) \cap \widetilde{F}(\widetilde{Q}, \varepsilon_0)$ holds for $(\varphi_x, x \in \mathbb{Z}^d)$. Let $A \subseteq \mathcal{C}^\vee$ and $\Gamma = \mathcal{C}^\vee - A$ be connected.*

(a) *Let $x \sim y$, $x \in \widetilde{A}$ and $y \in \widetilde{Q} - \widetilde{A}$ with $x, y \in \widetilde{\mathcal{O}}$. Then the set $T^+(x) \cap Q$ contains at least one edge in $\partial_E(A, \Gamma | \mathcal{O}_E)$.*

(b) *Suppose $x \in \widetilde{A} \cap \widetilde{\Gamma} \cap \widetilde{\mathcal{O}}$. Then the set $T^+(x) \cap Q$ contains at least one edge in $\partial_E(A, \Gamma | \mathcal{O}_E)$.*

PROOF. As in Lemma 2.9 we have that $\mathcal{C}^\vee$ is not contained in any one $T^+(x)$. Let $\mathcal{C}_x = \mathcal{C}^\vee(T^+(x))$ be as in Lemma 2.9.

(a) Since $x \in \widetilde{A}$ there exists $x' \in T(x) \cap A$, so since $\mathcal{C}^\vee$ is connected there exists an open path $\gamma \subseteq Q$ from $x'$ to $Q \cap T^+(x)^c$. This path must have diameter greater than $k_0/3$, so it is connected within $T^+(x)$ to $\mathcal{C}_x$. Hence $x'$ is connected within $T^+(x)$ to $\mathcal{C}_y$. Choose $y' \in \mathcal{C}_y \cap T(y)$; then $y' \notin A$ but $y \in \mathcal{C}^\vee$. There exists an open path $\gamma'$ from $x'$ to $y'$ that must contain at least one edge in $\partial_E(A, \Gamma | \mathcal{O}_E)$.

(b) Let $x' \in A \cap T'(x)$ and $y' \in \Gamma \cap T'(x)$. Since $x \in \widetilde{\mathcal{O}}$, both $x'$ and $y'$ are in $\mathcal{C}_x$ and there therefore exists an open path in $T^+(x)$ between $x'$ and $y'$. As in (a) this path must contain at least one edge in $\partial_E(A, \Gamma | \mathcal{O}_E)$. □

PROPOSITION 2.11. *Let $Q$ and $\widetilde{Q}$ be as above with $m = s(\widetilde{Q}) \geq m_0$. Suppose that the event $\widetilde{K}(\widetilde{Q}, \frac{7}{8}) \cap \widetilde{F}(\widetilde{Q}, \varepsilon_0)$ holds for $(\varphi_x, x \in \mathbb{Z}^d)$. Let $A$ be a connected open subset of $\mathcal{C}^\vee(Q)$ with $|A| \leq \frac{1}{2}|\mathcal{C}^\vee(Q)|$ and such that $\Gamma = \mathcal{C}^\vee(Q) - A$ is also connected. Then*

$$(2.11) \qquad |\partial_E(A, \Gamma | \mathcal{O}_E)| \geq c_1 n^{-1} |A|.$$

PROOF. We write $\mathcal{C}^\vee = \mathcal{C}^\vee(Q)$, $\widetilde{\mathcal{C}}^\vee = \widetilde{\mathcal{C}}^\vee(\widetilde{Q})$. Since $A \neq \mathcal{C}^\vee(Q)$ and $\mathcal{C}^\vee$ is connected, we have $|\partial_E(A, \Gamma | \mathcal{O}_E)| \geq 1$, so (2.11) is immediate if $|A| \leq c_1^{-1} n$. Also, using Lemma 2.10(a) we have

$$(2.12) \qquad |\partial_E(A, \Gamma | \mathcal{O}_E)| \geq c_2 |\partial_E(\widetilde{A} \cap \widetilde{\mathcal{O}}, (\widetilde{Q} - \widetilde{A}) \cap \widetilde{\mathcal{O}})|.$$

We consider four cases.

CASE 1. $\widetilde{A} \cap \widetilde{\mathcal{C}}^\vee = \varnothing$. Let $\widetilde{\Lambda}$ be the connected component of $\widetilde{Q} - \widetilde{\mathcal{C}}^\vee$ which contains $\widetilde{A}$. By Proposition 2.7(a), $|\widetilde{\Lambda}| \leq c m^{d\beta/(d-1)}$, so that

$$|A| \leq (\tfrac{3}{2} k_0)^d |\widetilde{\Lambda}| \leq c k_0^d m^{d\beta/(d-1)} \leq n,$$

provided $m_0$ is large enough. Hence (2.11) holds for $A$.

CASE 2. $\widetilde{A} \cap \widetilde{\mathcal{C}}^\vee \neq \varnothing$ and $|\widetilde{A}| \leq \frac{3}{4}|\widetilde{Q}|$. We apply Proposition 2.7(b) to see that

$$|\partial_E(\widetilde{A} \cap \widetilde{\mathcal{C}}^\vee, (\widetilde{Q} - \widetilde{A}) \cap \widetilde{\mathcal{C}}^\vee)| \geq c m^{-1} |\widetilde{A}| \geq c n^{-1} |A|,$$

and combining this with (2.12) proves (2.11).

CASE 3. $\widetilde{A} \cap \widetilde{\mathcal{C}}^{\vee} \neq \varnothing$, $|\widetilde{A}| \geq \frac{3}{4}|\widetilde{Q}|$ and $|\widetilde{\Gamma}| \leq \frac{3}{4}|\widetilde{Q}|$. Using Lemma 2.9 we have $|\widetilde{\Gamma}| \geq (\frac{3}{2}k_0)^{-d}|\Gamma| \geq \frac{1}{2}(\frac{3}{2}k_0)^{-d}|\mathcal{C}^{\vee}| \geq cn^d$. So by Proposition 2.7(a), $\widetilde{\Gamma} \cap \widetilde{\mathcal{C}}^{\vee} \neq \varnothing$ and we can therefore apply Proposition 2.7(b) to $\widetilde{\Gamma}$ to deduce

$$|\partial_E(\widetilde{\Gamma} \cap \widetilde{\mathcal{C}}^{\vee}, (\widetilde{Q} - \widetilde{\Gamma}) \cap \widetilde{\mathcal{C}}^{\vee})| \geq cm^{-1}|\widetilde{\Gamma}| \geq cm^{-1}n^d \geq cn^{-1}|A|.$$

Using (2.12) (with $A$ and $\Gamma$ interchanged) we deduce that

$$|\partial_E(\Gamma, A|\mathcal{O}_E)| \geq c_2|\partial_E(\widetilde{\Gamma} \cap \widetilde{\mathcal{O}}, (\widetilde{Q} - \widetilde{\Gamma}) \cap \widetilde{\mathcal{O}})|$$
$$\geq c_2|\partial_E(\widetilde{\Gamma} \cap \widetilde{\mathcal{C}}^{\vee}, (\widetilde{Q} - \widetilde{\Gamma}) \cap \widetilde{\mathcal{C}}^{\vee})| \geq n^{-1}|A|.$$

CASE 4. $\widetilde{A} \cap \widetilde{\mathcal{C}}^{\vee} \neq \varnothing$, $|\widetilde{A}| \geq \frac{3}{4}|\widetilde{Q}|$ and $|\widetilde{\Gamma}| \geq \frac{3}{4}|\widetilde{Q}|$. Since $\widetilde{A}$ and $\widetilde{\Gamma}$ are both connected and $\widetilde{F}(\widetilde{Q}, \varepsilon_0)$ holds, we have

$$|\widetilde{A} \cap \widetilde{\mathcal{O}}| \geq (1 - \varepsilon_0)|\widetilde{A}| > \tfrac{2}{3}|\widetilde{Q}|,$$
$$|\widetilde{\Gamma} \cap \widetilde{\mathcal{O}}| \geq (1 - \varepsilon_0)|\widetilde{\Gamma}| > \tfrac{2}{3}|\widetilde{Q}|.$$

Hence $|\widetilde{A} \cap \widetilde{\Gamma} \cap \widetilde{\mathcal{O}}| \geq \frac{1}{3}|\widetilde{Q}|$. So by Lemma 2.10(b),

$$|\partial_E(A, \Gamma|\mathcal{O}_E)| \geq c|\widetilde{A} \cap \widetilde{\Gamma} \cap \widetilde{\mathcal{O}}| \geq c'm^d \geq c''|A|,$$

which implies (2.11). $\square$

As in Section 1 we define

$$\nu_{xy} = \nu_{xy}(\omega) = \begin{cases} 1, & \text{if } \{x,y\} \text{ is an open edge,} \\ 0, & \text{otherwise,} \end{cases}$$

$\mu(x) = \mu(x)(\omega) = \sum_y \nu_{xy}$, $x \in \mathbb{Z}^d$, and we extend $\nu$ and $\mu$ to measures. If $f : A \to \mathbb{R}$, we write

$$\bar{f}_A = \mu(A)^{-1} \int_A f \, d\mu.$$

Let $Q$ be a special cube and let $\widetilde{Q}$ be the associated macroscopic cube. Define

$$H_0(Q) = \widetilde{K}(\widetilde{Q}, \tfrac{7}{8}) \cap \widetilde{F}(\widetilde{Q}, \varepsilon_0)$$

and extend the definition of $H_0(Q)$ to all cubes $Q$ by taking $H_0(Q) = H_0(Q')$, where $Q'$ is the largest special cube contained in $Q$. [We set $H_0(Q) = \Omega$ if there is no such special cube.]

Recall the definition of $\beta$ from (2.2).

PROPOSITION 2.12. *Consider bond percolation $\eta_e$ on $(\mathbb{Z}^d, \mathbb{E}_d)$ with $p > p_c$. Let $Q$ be a special cube of side $n$.*

(a) *There exists*

(2.13) $$\mathbb{P}_p(H_0(Q)^c) \leq c_1 \exp(-c_2 n^\beta).$$

(b) *If $\omega \in H_0(Q)$ and $f: \mathcal{C}^\vee(Q)(\omega) \to \mathbb{R}$, then*

$$\int_{\mathcal{C}^\vee(Q)} (f(y) - \bar{f}_{\mathcal{C}^\vee(Q)})^2 \, d\mu \leq c_3 n^2 \int_{E(\mathcal{C}^\vee(Q))} |\nabla f|^2 \, d\nu.$$

PROOF. (a) To prove (2.13) note that

(2.14) $$\mathbb{P}_p(H_0(Q)^c) \leq \mathbb{P}_p(\widetilde{K}(\widetilde{Q}, \tfrac{7}{8})^c) + \mathbb{P}_p(\widetilde{F}(\widetilde{Q}, \varepsilon_0)^c).$$

As the events $K(\cdot, \tfrac{7}{8})$ and $F(\cdot, \cdot)$ are increasing, the two probabilities on the right-hand side of (2.14) are bounded by the probabilities of these events with respect to a Bernoulli site percolation process with probability $q^* = q^*(k_0)$. Using Lemmas 2.2 and 2.5,

$$\mathbb{P}_p(\widetilde{K}(\widetilde{Q}, \tfrac{7}{8})^c) + \mathbb{P}_p(\widetilde{F}(\widetilde{Q}, \varepsilon_0)^c) \leq \mathbb{Q}_{q^*}(K(\widetilde{Q}, \tfrac{7}{8})^c) + \mathbb{Q}_{q^*}(F(\widetilde{Q}, \varepsilon_0)^c)$$
$$\leq c \exp(-cm^{d-1}) + c' \exp(-c'm^\beta)$$
$$\leq c \exp(-c'n^\beta).$$

(b) This is immediate from Propositions 1.4 and 2.11. □

REMARK. See Section 3.2 of Mathieu and Remy (2004) for similar bounds in the case when $p$ is close to 1.

Recall from Definition 2.1 that $Q \subseteq Q^\oplus \subseteq Q^+$. Let $\alpha \in (0, \tfrac{1}{2})$, let $Q$ be a cube with $s(Q) = n$ and set

(2.15) $$H(Q, \alpha) = R(Q) \cap \{R(Q') \cap H_0(Q') \text{ occurs for every cube } Q' \\ \text{with } (Q')^+ \subseteq Q^+, \, Q' \cap Q^\oplus \neq \varnothing \\ \text{and } n^\alpha \leq s(Q') \leq n\}.$$

LEMMA 2.13. (a) *For $p > p_c$,*

$$\mathbb{P}_p(H(Q, \alpha)^c) \leq c \exp(-c_1 n^{\alpha\beta}).$$

(b) *If $\omega \in H(Q, \alpha)$ and $Q_0 \subset Q$ satisfies the condition (2.15), then $\mathcal{C}^\vee(Q_0) \subset \mathcal{C}^\vee(Q^+)$.*

PROOF. (a) By (2.9) and (2.13),

(2.16) $$\mathbb{P}_p(H(Q, \alpha)^c) \leq \sum_{r=n^\alpha}^{n} cn^d \exp(-c_1 r^\beta) \leq c' \exp(-c_2 n^{\alpha\beta}),$$

proving (a).

(b) Define a sequence of cubes $Q_i$, $0 \le i \le k$, by $Q_{i+1} = Q_i^*$ and where we stop at the last cube $Q_k$ with $Q_k^+ \subset Q^+$. The events $R(Q_i)$ then force $\mathcal{C}^\vee(Q_i) \subset \mathcal{C}^\vee(Q_{i+1})$, so that $\mathcal{C}^\vee(Q_0) \subset \mathcal{C}^\vee(Q_k^+)$. Since $\operatorname{diam} \mathcal{C}^\vee(Q_k^+) = \operatorname{diam}(Q_k^+) > n/8$, the event $R(Q)$ implies that $\mathcal{C}^\vee(Q_k^+) \subset \mathcal{C}^\vee(Q^+)$. □

Proposition 2.12 and Lemma 2.13 complete our results on the Poincaré inequality in sets of the form $\mathcal{C}^\vee(Q)$. However, to be able to obtain bounds (particularly lower bounds) on transition densities on $\mathcal{C}_\infty$, we need to relate $|x - y|_1$ to the shortest path (or chemical) metric $d_\omega$ on $\mathcal{C}_\infty$. This was done in Theorem 1.1 of Antal and Pisztora (1996), but we need some minor extensions of their results and it is desired to make this paper as self-contained as possible, so we repeat some of their constructions.

Let $k \ge 17$ and recall the site process $\varphi(x) = \mathbb{1}_{R(T(x))}$ introduced in (2.8). Set $\varphi'(x) = 1 - \varphi(x)$; since, by Lemma 2.8(b), $\varphi$ dominates Bernoulli site percolation with parameter $q^*(k)$, $\varphi$ is dominated by Bernoulli site percolation with parameter $q'(k) = 1 - q^*(k)$. We consider the clusters of the process $\varphi'$ on the graph $(\mathbb{Z}^d, \mathbb{E}_d^*)$. Given a function $\varphi' : \mathbb{Z}^d \to \{0, 1\}$, we write $\mathcal{O}(\varphi') = \{x : \varphi'(x) = 1\}$. If $\varphi'(x) = 0$, we set $\mathcal{C}^*(\varphi', x) = \varnothing$, and if $\varphi'(x) = 1$, we let $\mathcal{C}^*(\varphi', x)$ be the $*$-connected component of $\mathcal{O}(\varphi')$ that contains $x$. Let

$$(2.17) \qquad \mathcal{D}(\varphi', x) = \begin{cases} \partial_e(\mathcal{C}^*(\varphi', x) | \mathbb{Z}^d), & \text{if } \varphi'(x) = 1, \\ \{x\}, & \text{if } \varphi'(x) = 0. \end{cases}$$

For $x, y \in \mathbb{Z}^d$ let $\gamma(x, y)$ be a shortest path [in $(\mathbb{Z}^d, \mathbb{E}^d)$] that connects $x$ and $y$. Note that if $x$ and $y$ are contained in a cube $Q$, then $\gamma(x, y) \subseteq Q$. Set

$$(2.18) \qquad W(\varphi', x, y) = \bigcup_{z \in \gamma(x, y)} \mathcal{D}(\varphi', z).$$

The importance of $W$ comes from the following result, which is Proposition 3.1 of Antal and Pisztora (1996).

PROPOSITION 2.14. *Let $k \ge 17$, $x, y \in \mathbb{Z}^d$ and $\tilde{x}, \tilde{y} \in \mathbb{Z}^d$ be such that $x \in T(\tilde{x})$ and $y \in T(\tilde{y})$. Suppose that [for the bond percolation process $\eta_e$ on $(\Omega_b, \mathbb{P}_p)$] $d_\omega(x, y) < \infty$. Then there exists an open path $\gamma'(x, y)$ that connects $x$ and $y$ contained in*

$$W'(x, y) = \bigcup_{\tilde{z} \in W(\varphi', \tilde{x}, \tilde{y})} T^+(\tilde{z}).$$

*In particular, $d_\omega(x, y) \le |W'(x, y)| \le (3k)^d |W(\varphi', \tilde{x}, \tilde{y})|$.*

PROPOSITION 2.15. *Let $p > p_c$. There exists $k_1 = k_1(p, d)$ and $C_{\text{AP}}$ such that if $k \ge k_1$, then the following hold.*

(a) *If $\tilde{x}$, $\tilde{y} \in \mathbb{Z}^d$, $x \in T(\tilde{x})$ and $y \in T(\tilde{y})$, then*

$$\mathbb{P}_p(|W(\varphi', x, y)| \geq C_{\text{AP}}|\tilde{x} - \tilde{y}|_1) \leq c_2 \exp(-c_3|\tilde{x} - \tilde{y}|_1).$$

(b) *If $\tilde{x}$ and $\tilde{y} \in \mathbb{Z}^d$, $\lambda \geq 0$, then*

$$\mathbb{P}_p\bigg(\max_{\tilde{z} \in \gamma(\tilde{x}, \tilde{y})} \operatorname{diam}(\mathcal{D}(\varphi', \tilde{z})) > \lambda\bigg) \leq c_4 |\tilde{x} - \tilde{y}|_1 \exp(-c_5 \lambda).$$

PROOF. (a) This is proved on page 1047 of Antal and Pisztora (1996).

(b) We choose $k_1$ large enough so that $q'(k) < q_c$ for all $k \geq k_1$. Since $\varphi'$ is dominated by Bernoulli site percolation with parameter $q'$, we have

$$\mathbb{P}_p(\operatorname{diam}(\mathcal{D}(\varphi', \tilde{z})) > \lambda) \leq \mathbb{Q}_{q'}(2 + \operatorname{diam}(\mathcal{C}^*(0)) > \lambda) \leq c \exp(-c' \lambda).$$

[For the second estimate above, see Theorem 5.4 of Grimmett (1999).] The bound in (b) is now immediate. $\square$

Now fix $k_1$ as in Proposition 2.15 and set

(2.19) $$C_H = dC_{\text{AP}}(3k_1)^d.$$

Let $Q$ be a cube with side $n$. For $x, y \in Q$ let

$$E(Q, x, y) = \{x, y \in \mathcal{C}^\vee(Q^+) : d_{\mathcal{C}^\vee(Q^+)}(x, y) > C_H |x - y|_\infty\}.$$

LEMMA 2.16. *Let $p > p_c$, let $Q$ be a cube with side $n$ and let $x, y \in Q$. Then*

$$\mathbb{P}_p(E(Q, x, y)) \leq c \exp(-c_5 |x - y|_\infty).$$

PROOF. Let $x \in T(\tilde{x})$ and $y \in T(\tilde{y})$. Suppose $x, y \in \mathcal{C}^\vee(Q^+)$. Then by Proposition 2.14 there exists a open path $\gamma'$ that connects $x$ and $y$ contained in $W'(\varphi', x, y)$. So if $E(Q, x, y)$ holds, then either $W'(\varphi', x, y)$ is not contained in $Q^+$ or

$$C_H |x - y|_\infty < d_{\mathcal{C}^\vee(Q^+)}(x, y) \leq (3k_1)^d |W(\varphi', \tilde{x}, \tilde{y})|.$$

Thus

$$\mathbb{P}_p(E(Q, x, y)) \leq \mathbb{P}_p(|W(\varphi', \tilde{x}, \tilde{y})| \geq C_{\text{AP}}|\tilde{x} - \tilde{y}|_1)$$
$$+ \mathbb{P}_p\bigg(\max_{\tilde{z} \in \gamma(\tilde{x}, \tilde{y})} \operatorname{diam}(\mathcal{D}(\varphi', \tilde{z})) > \frac{n}{8k_1}\bigg)$$
$$\leq c \exp(-c|\tilde{x} - \tilde{y}|_1) + c|\tilde{x} - \tilde{y}|_1 \exp\bigg(-\frac{cn}{k_1}\bigg)$$
$$\leq c \exp(-c|x - y|_\infty). \qquad \square$$

Let $Q$ be a cube with side $n$. Set

(2.20) $$D_0(Q) = R(Q) \cap \{d_{\mathcal{C}^\vee(Q^+)}(x,y) \leq C_H |x-y|_\infty \\ \text{if } x, y \in \mathcal{C}^\vee(Q^+) \cap Q, |x-y|_\infty \geq n/12\}$$

and

(2.21) $$D(Q, \alpha) = \{D_0(Q') \text{ occurs for every cube } Q' \text{ with } (Q')^+ \subseteq Q^+, \\ Q' \cap Q^\oplus \neq \emptyset \text{ and } n^\alpha \leq s(Q') \leq n\}.$$

Let also

$$B_\omega(y, r) = \{x : d_\omega(x, y) < r\}.$$

Since $\mathcal{C}_\infty$ is embedded in $\mathbb{Z}^d$ we have $\mu(B_\omega(y, r)) \leq C_0 r^d$ for some $C_0 = C_0(d)$.

PROPOSITION 2.17. *Let $p > p_c$.*

(a) *There exists $\mathbb{P}_p(D_0(Q)^c) \leq c_1 \exp(-c_2 n)$.*
(b) *There exists $\mathbb{P}_p(D(Q,\alpha)^c) \leq c_3 \exp(-c_4 n^\alpha)$.*
(c) *Let $\omega \in D_0(Q)$ and $x, y \in Q \cap \mathcal{C}^\vee(Q^+)$. Then $d_\omega(x, y) \leq C_H n$.*
(d) *Let $\omega \in D(Q, \alpha)$ and $x, y \in Q^\oplus \cap \mathcal{C}^\vee(Q^+)$. Then*

$$|x - y|_\infty \leq d_\omega(x, y) \leq C_H((1 + n^\alpha) \vee |x - y|_\infty).$$

(e) *Let $\omega \in D(Q, \alpha)$, let $x \in Q^\oplus$ and let $Q'$ satisfy the conditions in (2.21) with $s(Q') = r$. Then $Q' \cap \mathcal{C}^\vee(Q^+) \subset B_\omega(x, C_H r)$.*

PROOF. (a) By Lemma 2.16,

$$\mathbb{P}_p(D_0(Q)^c) \leq \mathbb{P}_p(R(Q)^c) + \sum_{x', y'} c \exp(-c_5 |x' - y'|_\infty),$$

where the sum is over $x', y' \in Q$ with $|x' - y'|_\infty \geq n/12$. Hence, using (2.9),

(2.22) $$\mathbb{P}_p(D_0(Q)^c) \leq c(n+1)^{2d} \exp(-c_5 n/8) \leq c_1 \exp(-c_2 n).$$

(b) This is immediate from (a), since

$$\mathbb{P}_p(D(Q, \alpha)^c) \leq \sum_{r=n^\alpha}^{n} (\tfrac{3}{2}n + 1)^d c_1 \exp(-c_2 r) \leq c \exp(-c' n^\alpha).$$

(c) This is immediate if $|x - y|_\infty \geq n/12$. If $|x - y|_\infty < n/8$, then choose $z' \in \mathcal{C}^\vee(Q^+) \cap Q$ with $n/12 \leq |z - x|_\infty \leq n/4$ and $n/12 \leq |z - y|_\infty \leq n/4$. Since $\mathcal{C}^\vee(Q^+)$ is crossing for $Q$, such a choice of $z$ is possible. Then $d_\omega(x, y) \leq d_\omega(x, z) + d_\omega(z, y) \leq 2C_H(n/4) \leq C_H n$.

(d) Since $D(Q, \alpha) \subseteq D_0(Q)$, this is immediate from (b) if $|x-y|_\infty \geq n/12$. Otherwise choose the smallest possible cube $Q'$ such that $s(Q') \geq n^\alpha \vee |x-y|_\infty$ and $x, y \in Q'$. We have $(Q')^+ \subset Q^+$. As in Lemma 2.13(b), we have $\mathcal{C}^\vee(Q') \subset \mathcal{C}^\vee(Q^+)$ and, by (c), $d_{\mathcal{C}^\vee(Q^+)}(x,y) \leq C_H s(Q') \leq C_H(|x-y|_\infty \vee (1+n^\alpha))$.

(e) Since $D_0(Q')$ occurs, this is immediate from (c). □

Recall from Definition 1.7 the definition of good and very good balls.

THEOREM 2.18. *Let $\alpha \in (0, \frac{1}{2})$, let $Q$ be a cube of side $n$, let $\omega \in H(Q, \alpha) \cap D(Q, \alpha)$ and let $C_H n^\alpha \leq r \leq n$. Write $Q(y, s) = \{z \in \mathbb{Z}^d : |z - y|_\infty \leq s\}$. Let $y \in \mathcal{C}^\vee(Q^+) \cap Q^\oplus$ with $Q(y, r + k_0)^+ \subseteq Q^+$.*

(a) *There exists $C_V = C_V(p, d)$ such that*

$$(2.23) \qquad C_V r^d \leq |B_\omega(y, r)| \leq C_0 r^d.$$

(b) *There exist constants $C_P(p, d)$ and $C_W(p, d)$ such that if $f : B_\omega(y, C_W r) \to \mathbb{R}$ and writing $\bar{f}_B = \bar{f}_{B_\omega(y,r)}$, then*

$$(2.24) \qquad \int_{B_\omega(y,r)} (f - \bar{f}_{B_\omega(y,r)})^2 \, d\mu \leq C_P r^2 \int_{E(B_\omega(y, C_W r))} |\nabla f|^2 \, d\nu.$$

(c) *If $(C_H n^\alpha)^{d+2} \leq R \leq n$ and $B_\omega(y, \frac{3}{2}R) \subseteq Q^\oplus$, then $B_\omega(y, R)$ is $(C_V, C_P, C_W)$-very good with $N_{B_\omega(y,R)} \leq C_H n^\alpha$.*

PROOF. Recall from Lemma 2.13(b) that $\omega \in H(Q, \alpha)$ implies that $\mathcal{C}^\vee(Q') \subseteq \mathcal{C}^\vee(Q^+)$ for every $Q'$ satisfying (2.15).

(a) Since $B_\omega(y, r) \subseteq Q(y, r)$, the upper bound in (2.23) is clear. Let $s = r/(2C_H)$, so that $2s \geq n^\alpha$. By Proposition 2.17(e), $\mathcal{C}^\vee(Q(y, s)) \subset B_\omega(y, C_H s)$, so that by Lemma 2.9,

$$|B_\omega(y, r)| \geq |\mathcal{C}^\vee(Q(y, s))| \geq c_1 s^d \geq c_2 r^d.$$

(b) Let $Q_1$ be the smallest special cube that contains $Q(y, r)$; we have $Q_1^+ \subseteq Q^+$. Let $r_1 = s(Q_1)$. By Proposition 2.17(e), $\mathcal{C}^\vee(Q_1) \subset B_\omega(y, C_H r_1)$ and so taking $C_W = 2C_H$, $\mathcal{C}^\vee(Q_1) \subseteq B_\omega(y, C_W r)$. By Proposition 2.12(b),

$$\int_{B_\omega(y,r)} (f - \bar{f}_B)^2 \, d\mu \leq \int_{B_\omega(y,r)} (f - \bar{f}_{\mathcal{C}^\vee})^2 \, d\mu$$

$$\leq \int_{\mathcal{C}^\vee(Q_1)} (f - \bar{f}_{\mathcal{C}^\vee})^2 \, d\mu$$

$$\leq c_3 r_1^2 \int_{E(\mathcal{C}^\vee(Q_1))} |\nabla f|^2 \, d\nu \leq c_4 r^2 \int_{E(B_\omega(y, C_W r))} |\nabla f|^2 \, d\nu.$$

(c) This is immediate from (a), (b) and the definition of very good balls. □

Using the estimates in Lemma 2.13(a), Proposition 2.17(b) and and Borel–Cantelli lemma, we obtain the following lemma.

LEMMA 2.19. *Let $p > p_c$. For each $x \in \mathbb{Z}^d$ there exists $M_x(\omega)$ with $\mathbb{P}_p(M_x \geq n) \leq c_1 \exp(-c_2 n^{\alpha\beta})$ such that whenever $n \geq M_x$, then $H(Q,\alpha) \cap D(Q,\alpha)$ holds for all cubes $Q$ of side $n$ with $x \in Q$.*

REMARKS. 1. An inequality of the form (2.24) is called a *weak Poincaré inequality*. In many situations (including this one) it is possible to derive a strong Poincaré inequality (i.e., with $C_W = 1$) from a family of weak ones; see Lemma 4.9.

2. Note that if $x \in Q$ and $s(Q) \geq M_x$, then $\mathcal{C}^\vee(Q) \subseteq \mathcal{C}_\infty$.

Theorem 2.18 and Lemma 2.19 are suitable for most of our needs, but they have the defect that the minimum size of ball inside a cube $Q$ of side $n$ for which the Poincaré inequality is certain to hold increases with $n$. Since (for a fixed $C_P$) the cluster $\mathcal{C}_\infty$ contains arbitrarily large balls in which the Poincaré inequality fails, we cannot do better than this as long as we require it for all balls of some size. However, we can improve Theorem 2.18 if we relax this condition, and in Section 5 we want to connect points by a chain of very good balls of some fixed size. To do this we need an additional percolation argument.

We consider again Bernoulli site percolation $\zeta_x$ on $(\mathbb{Z}^d, \mathbb{E}_d)$ with parameter $q$, where $q > q_c$ is close to 1. Let $Q$ be a cube of side $n$. For $x, y \in Q$, $\lambda > 1$, let

$$\begin{aligned}(2.25)\quad & S(\zeta, Q, \lambda, x, y) \\ & = K(Q, \tfrac{7}{8}) \cap \{\text{there exist } x', y' \in \mathcal{C}^\vee(Q^+) \\ & \qquad \text{with } |x - x'|_\infty \leq n^{1/9}, |y - y'|_\infty \leq n^{1/9}, \\ & \qquad \text{such that } d_{\mathcal{C}^\vee(Q^+)}(x', y') \leq \lambda |x - y|_\infty\}.\end{aligned}$$

Note that this event is increasing.

LEMMA 2.20. *Let $Q$ be a cube of side $n$ and let $x, y \in Q$. There exists $q_2 = q_2(d) \in (q_c, 1)$ and $\lambda_0 \geq 1$ so that if $q > q_2$, then*

$$(2.26) \qquad \mathbb{Q}_q(S(\zeta, Q, \lambda_0, x, y)^c) \leq c_1 \exp(-c_2 n^{1/9}).$$

PROOF. We follow the proof of Theorem 1.1 of Antal and Pisztora (1996) and consider the dual process to $\zeta$ given by $\zeta'_x = 1 - \zeta_x$. We view $\zeta'$ as site percolation (with parameter $q' = 1 - q$) on the lattice $(\mathbb{Z}^d, \mathbb{E}_d^*)$ and write $\mathcal{C}^*(\zeta', z)$ for the $*$-connected cluster of the process $\zeta'$ that contains $z$. Then by Theorem 5.4 of Grimmett (1999) we can choose $q_2$ large enough so that if $q \geq q_2$, then

$$\mathbb{Q}_q(|\mathcal{C}^*(\zeta', z)| \geq k) \leq \exp(-c_3 k), \qquad k \geq 1.$$

So, if

$$G = \{|\mathcal{C}^*(\zeta', x)| \leq n^{1/9} \text{ for all } x \in Q\},$$

then (using Lemma 2.2)

(2.27) $$\mathbb{Q}_q(K(Q, \tfrac{7}{8})^c \cup G^c) \leq \exp(-c_4 n^{1/9}).$$

If $|x - y|_\infty \leq n^{1/9}$ and $\omega \in G$, there exists $x' \in \mathcal{C}^\vee(Q^+)$ with $|x - x'|_\infty \leq n^{1/9}$. In this case we can take $y' = x'$.

So suppose $\omega \in G$ and $|x - y|_\infty > n^{1/9}$. Let $l = |x - y|_1$ and $\gamma = \{x = x_0, x_1, \ldots, x_l = y\}$ be a path in $(\mathbb{Z}^d, \mathbb{E}_d)$ of length $l$ that connects $x$ and $y$—note that $\gamma \subset Q$. Whereas each cluster $\mathcal{C}^*(\zeta', y)$, $y \in Q$, has diameter less than $n^{1/9}$, the path $\gamma$ must intersect $\mathcal{C}^\vee(Q^+)$. Let $V_x$ and $V_y$ be the first and last (resp.) points in $\gamma \cap \mathcal{C}^\vee(Q^+)$; we have $|x - V_x| \leq n^{1/9}$ and $|y - V_y| \leq n^{1/9}$.

We take $x' = V_x$, $y' = V_y$ and construct a path $\Gamma$ from $x'$ to $y'$ in $\mathcal{C}^\vee(Q^+)$. This path follows $\gamma$ whenever possible, and when it encounters a site $z$ with $\zeta_z = 0$ it "walks around" $\mathcal{C}^*(\zeta', z)$—this requires at most $3^d |\mathcal{C}^*(\zeta', z)|$ steps. Since $\omega \in G$ this path does not leave $Q^+$. Hence, recalling from (2.18) the definition of $W(\cdot)$,

(2.28) $$|\Gamma| \leq l + 3^d \sum_{i=0}^{l} |\mathcal{C}^*(\zeta', z)| \leq 3^d |W(\zeta', x', y')| \leq 3^d |W(\zeta', x, y)|.$$

By Proposition 2.15(a) we have

(2.29) $$\mathbb{Q}_q(|W(\zeta', x, y)| \geq c_5 |x - y|_1) \leq c_6 \exp(-c_7 |x - y|_1) \leq c_6 \exp(-c_8 n^{1/9}).$$

Taking $\lambda_0 = dc_5$ and combining the bounds (2.27)–(2.29) completes the proof. □

Let $m \geq k_0 \vee k_1$ and let $\{T^m(\tilde{x}), \tilde{x} \in \mathbb{Z}^d\}$ be the tiling of $\mathbb{Z}^d$ by disjoint cubes of side $m - 1$ given by (2.7). Let $\alpha_1 = 1/(4+d)$ and define [on the space $(\Omega, \mathbb{P}_p)$ carrying the bond percolation process $\eta_e$]

(2.30) $$\psi_{\tilde{x}}^{(m)} = \mathbb{1}_{H(T^m(\tilde{x}), \alpha_1) \cap D(T^m(\tilde{x}), \alpha_1)}, \qquad \tilde{x} \in \mathbb{Z}^d.$$

LEMMA 2.21. *There exists $C_E = C_E(d, p) \geq 1$ such that for any $k \geq C_E$ the process $\psi^{(m)}$, $\tilde{x} \in \mathbb{Z}^d$, under $\mathbb{P}_p$ dominates Bernoulli site percolation on $\mathbb{Z}^d$ with parameter $q_2$.*

PROOF. Note that $\psi_{\tilde{x}}^{(m)}$ and $\psi_{\tilde{y}}^{(m)}$ are independent if $|\tilde{x} - \tilde{y}|_\infty \geq 3$. Using this and the fact that $\mathbb{P}_p(H(Q, \alpha_1) \cap D(Q, \alpha_1)) \to 1$ as $p \uparrow 1$, this is an immediate consequence of Theorem 0.0 of Liggett, Schonmann and Stacey (1997). □

Let $\lambda_0$ be as in Lemma 2.20 and let $Q$ be a cube of side $n$. For $x_0, x_1 \in Q$ and $C_E \leq m \leq n^{1/9}$ set

(2.31)
$$\begin{aligned}L(Q, m, x_0, x_1) \\= \{\text{there exist } x_0', x_1' \in \mathcal{C}^\vee(Q^+) \text{ with } |x_j - x_j'| \leq n^{2/9}, j = 0, 1, \\k \text{ with } mk < 2\lambda_0 |x_0 - x_1|_\infty \text{ and a path } \{\tilde{y}_0, \ldots, \tilde{y}_k\} \text{ in } (\mathbb{Z}^d, \mathbb{E}_d) \\\text{such that } T^m(\tilde{y}_i) \subseteq Q^+, 0 \leq i \leq k, x_0' \in T^m(\tilde{y}_0), x_1' \in T^m(\tilde{y}_k) \\\text{and } H(T^m(\tilde{y}_i), \alpha_1) \cap D(T^m(\tilde{y}_i), \alpha_1) \text{ holds for each } i\}.\end{aligned}$$

LEMMA 2.22. *Let $Q$ be a cube of side $n$ and let $C_E \leq m \leq n^{1/9}$. Then if $|x_0 - x_1|_\infty \geq n^{2/9}$*

(2.32) $$\mathbb{P}_p(L(Q, m, x_0, x_1)^c) \leq c \exp(-c_1 n^{1/11}).$$

PROOF. Whereas $m$ is fixed in this argument, we write $T(\tilde{x})$ for $T^m(\tilde{x})$. Let $n'$ be such that $mn' \geq n \geq m(n' - 1)$. Let $\widetilde{Q}$ be a (macroscopic) cube of side $n'$ such that $Q \subset \bigcup\{T(\tilde{x}), \tilde{x} \in \widetilde{Q}\}$. Let $\tilde{x}_i$ be such that $x_i \in T(\tilde{x}_i)$ and let $s = |\tilde{x}_0 - \tilde{x}_1|_1$, so that $m(s - 1) \leq |x_0 - x_1|_\infty \leq m(s + 1)$ and $s \geq n^{1/9}$. Let $\widetilde{S} = \widetilde{S}(\psi^{(m)}, \widetilde{Q}, \lambda_0, \tilde{x}_0, \tilde{x}_1)$ be the event defined from the process $\psi^{(m)}$ in the same way as $S(\cdot)$ in (2.25) is for $\zeta$. Then, as $S$ is an increasing event, by Lemmas 2.20 and 2.21,

(2.33) $$\mathbb{P}_p(\widetilde{S}) \leq c \exp(-c'(n')^{1/9}) \leq c \exp(-c_2 n^{1/11}).$$

Let $\omega \in \widetilde{S}$ and let $\tilde{x}_0' = \tilde{y}_0, \ldots, \tilde{y}_k = \tilde{x}_1'$ be the open path (with respect to $\psi^{(m)}$) given by the event $\widetilde{S}$. Since $\widetilde{S}$ occurs we have $|\tilde{x}_i - \tilde{x}_i'|_\infty \leq (n')^{1/9}$ for $i = 0, 1$ and $k \leq \lambda |\tilde{x}_0 - \tilde{x}_1|_\infty$. Choose $x_i' \in T(\tilde{x}_i') \cap \mathcal{C}^\vee(Q^+)$; then

$$|x_i - x_i'|_\infty \leq m(1 + |\tilde{x}_i - \tilde{x}_i'|_\infty) \leq m(1 + (n')^{1/9}) < n^{2/9}.$$

Also, since $k \leq \lambda_0 s$, $mk \leq ms\lambda_0 \leq \lambda_0(m + |x_0 - x_1|_\infty) \leq 2\lambda_0 |x_0 - x_1|_\infty$. Thus $\omega \in L(Q, m, x_0, x_1)$ and using (2.33), this proves the lemma. □

Now let $\alpha_2 = (11(d + 2))^{-1}$ and let

$$\begin{aligned}L(Q) = H(Q, \alpha_2) \cap D(Q, \alpha_2) \cap \{L(Q, m, x, y) \text{ holds for every } x, y \in Q, \\\text{with } |x - y|_\infty \geq n^{2/9} \text{ and } C_E \leq m \leq n^{1/9}\}.\end{aligned}$$

THEOREM 2.23. *Let $Q$ be a cube of side $n$ and let $p > p_c$.*

(a) *There exists* $\mathbb{P}_p(L(Q)^c) \leq c_1 \exp(-c_2 n^{\alpha_2 \beta})$.

(b) *Let* $\omega \in L(Q)$ *and* $C_E \leq m \leq n^{1/9}$. *Then if* $x_0, x_1 \in Q \cap \mathcal{C}^\vee(Q^+)$ *with* $d_\omega(x_0, x_1) \geq \frac{1}{3} n^{1/4}$ *there exist* $x'_i \in \mathcal{C}^\vee(Q^+)$ *with* $d_\omega(x_i, x'_i) \leq \frac{1}{3} n^{1/4}$ *and a path* $\gamma = (z_0, \ldots, z_j)$ *in* $\mathcal{C}^\vee(Q^+)$ *that connects* $x'_0$ *and* $x'_1$ *such that:*

(i) *For each* $0 \leq l \leq j$, *the ball* $B_l = B_\omega(z_l, m/16)$ *is very good, with* $N_{B_l} \leq C_H m^{\alpha_1}$.

(ii) *There exists* $j \leq c_3 |x_0 - x_1|_\infty \leq C_F |x_0 - x_1|_1$.

PROOF. (a) This is immediate from the bounds in Lemma 2.13, Proposition 2.17 and Lemma 2.22.

(b) Since $D(Q, \alpha_2)$ occurs, $|x_0 - x_1|_\infty \geq c d_\omega(x_0, x_1) \geq n^{2/9}$ and so $\omega \in L(Q, m, x_0, x_1)$. Let $x'_i, \tilde{y}_0, \ldots, \tilde{y}_k$ be as in (2.31). Note that we can choose $x'_i$ to be within a distance $m/8$ of the center of the cubes $T(\tilde{y}_0)$ and $T(\tilde{y}_k)$. Then, by Proposition 2.17(d), $d_\omega(x_i, x'_i) \leq C_H((1 + n^{\alpha_2}) \vee |x_i - x'_i|_\infty) \leq c n^{2/9} \leq \frac{1}{3} n^{1/4}$.

We now show that the clusters $\mathcal{C}^\vee(Q_i)$ are all in $\mathcal{C}^\vee(Q^+)$. Consider first two adjacent cubes $T(y_i)$ and $T(y_{i+1})$. Since the event $R(T(y_i)) \cap R(T(y_{i+1}))$ occurs, the clusters $\mathcal{C}^\vee(T(y_i^+))$ and $\mathcal{C}^\vee(T(y_{i+1}^+))$ are connected. Thus there exists a $Q^+$ cluster $\mathcal{C}$ which contains each $\mathcal{C}^\vee(T(y_i^+))$, and so has diameter $D$ with $D \geq |x'_0 - x'_1|_\infty \geq |x_0 - x_1|_\infty - 2n^{1/9} \geq \frac{1}{2} n^{1/4}$. It now follows, as in Lemma 2.13(b), that $\mathcal{C} \subset \mathcal{C}^\vee(Q^+)$.

Since each event $D(Q_i, \alpha_1)$ holds, we can find a path $\gamma = \{z_0, \ldots, z_j\} \subseteq \bigcup T(y_i)$ that connects $x'_0$ and $x'_1$ with length $j \leq 2 C_H m k \leq c |x_0 - x_1|_\infty$. Each point $z_i$ is in a cube $Q_i$ for which $H(Q_i, \alpha_1) \cap D(Q_i, \alpha_1)$ occurs, so using Theorem 2.18(c), $B_i = B_\omega(z_i, m/16)$ is very good with $N_{B_i} \leq C_H m^{\alpha_1}$. □

LEMMA 2.24. *Let* $p > p_c$ *and, for each* $x \in \mathbb{Z}^d$, *let* $N_x$ *be the largest* $n$ *such that* $L(Q)$ *fails for some* $Q$ *with* $s(Q) = n$ *and* $x \in Q$. *Then*

(2.34) $$\mathbb{P}_p(N_x \geq n) \leq c_1 \exp(-c_2 n^{\alpha_2 \beta}).$$

PROOF. This is immediate from Theorem 2.23(a). □

**3. Upper bounds.** We now consider a connected graph $\mathcal{G} = (G, E)$ that satisfies the conditions (1.1) and (1.2). We use the notation of Section 1, and study the transition density $q_t(x, y)$ of the continuous time r.w. $Y_t$ on $\mathcal{G}$. Fix constants $C_V$, $C_P$ and $C_W$, and recall from Definition 1.7 the definition of good and very good balls, and of $N_B$. In this section the constants $c_i$ depend on the constants $d$, $C_0$, $C_V$, $C_P$ and $C_W$ in (1.1), (1.15) and (1.16).

As in Section 2, we assume without always stating it explicitly that the radius $R$ of a ball $B(x, R)$ is sufficiently large; that is, that $R \geq c_0 =$ $c_0(d, C_0, C_P, C_W)$. All the bounds in this section hold for balls $B(x, R)$ with $R \leq c_0$, with a suitable choice of the constants $c_i$ in the bounds, by elementary arguments.

We begin by investigating the on-diagonal decay of $q_t(x, x)$. We remark that a similar result was proved in Mathieu and Remy (2004), using an isoperimetric inequality directly. We give another proof here because it is quite short and also allows us to estimate the "initial time" $T_B$ directly in terms of $N_B$.

PROPOSITION 3.1 [see Mathieu and Remy (2004), Theorem 1.2]. *Let* $x_0 \in G$ *and let* $B = B(x_0, R)$ *be very good, so that* $N_B^{2+d} \leq R$. *Then*

$$(3.1) \quad q_t(x_1, x_1) \leq \frac{c_2}{t^{d/2}} \qquad \text{for } c_1 N_B^{2d} \leq t \leq \frac{R^2}{\log R} \text{ and } x_1 \in B(x_0, \tfrac{8}{9}R).$$

PROOF. Let $c_3 < c_{1.1.5}$, let $t_2 = c_3 R^2 / \log R$ and suppose that $t \leq t_2$. Then provided $R \geq c$, we have $t \log t \leq c_3 R^2$ and hence that $t \leq \exp(c_3 R^2/t)$. Fix $x_1 \in B(x_0, \tfrac{8}{9}R)$, write $f_t(x) = q_t(x_1, x)$ and let $\psi(t) = \int q_t(x_1, y)^2 \, d\mu = q_{2t}(x_1, x_1)$. Note that by (1.7),

$$(3.2) \qquad \frac{c_4}{(t \log t)^{d/2}} \leq \psi(t) = q_{2t}(x_1, x_1) \leq \frac{1}{c_5 t^{1/2}}, \qquad t \geq 1.$$

Using the discrete Gauss–Green formula,

$$\psi'(t) = 2 \sum_x f_t(x) \frac{\partial f_t(x)}{\partial t} = 2 \sum_x f_t(x) \mathcal{L} f_t(x) = -\sum_{x,y} a_{xy}(f_t(x) - f_t(y))^2$$

and, in particular, $\psi(t)$ is decreasing (and continuous).

Define $t_1$ so that $t_1^{1/2} = c_6 N_B^d$, where $c_6$ is chosen later, and choose $r(t)$ so that

$$(3.3) \qquad \frac{2}{c_5 c_6} \geq r(t)^d \psi(t) \geq \frac{1}{c_5 c_6}.$$

Then, if $t \geq t_1$, by (3.2), $\psi(t)^{-1} \geq c_5 c_6 N_B^d$ and so $r(t) \geq N_B$. If $t \leq t_2$, then by (3.2) and (3.3), $r(t)^d \leq c(t \log t)^{d/2} \leq c c_3^{d/2} R^2$, so if the constant $c_3$ is chosen small enough, we have $r(t) \leq R/18$.

Write $B' = B(x_0, 17R/18)$. Let $t \in [t_1, t_2]$, so that $r = r(t) \in [r_0, R/18]$, and let $B(y_i, r/2)$, $i = 1, \ldots, m$, be a maximal collection of disjoint balls with centers in $B'$. Set $B_i = B(y_i, r)$ and $B_i^* = B(y_i, C_W r)$. Note that $B' \subseteq \bigcup_i B_i \subseteq B$. If $x \in B \cap B_i^*$, we have $B(y_i, r/2) \subseteq B(x, r(1 + C_W))$, and so

$$C_0(1 + C_W)^d r^d \geq \mu(B(x, r(1 + C_W)))$$
$$\geq \sum_i \mathbb{1}_{\{x \in B_i^*\}} \mu(B(y_i, r/2)) \geq |\{i : x \in B_i^*\}| C_V 2^{-d} r^d.$$

Thus any $x \in B$ is in at most $c_7$ of the $B_i^*$.

The bounds on $r$ above imply that each $B(y_i, r)$ is good. So, applying the Poincaré inequality (1.16) to each $B_i$ and writing $\bar{f}_{t,i} = \mu(B_i)^{-1} \int_{B_i} f_t$, we have

$$-c_7 \psi'(t) \geq \sum_i \int_{E(B_i^*)} |\nabla f_t|^2$$

(3.4)
$$\geq C_P^{-1} r^{-2} \sum_i \int_{B_i} |f_t - \bar{f}_{t,i}|^2$$

$$= C_P^{-1} r^{-2} \sum_i \int_{B_i} f_t^2 - C_P^{-1} r^{-2} \sum_i \mu(B_i)^{-1} \left( \int_{B_i} f_t \right)^2.$$

By Lemma 1.1(b),

$$\int_{G-B'} f_t^2 \leq \left( \sup_{G-B'} f_t \right) \int_{G-B'} f_t \leq ct^{-d},$$

while, by (3.2), $\psi(t) \geq c_4 (t \log t)^{-d/2}$. So as $t \geq c_1$,

$$\sum_i \int_{B_i} f_t^2 \geq \int_{B'} f_t^2 = \psi(t) - \int_{G-B'} f_t^2 \geq \tfrac{1}{2} \psi(t).$$

Also, since $f_t$ has total mass 1,

$$\sum_i \mu(B_i)^{-1} \left( \int_{B_i} f_t \right)^2 \leq (C_V r^d)^{-1} \left( \sum_i \int_{B_i} f_t \right)^2 \leq c_7 (C_V r^d)^{-1} = c_8 r(t)^{-d}.$$

Combining these estimates, we obtain

(3.5) $\qquad -c_7 \psi'(t) \geq C_P^{-1} r(t)^{-2} (\tfrac{1}{2} \psi(t) - c_8 r(t)^{-d}).$

Now let $c_6 = (4c_5 c_8)^{-1}$ so that by the choice of $r = r(t)$ in (3.3),

$$-\psi'(t) \geq c_9 r(t)^{-2} \psi(t) \geq c_{10} \psi(t)^{1+2/d}.$$

Setting $\varphi(t) = \psi(t)^{-2/d}$ we have $\varphi'(t) \geq 2c_{10}/d$, from which it follows that

$$\varphi(t) \geq \varphi(t_1) + (2c_{10}/d)(t - t_1) \geq c_{11} t, \qquad 2t_1 \leq t \leq t_2.$$

Rearranging, this gives $\psi(t) \leq ct^{-d/2}$ for $2t_1 \leq t \leq t_2$. Since $\psi$ is decreasing it follows, by adjusting the constant $c$, that $q_t(x_1, x_1) = \psi(t/2) \leq c_{11} t^{-d/2}$ for $4t_1 \leq t \leq R^2/\log R$. $\square$

We need a bound for $y$ outside $B(x_0, R)$.

COROLLARY 3.2. *Let $x_0 \in G$ and let $B = B(x_0, R)$ be very good. Then*

(3.6) $\quad q_t(x_1, y) \leq \dfrac{c_1}{t^{d/2}} \qquad \text{for } c_2 N_B^{2d} \leq t \leq \dfrac{c_3 R^2}{\log R}, x_1 \in B(x_0, \tfrac{7}{9} R) \text{ and } y \in G.$

PROOF. If $y \in B(x_0, \frac{8}{9}R)$, then $q_t(x_1, y) \leq q_t(x_1, x_1)^{1/2} q_t(y, y)^{1/2} \leq ct^{-d/2}$ by Proposition 3.1. If $y \notin B(x_1, \frac{8}{9}R)$, then $d(x, y) \geq R/9$ and we use Lemma 1.1(a). □

For a very good ball $B$, let $T_B = c_{3.2.2} N_B^{2d}$ and $T_B' = c_{3.2.3} R^2 / \log R$.

REMARK. It is natural to ask if the bounds in (3.1) and (3.6) hold for $t \leq cR^2$ rather than just $t \leq cR^2 / \log R$. However, in this paper this restriction on $t$ does not matter, since we ultimately apply (3.6) in the situation where $B(x_0, R)$ is very good for all sufficiently large $R$.

We now use the method of Bass (2002) to obtain off-diagonal upper bounds. Following Nash (1958) and Bass (2002), we introduce the functions, for $x_1 \in G$, $t > 0$,

$$M(t) = M(x_1, t) = \sum_y d(x_1, y) q_t(x_1, y) \mu(y),$$

$$Q(t) = Q(x_1, t) = -\sum_y q_t(x_1, y) \log q_t(x_1, y) \mu(y).$$

We can extend $M$ and $Q$ to $t = 0$ by continuity: $M(0) = 0$, while since $q_t(x_1, x_1) \to \mu(x_1)^{-1}$ as $t \downarrow 0$, $Q(0) = \log \mu(x_1) \geq 0$.

LEMMA 3.3. *Let $B(x_0, R)$ be very good and let $x_1 \in B(x_0, \frac{7}{9}R)$.*

(a) *We have*

$$Q(x_1, t) \geq c_1 + \tfrac{1}{2} d \log t, \qquad T_B \leq t \leq T_B'.$$

(b) *We have*

$$M(x_1, t) \geq c_2 \exp(Q(x_1, t)/d), \qquad t \geq c_3.$$

PROOF. Fix $x_1 \in B(x_0, \frac{7}{9}R)$. Part (a) follows directly from the upper bound (3.6).

The proof of (b) is similar to that in Nash (1958) or Bass (2002). Let $0 < a < 1$, and set $D_0 = \{x_0\}$ and $D_n = B(x_0, 2^n) - B(x_0, 2^{n-1})$ for $n \geq 1$. Then using (1.1) to bound $\mu(D_n)$, we have, for $a \leq 2$,

$$\sum_{y \in G} \exp(-ad(x_1, y)) \mu(y) \leq \sum_{n=0}^{\infty} \sum_{y \in D_n} \exp(-a2^n) \mu(y)$$

$$\leq \sum_{n=0}^{\infty} C_0 2^{nd} \exp(-a2^n) \leq c_4 a^{-d}.$$

Now note that $u(\log u + \lambda) \geq -e^{-1-\lambda}$ for $u > 0$. So, setting $\lambda = ad(x_0, y) + b$, where $a \leq 2$,

$$-Q(x_1, t) + aM(x_1, t) + b = \sum_y q_t(x_1, y)(\log q_t(x_1, y) + ad(x_1, y) + b)\mu(y)$$

$$\geq -\sum_y \exp(-1 - ad(x_1, y) - b)\mu(y)$$

$$\geq -e^{-1-b} \sum_y \exp(-ad(x_1, y))\mu(y) \geq -c_5 e^{-b} a^{-d}.$$

Since $M(x_1, t) \geq P^{x_1}(Y_t \neq x_1)$, using (1.7) we have $M(x_1, t) \geq \frac{1}{2}$ when $t \geq c_5$. Setting $a = 1/M(x_1, t)$ and $e^b = M(x_1, t)^d = a^{-d}$, we obtain

$$-Q(x_1, t) + 1 + d \log M(x_1, t) \geq -c_4,$$

and rearranging gives (b). $\square$

PROPOSITION 3.4. *Let $x_0 \in G$ and let $B(x_0, R)$ be very good. Then*

(3.7) $$c_1 t^{1/2} \leq M(x_1, t) \leq c_2 t^{1/2} \qquad \text{for } x \in B(x_0, \tfrac{7}{9}R) \text{ and } T_B \log T_B \leq t \leq T'_B.$$

PROOF. For the moment we just write $Q(t)$ and $M(t)$. Set $f_t(x) = q_t(x_1, x)$ and let $b_t(x, y) = f_t(x) + f_t(y)$. We have

$$M'(t) = \sum_y d(x_1, y) \frac{\partial f_t(y)}{\partial t} \mu(y) = \sum_y d(x_1, y) \mathcal{L} f_t(y) \mu(y)$$

(3.8) $$= -\frac{1}{2} \sum_x \sum_y a_{xy}(d(x_1, y) - d(x_1, x))(f_t(y) - f_t(x))$$

$$\leq \frac{1}{2} \sum_x \sum_y a_{xy} |f_t(y) - f_t(x)|$$

$$\leq \frac{1}{2} \left( \sum_x \sum_y a_{xy} b_t(x, y) \right)^{1/2} \left( \sum_x \sum_y a_{xy} \frac{(f_t(y) - f_t(x))^2}{b_t(x, y)} \right)^{1/2}$$

(3.9) $$\leq c \left( \sum_x \sum_y a_{xy} \frac{(f_t(y) - f_t(x))^2}{f_t(x) + f_t(y)} \right)^{1/2}.$$

In the calculation above the use of the discrete Gauss–Green formula to obtain (3.8) is valid since, by (1.5), $q_t(x_1, \cdot)$ decays exponentially. Since we have, for $u, v > 0$,

$$\frac{(u-v)^2}{u+v} \leq (u-v)(\log u - \log v),$$

we deduce
$$M'(t)^2 \leq \sum_x \sum_y a_{xy}(f_t(y) - f_t(x))(\log f_t(y) - \log f_t(x)).$$

On the other hand [again using (1.5) and the discrete Gauss–Green formula],

(3.10)
$$\begin{aligned}Q'(t) &= -\sum_y (1 + \log f_t(y))\mathcal{L}f_t(y) \\ &= \tfrac{1}{2}\sum_x \sum_y a_{xy}(\log f_t(y) - \log f_t(x))(f_t(y) - f_t(x)) \geq \tfrac{1}{2}M'(t)^2.\end{aligned}$$

The remainder of this proof is similar to that in Nash (1958) or Bass (2002), except that we have to control the growth of $M$ for small $t$. Set $R(t) = d^{-1}(Q(t) - c_{3.3.1} - \tfrac{1}{2}d\log t)$, so that $R(t) \geq 0$ if $T_B \leq t \leq T'_B$. Define
$$T_0 = \begin{cases} 1, & \text{if } R(t) \geq 0 \text{ on } [1, T'_B], \\ \sup\{t \leq T'_B : R(t) < 0\}, & \text{otherwise.} \end{cases}$$

If $T_0 > 1$, then $T_0 \leq T_B$ and, by (3.10),
$$\begin{aligned}M(T_0) &= \int_0^{T_0} M'(s)\,ds \leq 2^{1/2}\int_0^{T_0} Q'(s)^{1/2}\,ds \\ &\leq 2^{1/2}\left(\int_0^{T_0} Q'(s)\,ds\right)^{1/2} T_0^{1/2} \\ &\leq c_3 T_0^{1/2}(Q(T_0) - Q(0))^{1/2} \\ &\leq c_3 T_0^{1/2}(c_{3.3.1} + \tfrac{1}{2}d\log T_0)^{1/2} \leq c_4(T_B \log T_B)^{1/2}.\end{aligned}$$

If $T_0 = 1$, then $M(T_0) = E^x d(x, Y_1) \leq c_5$ by elementary arguments.

By Lemma 3.3(b) and (3.10), if $T_0 < t < T'_B$, then
$$ct^{1/2}e^{R(t)} = e^{Q(t)/d} \leq M(t) \leq M(T_0) + 2^{1/2}\int_{T_0}^t Q'(s)^{1/2}\,ds$$
$$\leq M(T_0) + (2d)^{1/2}\int_{T_0}^t \left(R'(s) + \frac{1}{2s}\right)^{1/2} ds.$$

Using the inequality $(a+b)^{1/2} \leq b^{1/2} + a/(2b)^{1/2}$ gives

(3.11) $\quad ct^{1/2}e^{R(t)} \leq M(t) \leq M(T_0) + ct^{1/2} + c\int_{T_0}^t s^{1/2}R'(s)\,ds.$

Integrating by parts, and using the fact that $R \geq 0$ on $[T_0, T'_B]$, the final term in (3.11) is bounded by $c(1 + R(t)t^{1/2})$. Combining these estimates, for $T_B \leq t \leq T'_B$,
$$ct^{1/2}e^{R(t)} \leq M(t) \leq c(1 + R(t))t^{1/2} + c_4(T_B \log T_B)^{1/2}.$$

So for $T_B \log T_B \leq t \leq T'_B$,

$$ct^{1/2} e^{R(t)} \leq M(t) \leq c(1 + R(t))t^{1/2}.$$

Thus $R(t)$ is bounded and this implies (3.7). □

As in Bass (2002) we can use the moment bounds in (3.7) to obtain off-diagonal upper bounds on $q_t$ by the method of Barlow and Bass (1989, 1992). We define

$$\tau(x, r) = \inf\{t : Y_t \notin B(x, r)\}, \qquad x \in G, r > 0,$$

and begin by controlling the probability that $\tau(x, r)$ is small.

LEMMA 3.5. *Let $x_0 \in G$ and let $B(x_0, R)$ be very good. Let $c_1 N_B^d \times (\log N_B)^{1/2} \leq r \leq R$. Then*

$$(3.12) \quad P^x(\tau(x, r) < t) \leq \frac{1}{2} + \frac{c_2 t}{r^2} \qquad \text{for } x \in B(x_0, \tfrac{6}{9}R) \text{ and } 0 \leq t \leq \tfrac{1}{2} T'_B.$$

PROOF. Suppose first that $r < R/9$. Let $x \in B(x_0, \tfrac{6}{9}R)$, $A = B(x, r) \cup \partial B(x, r)$ and $\tau = \tau(x, r)$. Then if $T_B \log T_B \leq t \leq \tfrac{1}{2} T'_B$, since $A \subseteq B(x_0, \tfrac{7}{9}R)$,

$$\begin{aligned}
c_3 t^{1/2} &\geq E^x d(x, Y_{2t}) \geq E^x (d(x, Y_{t \wedge \tau}) - d(Y_{t \wedge \tau}, Y_{2t})) \\
&\geq E^x \mathbb{1}_{(\tau < t)} d(x, Y_\tau) - E^x (E^{Y_{t \wedge \tau}} d(Y_{t \wedge \tau}, Y_{2t - t \wedge \tau})) \\
&\geq P^x(\tau < t) r - \sup_{z \in A, s \leq t} E^z d(z, Y_{2t - s}) \\
&\geq P^x(\tau < t) r - c_3 t^{1/2}.
\end{aligned}$$

Thus

$$(3.13) \qquad\qquad P^x(\tau < t) \leq 2c_3 t^{1/2}/r.$$

Since $\lambda \leq \tfrac{1}{2}(1 + \lambda^2)$, (3.12) is immediate. If $t \leq T_B \log T_B$, then

$$P^x(\tau(x, r) < t) \leq P^x(\tau(x, r) < T_B \log T_B) \leq 2c_3 (T_B \log T_B)^{1/2} r^{-1} \leq \tfrac{1}{2},$$

provided $r \geq c(T_B \log T_B)^{1/2} = c' N_B^d (\log N_B)^{1/2}$. Finally if $R/9 \leq r \leq R$, we have $\tau(x, r) \geq \tau(x, R/9)$, so (adjusting the constant $c_2$) we deduce (3.12). □

REMARK. In the end we gain nothing useful by using the stronger bound (3.13).

We need the following estimate.

LEMMA 3.6 [Barlow and Bass (1989), Lemma 1.1]. *Let $\xi_1, \xi_2, \ldots, \xi_n$, $V$ be nonnegative r.v. such that $V \geq \sum_1^n \xi_i$. Suppose that for some $p \in (0,1)$, $a > 0$,*

$$P(\xi_i \leq t | \sigma(\xi_1, \ldots, \xi_{i-1})) \leq p + at, \qquad t > 0.$$

*Then*

$$\log P(V \leq t) \leq 2\left(\frac{ant}{p}\right)^{1/2} - n \log \frac{1}{p}.$$

PROPOSITION 3.7. *Let $x_0 \in G$ and let $B(x_0, R)$ be very good. If $x \in B(x_0, \frac{5}{9}R)$ and $t > 0$, $\rho > 0$ satisfy*

(3.14) $$\rho \leq R, \qquad c_1 N_B^d (\log N_B)^{1/2} \rho \leq t \quad \text{and} \quad t \leq T_B',$$

*then*

(3.15) $$P^x(\tau(x, \rho) < t) \leq c_2 \exp(-c_3 \rho^2 / t).$$

PROOF. Let $r_1 = c_{3.5.1} N_B^d (\log N_B)^{1/2}$. Suppose first that, in addition, $\rho < R/9$. Let $m \geq 1$ be chosen later, and let $s = t/m$ and $r = \lfloor \rho/m \rfloor$. Define stopping times

$$S_0 = 0, \qquad S_i = \inf\{t \geq S_{i-1} : d(Y_{S_{i-1}}, Y_t) = r\}, \qquad i \geq 1.$$

Set $\xi_i = S_i - S_{i-1}$ and write $\mathcal{F}_t = \sigma(Y_s, s \leq t)$ for the filtration of $Y$. By Lemma 3.5,

(3.16) $$P^x(\xi_i > u | \mathcal{F}_{S_{i-1}}) \leq \frac{1}{2} + \frac{c_4 u}{r^2}, \qquad u > 0,$$

provided $r_1 \leq r \leq R$, $u \leq T_B'$ and $Y_{S_{i-1}} \in B(x_0, \frac{6}{9}R)$. Whereas $d(Y_0, Y_{S_m}) \leq mr \leq \rho < R/9$, we have $S_m \leq \tau(x, \rho)$ and $Y_{S_j} \in B(x_0, \frac{6}{9}R)$ for $0 \leq j \leq m$. Using Lemma 3.6 and writing $p = \frac{1}{2}$, $a = c_4/r^2$, we deduce that

$$\log P^x(\tau(x, \rho) < t) \leq \log P^x(S_m < t) \leq 2(amt/p)^{1/2} - m \log p^{-1}.$$

Simplifying this expression, we obtain

(3.17) $$\log P^x(\tau(x, \rho) < t) \leq -c_5 m \left(1 - \left(\frac{c_6 t m}{\rho^2}\right)^{1/2}\right).$$

Let $\lambda = \rho^2/(2c_6 t)$. If we can choose $m \in \mathbb{N}$ with $\frac{1}{2}\lambda \leq m < \lambda$ and so that the estimate (3.16) is valid, then (3.17) implies (3.15).

If $\lambda \leq 1$, then, adjusting the constant $c_2$ appropriately, (3.15) is immediate. If $\lambda > 1$, then let $m = \lfloor \frac{1}{2}\lambda \rfloor + 1$. Since then $m \geq 1$, we have $s \leq t \leq T_B'$ and $r \leq \rho \leq R$, while the condition $\rho \leq c_6 t/r_1$ ensures that $r \geq r_1$.

Finally let $\rho$ satisfy (3.14) but with $\rho \geq R/9$. Then (adjusting $c_1$ if necessary) we can apply the argument above to $\rho_0 = R/9$ and adjust the constant $c_3$ to obtain (3.15). □

THEOREM 3.8. *Let $x_0 \in G$ and let $B(x_0, R)$ be very good. Let $x \in B(x_0, \frac{1}{2}R)$, let $y \in G$ and assume that*

$$N_B^{2d+1} \vee d(x,y) \leq t \leq \frac{c_1 R^2}{\log R}. \tag{3.18}$$

*Then*

$$q_t(x,y) \leq c_2 t^{-d/2} \exp(-c_3 d(x,y)^2/t). \tag{3.19}$$

PROOF. Let $D = d(x,y)$. Using (1.5) we have, since $D \leq t$,

$$q_t(x,y) \leq c_5 \exp(-2c_4 D^2/t).$$

If $t \log t \leq 2c_4 d^{-1} D^2$, then $\exp(-c_4 D^2/t) \leq t^{-d/2}$ and we deduce that

$$q_t(x,y) \leq c_5 t^{-d/2} \exp(-c_4 D^2/t).$$

Suppose therefore that $t \log t \geq 2c_4 d^{-1} D^2$. Note that this implies that $y \in B(x, \frac{5}{9}R)$, provided $R \geq c$ and $c_1$ in (3.18) is chosen small enough. Let $A_x = \{z : d(x,z) \leq d(y,z)\}$, $A_y = G - A_x$, $s = t/2$ and $\rho = D/2$. Note that $B(x, \rho) \subseteq A_x$. Then

$$\begin{aligned}\mu(x) P^x(Y_t = y) \\ = \mu(x) P^x(Y_t = y, Y_s \in A_y) + \mu(x) P^x(Y_t = y, Y_s \in A_x).\end{aligned} \tag{3.20}$$

To bound the first term in (3.20) we write

$$\begin{aligned}P^x(Y_t = y, Y_s \in A_y) &= P^x(\tau(x,\rho) < s, Y_s \in A_y, Y_t = y) \\ &\leq E^x(\mathbb{1}_{\{\tau(x,\rho) < s\}} P^{Y_\tau}(Y_{t-\tau} = y)) \\ &\leq P^x(\tau(x,\rho) < s) \sup_{z \in \partial B(x,\rho), s \leq t} q_{2t-s}(z,y) \mu(y).\end{aligned} \tag{3.21}$$

Since $2T_B < \frac{1}{2} N_B^{2d+1} \leq s < T_B'$, by Corollary 3.2 the second term in (3.21) is bounded by $ct^{-d/2}$. To control the first term we use Proposition 3.7. We have $\rho < D < R$ and $s < T_B'$, while, since $t \geq N_B^{2d+1}$,

$$c_{3.7.1} N_B^d (\log N_B)^{1/2} \rho \leq c N_B^d (\log N_B)^{1/2} (t \log t)^{1/2} \leq \tfrac{1}{2} t = s,$$

the three conditions in (3.14) are satisfied and, by (3.15),

$$P^x(Y_t = y, Y_s \in A_y) \leq c t^{-d/2} \exp(-c' D^2/s).$$

By symmetry the second term in (3.20) equals

$$\mu(y) P^y(Y_t = x, Y_s \in A_x) \tag{3.22}$$

and so can be bounded in the same way as the first term. Combining these estimates completes the proof. □

**4. A weighted Poincaré inequality.** While the weak Poincaré inequality of Section 2 is enough for upper bounds on the transition density, to obtain lower bounds we need a weighted Poincaré inequality, which we derive using the methods of Jerison (1986) and Saloff-Coste and Stroock (1991).

We continue with the notation and assumptions of the previous section. Fix $x_0 \in G$, fix $R \in \mathbb{N}$ and let $B = B(x_0, R)$ be a very good ball with $R_0 = N_B \leq R^{1/(1+d)}$. For each $x, y \in G$ we write $\gamma(x, y)$ for a shortest path $x = z_0, \ldots, z_{d(x,y)} = y$ between $x$ and $y$.

We begin with a Whitney decomposition of $B$, which we need to adapt to our situation. We have two differences from Jerison (1986), which both arise on small length scales. The first—minor—difficulty is that in our discrete setting we cannot use balls of size smaller than 1. The second difficulty is that we do not have any volume doubling estimate for balls smaller than $R_0$.

Let $(X, d)$ be the metric space obtained as the "cable system" of $\mathcal{G}$. This is the metric space obtained by replacing each edge $e$ by a copy of $(0, 1)$, linked in the obvious way at the vertices $x \in G$. We define a measure $\tilde{\mu}$ on $X$ by taking $\tilde{\mu}$ to be Lebesgue measure on each cable. See, for example, Barlow and Bass (2004) for further details of this construction. We write $\widetilde{B}(x, r)$ for balls in $X$. Since $R \in \mathbb{N}$, the boundary of $\widetilde{B}$ is contained in $G$. For $x \in \widetilde{B} = \widetilde{B}(x_0, R)$ we write $\rho(x) = d(x, \widetilde{B}^c)$. Note that if $x \in G$, then $\rho(x) = d(x, G - B)$. We frequently use the inequality

$$|\rho(x) - \rho(y)| \leq d(x, y). \tag{4.1}$$

Let $\lambda \geq 10^3 \vee (21 C_W)$ and let $10 \leq K \leq \lambda/10$ be fixed constants. We can assume that $R_0 > \lambda$.

LEMMA 4.1. *There exists a sequence of disjoint balls $\widetilde{B}_i = \widetilde{B}(x_i, r_i)$, $i \geq 1$, such that $r_1 \geq r_2 \geq \cdots$ and:*

(a) *There exists $\widetilde{B} = \bigcup_{i=1}^\infty \widetilde{B}(x_i, 2r_i)$.*
(b) *For each $i$, $\rho(x_i) = \lambda r_i$.*
(c) *If $y \in \widetilde{B}(x_i, K r_i)$, then*

$$(\lambda - K) r_i \leq \rho(y) \leq (\lambda + K) r_i. \tag{4.2}$$

PROOF. This is standard. We start by choosing a ball $\widetilde{B}_1$ of maximal size that satisfies (b) and continue, so that $\widetilde{B}_n$ is chosen to be the largest ball contained in $\widetilde{B} - \bigcup_1^{n-1} \widetilde{B}_i$ that satisfies (b). To prove (a), suppose $y \notin \bigcup \widetilde{B}(x_i, 2r_i)$. Since $\tilde{\mu}(\widetilde{B}) < \infty$ and $\tilde{\mu}(\widetilde{B}_i) \geq r_i$, we must have $r_i \to 0$. Let $t_i = d(y, x_i) \geq 2r_i$, $t = \inf t_i$. Then by (4.1),

$$\rho(y) \leq \rho(x_i) + t_i = \lambda r_i + t_i \leq (1 + \tfrac{1}{2}\lambda) t_i,$$

so that $t \geq c\rho(y) = t' > 0$, contradicting the definition of $\widetilde{B}_i$ if $r_i < t$.

Part (c) follows immediately from (4.1) and (b). □

We now adapt this construction to our discrete setup. Let $N$ be defined by $r_N \geq R_0 + 1 > r_{N+1}$. For each $i \leq N$ the center $x_i$ of $\widetilde{B}_i = \widetilde{B}(x_i, r_i)$ lies on a cable $[y_i, y'_i]$, where $y_i, y'_i \in G$. We label these so that $y_i$ is the point in $G$ closest to $x_i$ and we set $s_i = r_i - d(x_i, y_i)$. Then we have $B(y_i, s_i) \subseteq \widetilde{B}(x_i, r_i) \cap G$ and $r_i \geq s_i \geq r_i - \frac{1}{2} > R_0$.

We set $\lambda_1 = \lambda - 2K$ and $\lambda_2 = \lambda + 2K$.

LEMMA 4.2. *The sequence of disjoint balls $B_i = B(y_i, s_i)$, $1 \leq i \leq N$, defined above satisfies the following statements.*

(a) *For each $i \leq N$,*

$$(4.3) \qquad \lambda s_i - \tfrac{1}{2} \leq \rho(y_i) \leq \tfrac{1}{2}(1 + \lambda) + \lambda s_i.$$

(b) *If $x \in B - \bigcup_{i=1}^{N} B(y_i, 3s_i)$, then $\rho(x) < \lambda_2 R_0$. Furthermore,*

$$(4.4) \quad B(x_0, R - \lambda_2 R_0) \subseteq \bigcup_{i=1}^{N} B(y_i, 3s_i - 1) \subseteq \bigcup_{i=1}^{N} B(y_i, \lambda_1 s_i) \subseteq B(x_0, R).$$

(c) *If $x \in B(y_i, K s_i)$, then*

$$(4.5) \quad \lambda_1 s_i \leq (\lambda - K) s_i - \tfrac{1}{2} \leq \rho(x) \leq (\lambda + K) s_i + \tfrac{1}{2}(1 + \lambda) \leq \lambda_2 s_i.$$

(d) *There exists a constant $c_1$ such that*

$$(4.6) \qquad |\{i \leq N : x \in B(y_i, K s_i)\}| \leq c_1.$$

PROOF. Since $\rho(x_i) = \lambda r_i$ and $|\rho(y_i) - \rho(x_i)| \leq \frac{1}{2}$, (a) is immediate.

(b) Let $x \in B$. Then $x \in \widetilde{B}(x_i, 2r_i)$ for some $i$ and so $\rho(x) \leq (2 + \lambda)r_i$. If $i > N$, then $r_i < 1 + R_0$ and so $\rho(x) \leq (2 + \lambda)(1 + R_0) < \lambda_2 R_0$, which implies that $d(x_0, x) > R - \lambda_2 R_0$. So, if $\rho(x) \geq \lambda_2 R_0$, then $x \in \widetilde{B}(x_i, 2r_i)$ for some $i \leq N$. We then have $d(x, y_i) \leq d(x, x_i) + d(x_i, y_i) < 2s_i + \frac{3}{2}$. Since each $s_i > 3$, this implies that $x \in \bigcup B(y_i, 3s_i - 1)$.

The final inclusion in (4.4) is immediate from (a) and the first follows from the inequality $d(x_0, x) + \rho(x) \geq R$.

(c) This is immediate from (a) and (4.1).

(d) If $x \in B(y_i, K s_i)$, then by (c), $\rho(x) \geq \lambda_1 s_i$ and $B_i \subseteq B(x, (1+K)c_1\rho(x))$, where $c_2 = (1+K)/\lambda_1$. Also, $\mu(B_i) \geq C_V s_i^d \geq C_V \lambda_2^{-1} \rho(x)^d$. So writing $I = \{i : x \in \widetilde{B}(y_i, K s_i)\}$, we have

$$(4.7) \quad C_0 c_2^d \rho(x)^d \geq \mu(B(y, c_2 \rho(x))) \geq \sum_{i \in I} \mu(B_i) \geq |I| C_V \lambda_2^{-d} \rho(x)^d,$$

which proves (4.6). □

Let
$$B'_i = B(y_i, 3s_i), \qquad 1 \le i \le N.$$

Let $\eta = 2\lambda_2$ and set
$$B''_i = B(y_i, 10s_i) \qquad \text{if } s_i \ge \eta R_0.$$

If $s_i < \eta R_0$, we call $B_i$ a *boundary ball* and define $B''_i$ to be the connected component of $B(y_i, 2\lambda s_i) \cap B$ which contains $y_i$. (While balls are connected, the intersection of two balls need not be.) We relabel the balls $B_i$ so that $x_0 \in B_1$, and $B_i$ is a boundary ball if and only if $M+1 \le i \le N$.

LEMMA 4.3. (a) *There exists* $B = (\bigcup_{i=1}^M B'_i) \cup (\bigcup_{i=M+1}^N B''_i)$.
(b) *There exists a constant $c_1$ such that for any* $x \in B$,

(4.8) $$|\{i \le N : x \in B''_i\}| \le c_1.$$

PROOF. (a) Suppose $x \in B$ but $x \notin \bigcup_i B'_i$. Then $\rho(x) < \lambda_2 R_0 = t$. Now choose $x' \in \gamma(x, x_0)$ with $1 + t \ge \rho(x') > t$ and choose $x'' \in \gamma(x_0, x')$ with $d(x', x'') = 1$. Then $x' \in B(y_i, 3s_i - 1)$ for some $i$, so $\lambda s_i < \rho(x'') < t < \lambda_2 R_0$. Hence
$$R_0 < s_i \le \frac{\lambda_2 R_0}{\lambda_1} < \eta R_0,$$

so that $B_i$ is a boundary ball. Now $d(x, y_i) \le d(x, x') + 3s_i - 1 \le t + 3s_i < 2\lambda s_i$, which proves that $x \in B(y_i, 2\lambda s_i)$. The same argument proves that each $y \in \gamma(x, x')$ is also in $B(y_i, 2\lambda s_i)$. Hence $x$ is connected to $x'$ (and so $y_i$) by a path in $B(y_i, 2\lambda s_i) \cap B$, and so $x \in B''_i$.

(b) Since $K \ge 10$, we have a bound on $|\{i : x \in B(y_i, 10s_i)\}|$. So it is enough to control $|I'|$, where $I' = \{i : x \in B(y_i, 2\lambda s_i), s_i < \eta R_0\}$. The argument is almost exactly the same as in Lemma 4.2(d): If $i \in I'$, then $s_i < \eta R_0$ and $d(x, y_i) \le 2\lambda_2 s_i < 2\lambda_2 \eta R_0$. So $B_i \subseteq B(x, cR_0)$ and we use volume bounds as in (4.7). □

For each $i$ define
$$\mathcal{F}(i) = \{j : \gamma(x_0, y_i) \cap B(y_j, Ks_j) \ne \varnothing\},$$
$$\mathcal{F}(i, r) = \{j \in \mathcal{F}(i) : r \le s_j \le 2r\}.$$

LEMMA 4.4. (a) *If* $j \in \mathcal{F}(i)$, *then* $s_j \ge \frac{1}{2}\lambda_1 \lambda_2^{-1} s_i \ge \frac{1}{4} s_i$.

(b) *If* $j \in \mathcal{F}(i)$, *then* $d(y_i, y_j) \le (K + \lambda_2) s_j$.
(c) *There exists* $|\mathcal{F}(i, r)| \le c_1 C_0 / C_V$.

PROOF. Let $j \in \mathcal{F}(i)$, so that there exists $x \in B(y_j, Ks_j) \cap \gamma(x_0, y_i)$. Since $x \in \gamma(x_0, y_i)$, $d(x, y_i) = d(x_0, y_i) - d(x_0, x) < R - d(x_0, x) \leq \rho(x)$. Thus using (4.5),

$$\lambda_1 s_i \leq \rho(y_i) \leq d(x, y_i) + \rho(x) \leq 2\rho(x) \leq 2\lambda_2 s_j,$$

which proves (a).

For (b) note that $d(y_i, y_j) \leq d(y_i, x) + Ks_j \leq \rho(x) + Ks_j$.

From the estimates above, if $j \in \mathcal{F}(i, r)$, then $y_j \in B(y_i, 2(K + \lambda_2)r)$, so that $B(y_j, s_j) \subseteq B(y_i, 3\lambda r)$. Hence

$$c_2 C_0 r^d \geq \sum_{j \in \mathcal{F}(i,r)} \mu(B_j) \geq |\mathcal{F}(i, r)| C_V r^d,$$

proving (c). □

COROLLARY 4.5. *There exists $|\mathcal{F}(i)| \leq c_1 \log(R/s_i)$.*

PROOF. The proof follows easily from Lemma 4.4(c). □

Now let

$$\mathcal{F}^*(j) = \{i : j \in \mathcal{F}(i)\},$$
$$\mathcal{F}^*(j, r) = \{i : j \in \mathcal{F}(i), r \leq s_i \leq 2r\}.$$

LEMMA 4.6. *There exists $\alpha = \alpha(d, C_0, C_V) > 0$ such that for each $1 \leq j \leq N$ we have, for $r > R_0$,*

$$\sum_{i \in \mathcal{F}^*(j,r)} \mu(B_i') \leq c\mu(B_j') \left(\frac{r}{s_j}\right)^\alpha.$$

PROOF. This argument runs along the same lines as the proof of Lemma 5.9 in Jerison (1986). Note first that we can assume that $r \leq 4s_j$, since if $i \in \mathcal{F}^*(j, r)$, then by Lemma 4.4(a) we have $s_i \leq 4s_j$.

Write $\partial B = \{y : d(x_0, y) = R\}$. Fix $j$. Using (4.5) we have $\rho(y_j) \leq \lambda_2 s_j$, so we can choose $z \in \partial B$ with $d(y_j, z) \leq \lambda_2 s_j$. Set $t = (4 + K + 4\lambda_2)s_j$ and for $u > 0$, let

$$\Lambda(u) = B(z, t + 2u) \cap \{x \in B : \rho(x) \leq u\}.$$

Suppose that $i \in \mathcal{F}^*(j, r)$. By Lemma 4.4(b), $d(y_i, y_j) \leq (K + \lambda_2)s_j$, so that $d(y_i, z) \leq (K + 2\lambda_2)s_j$ and

$$B_i' \subseteq B(z, (K + 2\lambda_2 + 4)s_j).$$

By (4.5), $\rho(x) \le 2\lambda_2 r$ on $B_i'$, so $B_i' \subseteq \Lambda(2\lambda_2 r)$. Whereas the balls $B_i'$ are disjoint,

$$\sum_{i \in \mathcal{F}^*(j,r)} \mu(B_i') \le \mu(\Lambda(2\lambda_2 r)). \tag{4.9}$$

Now fix $u > 0$ and choose a maximal set of points $\{z_1, \ldots, z_m\} \subseteq \partial B$ such that $B(z_k, u)$ are disjoint and $d(z, z_k) \le t + u$ for each $k$. We next show that

$$\Lambda(u/4) \subseteq \bigcup_{k=1}^m B(z_k, 3u). \tag{4.10}$$

Let $x \in \Lambda(u/4)$, so that $d(x, z) \le t + u/2$ and there exists $z' \in \partial B$ with $d(x, z') \le u/4$. Hence $d(z', z) \le t + \frac{3}{4}u < t + u$. Whereas $\{z_1, \ldots, z_m\}$ is maximal, we must have $d(z', z_k) < 2u$ for some $k$. Thus $d(x, z_k) < 2u + \frac{1}{4}u < 3u$, proving (4.10).

For each $z_j$ we have $d(x_0, z_j) = R$ by construction. Choose $w_j$ on $\gamma(x_0, z_j)$ such that

$$\tfrac{1}{2}u < d(z_j, w_j) \le \tfrac{2}{3}u;$$

this is possible provided $6 < u < R$. We have $d(w_k, w_l) > d(z_k, z_l) - \frac{4}{3}u \ge \frac{2}{3}u$, so the balls $B(w_k, \frac{1}{4}u)$ are disjoint. The choice of $w_k$ implies that $\rho(w_k) = d(w_k, z_k) > \frac{1}{2}u$ and therefore

$$B(w_k, \tfrac{1}{4}u) \cap \Lambda(u/4) = \varnothing. \tag{4.11}$$

We also have

$$B(w_k, \tfrac{1}{4}u) \subseteq \Lambda(u). \tag{4.12}$$

To check this, if $x \in B(w_k, \frac{1}{4}u)$, then $\rho(x) \le d(x, z_k) \le \frac{3}{4}u$, while $d(x, z) \le d(x, z_k) + d(z_k, z) \le \frac{3}{4}u + t + u < t + 2u$. By (4.10), we deduce that

$$\mu(\Lambda(u/4)) \le \sum_{k=1}^m \mu(B(z_k, 3u)) \le m C_0 (3u)^d, \tag{4.13}$$

while by (4.11) and (4.12),

$$\mu(\Lambda(u)) \ge \mu(\Lambda(u/4)) + \sum_{k=1}^m \mu(B(w_k, \tfrac{1}{4}u)). \tag{4.14}$$

So, provided $\frac{1}{4}u \ge R_0$,

$$\mu(\Lambda(u)) \ge \mu(\Lambda(u/4)) + mC_V(\tfrac{1}{4}u)^d \ge (1+c_1)\mu(\Lambda(u/4)). \tag{4.15}$$

Note that the constant $c_1$ here depends only on $C_0$, $C_V$ and $d$.

Now let $R_0 \leq r \leq 4s_j$. Choose $n \in \mathbb{Z}_+$ as large as possible so that $4^n r \leq 4s_j$. Then

(4.16) $\quad \mu(\Lambda(2\lambda_2 r)) \leq (1+c_1)^{-n} \mu(\Lambda(2\lambda_2 4^n r)) \leq (1+c_1)^{-n} \mu(\Lambda(2t)).$

Set $\alpha = \log(1+c_1)/\log 4$. Then $(1+c_1)^n \geq (s_j/r)^\alpha$. We have $\mu(\Lambda(2t)) \leq \mu(B(z,5t)) \leq C_0 (5t)^d$ and also $\mu(B'_j) \geq C_V s_j^d \geq ct^d$. Combining this with (4.16) and (4.9) completes the proof of the lemma. $\square$

PROPOSITION 4.7. *For each $1 \leq j \leq N$ we have*
$$\sum_{i \in \mathcal{F}^*(j)} |\mathcal{F}(i)| \frac{\mu(B'_i)}{\mu(B'_j)} \leq c_1 \log\left(\frac{R}{s_j}\right).$$

PROOF. We can write
$$\mathcal{F}^*(j) = \bigcup_{n=-1}^\infty \mathcal{F}^*(j, 2^{-n} s_j).$$

Hence using Corollary 4.5 and Lemma 4.6,

$$\sum_{i \in \mathcal{F}^*(j)} |\mathcal{F}(i)| \mu(B'_i) \leq \sum_{n=-1}^\infty \sum_{i \in \mathcal{F}^*(j, 2^{-n} s_j)} |\mathcal{F}(i)| \mu(B'_i)$$
$$\leq c \sum_{n=-1}^\infty \log(2^n R/s_j) 2^{-n\alpha} \mu(B_j) \leq c\mu(B_j) \log(R/s_j). \quad \square$$

Set

(4.17) $\quad \varphi(y) = \left(\frac{R \wedge \rho(y)}{R}\right)^2, \qquad y \in G.$

For any set $A$, let $\hat\mu(A) = \int_A \varphi\, d\mu$ and $\bar{f}_A = \hat\mu(A)^{-1} \int_A f\varphi\, d\mu$. For an edge $e = \{x,y\}$ define $\tilde\varphi(e) = \varphi(x) \wedge \varphi(y)$. Note that if $e \in E(B)$, then $\tilde\varphi(e) \geq R^{-2}$.

THEOREM 4.8. *Let $B = B(x_0, R)$ be very good and let $R_0 = N_B \leq R^{1/(d+2)}$. Then*

(4.18) $\quad \int_B (f(x) - \bar{f}_B)^2 \varphi(x)\, d\mu \leq c_1 R^2 \int_{E(B)} |\nabla f|^2 \tilde\varphi\, d\nu.$

PROOF. We follow the proof in Saloff-Coste and Stroock (1991), but need some extra care close to the boundary of $B$.

For $1 \leq i \leq N$ set
$$B_i^* = \begin{cases} B(y_i, 10C_W s_i), & 1 \leq i \leq M, \\ B''_i, & M+1 \leq i \leq N. \end{cases}$$

Whereas $10C_W s_i \leq \frac{1}{2}\lambda s_i - 1$, we have $B_i^* \subseteq B$ for each $i \leq M$, while $B_i^* \subseteq B$ by definition if $M+1 \leq i \leq N$. Note that for any ball $B(y_i, cs_i)$ with $c \leq \frac{1}{2}\lambda$, we have, from (4.17),

$$(4.19) \qquad \varphi(x) \leq \frac{3\lambda_2}{2\lambda_1}\varphi(y) \qquad \text{for all } x, y \in B(y_i, cs_i).$$

Let $P_j$ be the best constant in the weighted Poincaré inequality

$$(4.20) \qquad \int_{B_j''} (f(x) - \bar{f}_{B_j})^2 \varphi(x)\,d\mu \leq P_j \int_{E(B_j^*)} |\nabla f|^2 \tilde\varphi\,d\nu.$$

Then for $j \leq M$, as each $B_j''$ is good, we have, using (4.19), $P_j \leq cC_P s_j^2$. If $M < j \leq N$, then, using Corollary 1.5(a),

$$(4.21) \qquad P_j \leq 2\mu(B_j'')^2 \sup_{x,y \in B_j''} \frac{\varphi(x)}{\varphi(y)} \leq c(C_0 \lambda R_0^d)^2 (\lambda R_0)^2 \leq c' R_0^{2d+2}.$$

Fix for the moment a ball $B_n$. We define a sequence of balls $D_i = B(w_i, t_i)$, $1 \leq i \leq L_n$, with $D_1 = B_1$ and $D_{L_n} = B_n$, as follows. Suppose $D_1, \ldots, D_{k-1}$ have been defined. Let $z_k$ be the point in $\gamma(x_0, y_n) \cap \bigcup_{i=1}^{k-1} D_i'$ which is furthest from $x_0$. (If $D_i = B_j$, then $D_i' = B_j'$ and $D_i'' = B_j''$.) If $z_k = y_n$, then we let $D_k = B_n$, $L_n = k$ and stop. Suppose that $z_k \neq y_n$. If $z_k \in B(y_i, 3s_i - 1)$ for some $i$, then we take $D_k = B_i$ and continue. (We choose the largest such $i$ if this $i$ is not unique.) Note that $D_k$ must be distinct from $D_1, \ldots, D_{k-1}$.

Finally, suppose $z_k \neq y_n$ and $z_k \notin \bigcup B(y_i, 3s_i - 1)$. Then $\rho(z_k) < \lambda_2 R_0$ and $\rho(y_n) < d(y_n, z_k) + \rho(z_k) < 2\lambda_2 R_0$. Hence $s_n < 2\lambda_2 \lambda_1^{-1} R_0 < \eta R_0$, so that $B_n$ is a boundary ball. We also have $\rho(w_{k-1}) \leq 3t_{k-1} + \rho(z_k)$, so that $D_{k-1}$ is also a boundary ball. In this case we take $D_k = B_n$, $L_n = k$ and stop. Note that each $D_k$ is a ball in $\{B_j : j \in \mathcal{F}(n)\}$, so that $L_n \leq |\mathcal{F}(n)|$.

We now show that

$$(4.22) \qquad D_{k-1}' \cup D_k' \subseteq D_{k-1}'' \cap D_k'', \qquad 2 \leq k \leq L_n.$$

First, if $z_k \neq y_n$ and $z_k \notin \bigcup B(y_i, 3s_i - 1)$, then both $D_{k-1}$ and $D_k$ are boundary balls, and $d(w_{k-1}, w_k) \leq 3t_{k-1} + \lambda_2 R_0$, from which (4.22) follows easily. Now let $z_k \in B(w_k, 3t_k - 1)$. Then $d(w_{k-1}, z_k) < 3t_{k-1}$ and so $d(w_{k-1}, w_k) < 3t_{k-1} + 3t_k$. Since $z_k \in B(w_j, 4t_j)$ for $j = k-1, k$, by (4.5),

$$\lambda_2(t_k \vee t_{k-1}) \geq \rho(z_k) \geq \lambda_1(t_k \wedge t_{k-1}).$$

Since $\lambda_2/\lambda_1 < 10/9$ this implies (4.22).

Let $\bar f_1 = \bar f(B'_1)$. Then

$$\int_{B''_n} (f(x)-\bar f_1)^2 \varphi(x)\, d\mu$$

(4.23)
$$= \int_{B''_n} \left( f(x) - \bar f_1 + \sum_{k=1}^{L_n-1} (\bar f(D''_k) - \bar f(D''_{k+1})) \right)^2 \varphi(x)\, d\mu$$

$$\le L_n \int_{B''_n} (f(x)-\bar f_1)^2 \varphi(x)\, d\mu + L_n \sum_{i=1}^{L_n-1} (\bar f(D''_k) - \bar f(D''_{k+1}))^2 \hat\mu(B''_n).$$

To control the terms in (4.23) note that

$$\hat\mu(B''_n)(\bar f(D''_k) - \bar f(D''_{k+1}))^2$$

$$= \frac{\hat\mu(B''_n)}{\hat\mu(D''_k \cap D''_{k+1})} \int_{D''_k \cap D''_{k+1}} (\bar f(D''_k) - \bar f(D''_{k+1}))^2 \varphi\, d\mu$$

$$\le 2\frac{\hat\mu(B''_n)}{\hat\mu(D''_k \cap D''_{k+1})} \left( \int_{D''_k} (f - \bar f(D''_k))^2 \varphi\, d\mu + \int_{D''_{k+1}} (f - \bar f(D''_{k+1}))^2 \varphi\, d\mu \right)$$

$$\le 2\frac{\hat\mu(B''_n)}{\hat\mu(D'_k)} \left( \int_{D''_k} (f - \bar f(D''_k))^2 \varphi\, d\mu + 2\frac{\hat\mu(B''_n)}{\hat\mu(D'_{k+1})} \int_{D''_{k+1}} (f - \bar f(D''_{k+1}))^2 \varphi\, d\mu \right).$$

Using (4.19) and Lemma 4.4(a),

$$\frac{\hat\mu(B''_n)}{\hat\mu(D'_k)} \le c\frac{\varphi(y_n)^2 \mu(B''_n)}{\varphi(w_k)^2 \mu(D'_k)} \le c'\frac{\mu(B''_n)}{\mu(D'_k)} \le c''\frac{\mu(B'_n)}{\mu(D'_k)}.$$

Combining these estimates we obtain, from (4.23),

(4.24)
$$\int_{B''_i} (f(x)-\bar f_1)^2 \varphi(x)\, d\mu \le c|\mathcal F(i)| \sum_{j\in\mathcal F(i)} \frac{\mu(B'_i)}{\mu(B'_j)} \int_{B''_j} (f - \bar f(B''_j))^2 \varphi\, d\mu$$

$$\le c|\mathcal F(i)| \sum_{j\in\mathcal F(i)} \frac{\mu(B'_i)}{\mu(B'_j)} P_j \int_{E(B^*_j)} |\nabla f|^2 \varphi\, d\mu.$$

Summing inequalities (4.24) gives

$$\int_B (f(x)-\bar f_1)^2 \varphi(x)\, d\mu$$

$$\le c\sum_{i=1}^N |\mathcal F(i)|\mu(B'_i) \sum_{j\in\mathcal F(i)} \mu(B'_j)^{-1} P_j \int_{E(B^*_j)} |\nabla f|^2 \tilde\varphi\, d\nu$$

$$= c\sum_{j=1}^N \left( \sum_{i\in\mathcal F^*(j)} |\mathcal F(i)|\frac{\mu(B'_i)}{\mu(B'_j)} \right) P_j \int_{E(B^*_j)} |\nabla f|^2 \tilde\varphi\, d\nu$$

$$\le c\sum_{j=1}^N \log\left(\frac{R}{s_j}\right) P_j \int_{E(B^*_j)} |\nabla f|^2 \tilde\varphi\, d\nu.$$

If $j \leq M$, then $P_j \leq cs_j^2 C_P$ and so since $t^p \log(R/t) \leq cR^p$ for $t > 0$ we have $P_j \log(R/s_j) \leq cC_P R^2$. If $M + 1 \leq j \leq N$, then $R_0 \leq s_j \leq \eta R_0$ and $P_j \leq cR_0^{2d+2}$, and so since $R_0^{d+2} \leq R$,

$$\log(R/s_j)P_j \leq cR_0^{2d+2}\log(R/R_0) \leq cR^2.$$

So since any edge in $B$ is contained in at most $c$ of the $B_i^*$,

$$\int_B (f(x) - \bar{f}_1)^2 \varphi(x)\, d\mu \leq cC_P R^2 \sum_{j=1}^N \int_{E(B_j^*)} |\nabla f|^2 \tilde{\varphi}\, d\nu$$

$$\leq c'C_P R^2 \int_{E(B)} |\nabla f|^2 \tilde{\varphi}\, d\nu,$$

which completes the proof of the theorem. □

We can also take $\varphi = 1$ in (4.17) and use the same argument to obtain a (strong) Poincaré inequality from the family of weak ones. The condition on $N_B$ is slightly weaker than in Theorem 4.8, since we have $P_j \leq cR_0^{2d}$ in (4.21).

LEMMA 4.9. *Let $B = B(x_0, R)$ be very good and let $R_0 = N_B \leq R^{1/(d+2)}$. Then*

$$\min_a \int_B (f(x) - a)^2\, d\mu \leq c_1 R^2 \int_{E(B)} |\nabla f|^2\, d\nu.$$

REMARK 4.10. The weight function $\varphi$ in (4.17) is similar to that in Saloff-Coste and Stroock (1991). Fabes and Stroock (1986) and Stroock and Zheng (1997) used weight functions which are supported on the whole space. In particular, Stroock and Zheng (1997) used

$$\psi_R(x) = \exp(-d(x_0, x)/R), \qquad x \in \mathbb{Z}^d,$$

and proved a weighted Poincaré inequality of the form

(4.25) $$\min_a \int_{\mathbb{Z}^d} (f - a)^2 \psi_R\, d\mu \leq c_1(d) R^2 \int_{\mathbb{E}_d} |\nabla f|^2 \widetilde{\psi}_R\, d\nu.$$

It is interesting to note that this fails for percolation clusters when $d \geq 3$.

To see this, fix a point $x_0 \in \mathcal{C}_\infty$ and $R \geq 1$ large enough so that $B(x_0, R)$ is good. If we look far enough from $x_0$ we can, $\mathbb{P}_p$-a.s., find a cube $Q$ of side $R$ with $Q \subseteq \mathcal{C}_\infty$ and such that $Q$ is only connected to the rest of $\mathcal{C}_\infty$ by one edge $\{x_1, x_2\}$. We take $x_1 \in Q^c$, $x_2 \in Q$ and write $s = d(x_0, x_2)$; we can assume $s \gg R$. We have $e^{-(s+dR)/R} \leq \psi \leq e^{-s/R}$ on $Q$.

Let $f = \mathbb{1}_Q$. Then as $\int \psi_R \asymp R^d$,

$$\bar{f}_R = \frac{\int_Q \psi_R}{\int \psi_R} \leq c e^{-s/R} \leq \frac{1}{4}$$

and

$$\min_a \int_{\mathbb{Z}^d} (f-a)^2 \psi_R \, d\mu = \int (f - \bar{f}_R)^2 \psi_R \, d\mu \geq \tfrac{1}{2} \int_Q \psi_R \, d\mu \geq \tfrac{1}{2} e^{-d} R^d e^{-s/R}.$$

On the other hand,

$$R^2 \int |\nabla f|^2 \widetilde{\psi} \, d\nu = R^2 \psi(x_2) = R^2 e^{-s/R}.$$

Thus (4.25) cannot hold with a constant $c_1$ independent of $R$.

**5. Lower bounds.** In this section we use the weighted Poincaré inequality and the technique of Fabes and Stroock (1986) to prove lower bounds for $q_t(x, y)$. We continue with the notation and assumptions of Section 3.

PROPOSITION 5.1. *Let $x_0 \in G$ and let $B_1 = B(x_0, R \log R)$ be a very good ball with $N_B \leq R^{1/(d+2)}$. Then if $x_1 \in B(x_0, \tfrac{1}{2} R \log R)$,*

(5.1) $\quad q_t(x_1, x_2) \geq c_1 t^{-d/2} \quad \text{for } x_2 \in B(x_1, R) \text{ and } \tfrac{1}{8} R^2 \leq t \leq R^2.$

PROOF. Let $x_3$ be such that $x_1, x_2 \in B(x_3, R/2)$. Write $R_1 = R \log R$, let $\rho = R/6$ and let $T = c_2 R^2$. Let $x_4 \in B(x_3, R/2)$. We apply Proposition 3.7 to $B_1$ with $t = T$. Since $T'_{B_1} = c R_1^2 \log R_1 \geq c' R^2 \log R$, the third condition in (3.14) holds, while the other two are evident. So if $t \leq T$, then

(5.2)
$$\sum_{x \in B(x_3, 2R/3)^c} q_t(x_4, x) \mu(x) = P^{x_4}(Y_t \notin B(x_1, 2R/3))$$
$$\leq P^{x_4}(\tau(x_4, R/6) < t)$$
$$\leq P^{x_4}(\tau(x_4, R/6) < T)$$
$$\leq c \exp(c' R^2 / T) \leq \tfrac{1}{2},$$

provided $c_2$ is chosen small enough. We can assume that $c_2 \leq \tfrac{1}{8}$.

Let $B = B(x_3, R)$, set $\rho(x) = d(x, B^c)$ for $x \in B$, and set

$$\varphi(y) = \left( \frac{R \wedge \rho(y)}{R} \right)^2, \quad y \in G \text{ and } V_0 = \sum_{x \in B} \varphi(x) \mu(x).$$

Then we have

(5.3) $\quad\quad\quad\quad c_3 R^d \leq V_0 \leq \mu(B) \leq c_4 R^d.$

Write $u(s, x) = u_s(x) = q_s(x_4, x)$, $s \geq 0$ and $x \in G$. Set

$$w(s, y) = w_t(x) = \log V_0 u(s, y), \quad H(t) = H(x_4, t) = V_0^{-1} \int_B \varphi w_t \, d\mu.$$

Then

$$V_0 H'(t) = \int_B \varphi \frac{\partial w_t}{\partial t} \, d\mu = \int_B \varphi \frac{1}{u_t} \frac{\partial u_t}{\partial t} \, d\mu = \sum_{x \in B} \frac{\varphi(x)}{u_t(x)} \mathcal{L} u_t(x) \mu(x).$$

Hence, writing $f = \mathbb{1}_B \varphi / u_t$, we have

$$V_0 H'(t) = -\tfrac{1}{2} \sum_{x \in G} \sum_{y \in G} a_{xy}(f(x) - f(y))(u_t(x) - u_t(u)).$$

Now we use the elementary inequality [see Stroock and Zheng (1997)]

$$-\left(\frac{d}{b} - \frac{c}{a}\right)(b - a) \geq \frac{1}{2}(c \wedge d)(\log b - \log a)^2 - \frac{(d-c)^2}{2(c \wedge d)},$$

which holds for any strictly positive $a, b, c, d$. If $f(x) > 0$ and $f(y) > 0$, then $x, y \in B$ and

$$-(f(x) - f(y))(u_t(x) - u_t(u))$$
$$= -\left(\frac{\varphi(x)}{u_t(x)} - \frac{\varphi(y)}{u_t(y)}\right)(u_t(x) - u_t(u))$$
$$\geq \frac{1}{2}(\varphi(x) \wedge \varphi(y))(\log u_t(x) - \log u_t(y))^2 - \frac{(\varphi(x) - \varphi(y))^2}{2(\varphi(x) \wedge \varphi(y))}.$$

If both $x \in B^c$ and $y \in B^c$, then $f(x) = f(y) = 0$, while if $x \in B$ and $y \in B^c$, then

$$-(f(x) - f(y))(u_t(x) - u_t(u)) = -\varphi(x)\left(1 - \frac{u_t(y)}{u_t(x)}\right).$$

We therefore have

$$(5.4) \quad V_0 H'(t) \geq \frac{1}{4} \sum_{x \in B} \sum_{y \in B} a_{xy}(\varphi(x) \wedge \varphi(y))(w_t(x) - w_t(y))^2$$

$$(5.5) \quad -\frac{1}{4} \sum_{x \in B} \sum_{y \in B} a_{xy} \frac{(\varphi(x) - \varphi(y))^2}{\varphi(x) \wedge \varphi(y)}$$

$$(5.6) \quad -\frac{1}{2} \sum_{x \in B} \sum_{y \in B^c} a_{xy} \varphi(x)\left(1 - \frac{u_t(y)}{u_t(x)}\right).$$

The sums in (5.4), (5.5) and (5.6) are called $S_1$, $S_2$ and $S_3$, respectively. To bound $S_2$ note that if $x \sim y$ with $x, y \in B$ and $k = \rho(x) \wedge \rho(y)$, then $k \geq 1$ and

$$\frac{(\varphi(x) - \varphi(y))^2}{\varphi(x) \wedge \varphi(y)} = R^{-2} \frac{(2k+1)^2}{k^2} \leq \frac{9}{R^2}.$$

So
$$S_2 \geq -\tfrac{9}{4}R^{-2}\sum_{x\in B}\mu(x) \geq -\tfrac{9}{4}R^{-2}\mu(B) \geq -cR^{-2}V_0.$$

Also, if $x \in \partial_i(B)$, then $\varphi(x) = R^{-2}$, so that
$$S_3 \geq -\sum_{x\in B}\sum_{y\in B^c} a_{xy}R^{-2} \geq -R^{-2}\sum_{x\in B}\mu(x) \geq -cR^{-2}V_0.$$

So the terms $S_2$ and $S_3$ are controlled by bounds of the same size, and we deduce, using Theorem 4.8,

$$\begin{aligned}
H'(t) &\geq \tfrac{1}{4}V_0^{-1}\sum_{x\in B}\sum_{y\in B} a_{xy}(\varphi(x)\wedge\varphi(y))(w_t(x)-w_t(y))^2 - c_5 R^{-2} \\
&\geq -c_5 R^{-2} + c_6 R^{-2}V_0^{-1}\sum_{x\in B}(w_t(x)-H(t))^2\varphi(x)\mu(x).
\end{aligned} \tag{5.7}$$

Let $I = [\tfrac{1}{2}T, T]$. We use only (5.7) for $t \in I$. Note that by Theorem 3.8 we have $V_0 u_t(x) \leq c_7$ for $t \in I$. Since $v \to (\log v - h)^2/v$ is decreasing on $[e^{2+h}, \infty)$, we have

$$\begin{aligned}
&\sum_{x\in B}(w_t(x)-H(t))^2\varphi(x)\mu(x) \\
&= \sum_{x\in B}\frac{(\log V_0 u_t(x)-H(t))^2}{u_t(x)}\varphi(x)u_t(x)\mu(x) \\
&\geq \frac{(\log c_7 - H(t))^2}{c_7}\sum_{x\in B:\, V_0 u_t(x)>e^{2+H(t)}}\varphi(x)V_0 u_t(x)\mu(x).
\end{aligned} \tag{5.8}$$

Then since $\varphi(x) \geq \tfrac{1}{9}$ on $B(x_1, 2R/3)$,

$$\begin{aligned}
&\sum_{x\in B:\, V_0 u_t(x)>e^{2+H(t)}} \varphi(x)u_t(x)\mu(x) \\
&= \sum_{x\in B}\varphi(x)u_t(x)\mu(x) - \sum_{x\in B:\, V_0 u_t(x)\leq e^{2+H(t)}}\varphi(x)u_t(x)\mu(x) \\
&\geq \tfrac{1}{9}\sum_{x\in B(x_1,2R/3)} u_t(x)\mu(x) - \sum_{x\in B}\varphi(x)V_0^{-1}e^{2+H(t)}\mu(x) \\
&\geq \tfrac{1}{9}\Big(1 - \sum_{x\in B(x_1,2R/3)^c} u_t(x)\mu(x)\Big) - e^{2+H(t)} \geq \tfrac{1}{18} - e^{2+H(t)},
\end{aligned} \tag{5.9}$$

where we used (5.2) in the final line.

Combining the estimates (5.7), (5.8) and (5.9) we deduce that
$$TH'(t) \geq -c_5 + c_8(\log c_7 - H(t))^2(\tfrac{1}{18} - e^{2+H(t)}).$$

Since $(a-h)^2 \geq \tfrac{1}{2}h^2$ if $h < -a$, this implies that there exists $c_9$ such that

$$TH'(t) \geq c_{10}H(t)^2 \text{ provided } H(t) < -c_9, \qquad t \in I. \tag{5.10}$$

If $\sup_{t \in I} H(t) < -c_9$, then (5.10) implies that $H(T) \geq -c_{11}$, while since $H(t) + c_5 T^{-1} t$ is increasing, if $\sup_{t \in I} H(t) \geq -c_9$, then $H(T) \geq -c_9 - c_2 c_5$. We therefore deduce that $H(T) \geq -c_{12}$ and it follows that

$$(5.11) \qquad H(t) = H(x_4, t) \geq -c_{13}, \qquad T \leq t \leq R^2.$$

To conclude the argument, we have, for $x_1, x_2 \in B(x_3, R/2)$, $t \in [T, R^2]$,

$$V_0 q_{2t}(x_1, x_2) = V_0^{-1} \sum_y V_0 q_t(x_1, y) V_0 q_t(x_2, y) \mu(dy)$$

$$\geq V_0^{-1} \sum_{y \in B} V_0 q_t(x_1, y) V_0 q_t(x_2, y) \varphi(y) \mu(dy).$$

So, using Jensen's inequality,

$$\log(V_0 q_{2t}(x_1, x_2)) \geq V_0^{-1} \sum_{y \in B} (\log(V_0 q_t(x_1, y)) + \log(V_0 q_t(x_2, y))) \varphi(y) \mu(dy)$$

$$= H(x_1, t) + H(x_2, t) \geq -2c_{13}.$$

Using (5.3) completes the proof of (5.1). □

LEMMA 5.2. *Let $x, y \in G$. Suppose there exist $r \geq 1$ and a path $x = z_0, \ldots, z_m = y$ such that for each $i = 0, \ldots, m$, $B(z_i, r \log r)$ is very good with $N_{B(z_i, r \log r)}^{d+2} \leq r$. Then*

$$(5.12) \qquad q_{mr}(x, y) \geq c(mr)^{-d/2} \exp(-c_1 m / r).$$

PROOF. This uses a chaining argument similar to that in Proposition 3.7. Choose points $x = w_0, w_1, \ldots, w_k = y$ along the path $\{z_0, \ldots, z_m\}$, such that $d(w_{i-1}, w_i) < r/3$ and $3m/r \leq k \leq 4m/r$. Let $s = mr/k$, so that $\frac{1}{4} r^2 \leq s \leq \frac{1}{3} r^2$. Let $B_j = B(w_j, r/3)$. We have $d(x', y') < r$ whenever $x' \in B_{j-1}$ and $y' \in B_j$. So by Proposition 5.1,

$$q_s(x', y') \geq c_2 s^{-d/2} \geq c_3 \mu(B_j)^{-1}, \qquad x' \in B_{j-1}, y' \in B_j.$$

So $P^{x'}(Y_s \in B_j) \geq c_3$ and therefore

$$P^x(Y_{ks} = y) \geq P^x(Y_{js} \in B_j, 1 \leq j \leq k-1, Y_{ks} = y) \geq \left( \prod_{j=1}^{k-1} c_3 \right) c_3 s^{-d/2}. \qquad \square$$

THEOREM 5.3. *Let $B = B(x_0, R \log R)$ be a very good ball with $N_B \leq R^{1/(d+2)}$. Let $d(x_0, x_1) \leq \frac{1}{2} R \log R$ and $x, y \in B(x_1, R)$. Then*

$$(5.13) \qquad q_t(x, y) \geq c_1 t^{-d/2} \exp(-c_2 d(x,y)^2 / t),$$

*provided*

$$(5.14) \qquad N_B^{2(2+d)} \leq t \leq R^2$$

*and*

(5.15) $$N_B^{2+d} d(x,y) \leq t.$$

PROOF. Let $D = d(x,y)$. If $D^2/t \leq 1$, then set $r = t^{1/2}$ and apply Proposition 5.1 to $B_1 = B(x, r \log r)$. Since $D \leq r \leq R$, $B_1 \subseteq B$ and $N_{B_1} \leq N_B$. So $N_{B_1}^{d+2} \leq N_B^{d+2} \leq t^{1/2} = r$, the hypotheses of Proposition 5.1 hold and we deduce that $q_t(x,y) \geq c t^{-d/2}$.

For the case $D^2/t > 1$ we use Lemma 5.2. Let $x = z_0, \ldots, z_D = y$ be a shortest path between $x$ and $y$, and let $r = t/D$. By (5.15), $r \geq N_B^{d+2}$, while $r \leq t^{1/2} \leq R$, so that $B_i = B(z_i, r \log r) \subseteq B$ and hence $N_{B_i}^{d+2} \leq N_B^{d+2} \leq r$. So the hypotheses of Lemma 5.2 are satisfied and we obtain (5.13). □

REMARK. The restriction in (5.15) is weaker than the hypotheses in the upper bound of Theorem 3.8, where we were able to use global upper bounds on $q_t$ to restrict to the case when $t$ was close to $D^{1/2}$. The lower bound argument for $cD \leq t \leq cDN_B^{2+d}$ requires the existence of a chain of small balls (of size roughly $r = t/D$) on which the lower bounds of Proposition 5.1 are valid. If $r = O(1)$ so that $t \simeq D$, then we can just use a path in the graph and (1.4) to deduce that $q_t(x,y) \geq e^{-ct} \simeq \exp(-c' D^2/t)$. However, if $D = t^{1-\varepsilon}$, then we need $r \simeq t^{\varepsilon} \gg 1$ and elementary bounds are not enough.

For a very good ball $B$, the volume condition or the Poincaré inequality fails for some subballs of size $r < N_B$. However, these may hold for enough small subballs so that we can still use Lemma 5.2, and in Section 2 it is proved, in the percolation context, that it is possible to find such a chain. Fix $C_E \geq 1$ and $C_F > 1$.

DEFINITION 5.4. A ball $B = B(x_0, R_1)$ is *exceedingly good* if:

1. We have $B$ is very good with $N_B^{10(d+2)} \leq R_1$.
2. For each $x_1, x_2 \in B(x_0, R_1)$ with $d(x_1, x_2) \geq R_1^{1/4}$ and $C_E \leq r \leq N_B^{2+d}$ there exists a path $y_1 = z_0, \ldots, z_k = y_2$ with the following properties:

    (a) $B_i = B(z_i, r \log r)$ is very good with $N_{B_i}^{2+d} \leq r$;
    (b) $k \leq C_F d(x,y)$;
    (c) $d(x_j, y_j) \leq R_1^{1/4}$, $j = 1, 2$.

REMARK 5.5. If $B$ is very good and $N_B^{2+d} \leq r \leq R/\log R$, then taking $m = d(x,y)$ and $z_0, \ldots, z_m$ to be a shortest path between $x_1$ and $x_2$, we get a path satisfying 2(a)–(c) above, with $y_j = x_j$.

THEOREM 5.6. *Let $B = B(x_0, R \log R)$ be exceedingly good and let $x_1, x_2 \in B(x_0, R)$. Then there exist constants $c_i$ (depending on $C_E$ and $C_F$) such that*

(5.16)
$$q_t(x_1, x_2) \geq c_1 t^{-d/2} \exp(-c_2 d(x_1, x_2)^2/t),$$
$$R^{1/2} \vee d(x, y) \leq t \leq R^2.$$

PROOF. Let $R_1 = R \log R$ and $D = d(x_1, x_2)$. By Theorem 5.3 it is enough to consider the case when $t$ satisfies the condition in (5.16) but fails (5.14) or (5.15). Since $t \geq R^{1/2} \geq R_1^{1/3} \geq N_B^{3(d+2)}$, (5.14) must hold. So we just need to consider the case when $D \leq t \leq N_B^{2+d} D$. Note that this implies that

(5.17)
$$D \geq \frac{D^2}{t} \geq \frac{t}{N_B^{2(d+2)}} \geq \frac{R^{1/2}}{R_1^{1/5}} \geq R_1^{1/4}.$$

In particular, $D^2/t \gg \log t$ so that, for this range of $t$ and $D$, terms of the order $t^{-d/2}$ can be absorbed into the constant $c_2$ in (5.16).

Let $r = t/(2C_F D)$, so that $(2C_F)^{-1} \leq r \leq (2C_F)^{-1} N_B^{2+d}$. We have to consider two cases. If $1 \leq r \leq C_E$, then we use (1.4) directly. Set $s = t/D$, so that $c_3 \leq s \leq c_4$, and we obtain $q_{Ds} \geq e^{-cD} \geq \exp(-cD^2/t)$.

If $C_E \leq r$, then we use the fact that $B$ is exceedingly good so that there exists a path $z_0, \ldots, z_k$ that satisfies conditions 2(a)–(c) of Definition 5.4. We have $D \leq k \leq C_F D$, so that $(2C_F)^{-1} t \leq kr \leq \frac{1}{2} t$ and $k/r \leq 2C_F^2 D^2/t$. Applying Lemma 5.2,

(5.18) $\quad q_{kr}(z_0, z_k) \geq c(kr)^{-d/2} \exp(-c_5 k/r) \geq c t^{-d/2} \exp(-cD^2/t).$

By (5.17), $D^2/t \geq R_1^{1/4} \geq d(x_1, y_1) \vee d(x_2, y_2)$. Let $D_i = d(x_i, z_i)$. By (5.17), $D_i \leq R_1^{1/4} \leq t/8$. Using (1.4),

(5.19) $\qquad q_{m_j}(x_j, z_j) \geq \exp(-c_6 m_j) \geq \exp(-c_6 D^2/t).$

Let $u = D_1 + D_2 + kr$. Then $\frac{1}{2} t \leq u \leq \frac{3}{4} t$. By (1.4),

$$q_t(x_0, x_1) \geq q_u(x_0, x_1) q_{t-u}(x_1, x_1)$$
$$\geq q_{D_0}(x_0, z_0) q_{kr}(z_0, z_k) q_{D_1}(z_k, x_1) q_{t-u}(x_1, x_1).$$

Using the bounds (5.18) and (5.19), and Proposition 5.1 to control the final term, we obtain (5.16). □

THEOREM 5.7. *Let $d \geq 2$, and let $C_0$, $C_V$, $C_P$ and $C_W$ be constants. Let $\mathcal{G} = (G, E)$ be an infinite graph that satisfies (1.1) and let $x \in G$.*

(a) *Suppose that there exists $R_0 = R_0(x)$ such that $B(x, R)$ is very good with $N_{B(x,R)}^{3(d+2)} \leq R$ for each $R \geq R_0$. There exist constants $c_i = c_i(d, C_0, C_V, C_P, C_W)$ such that if $t$ satisfies*

(5.20)
$$t \geq R_0^{2/3},$$

*then*

(5.21) $\quad q_t(x,y) \leq c_1 t^{-d/2} \exp(-c_2 d(x,y)^2/t), \qquad d(x,y) \leq t,$

*and*

(5.22) $\quad q_t(x,y) \geq c_3 t^{-d/2} \exp(-c_4 d(x,y)^2/t), \qquad d(x,y)^{3/2} \leq t.$

(b) *If in addition there exists $R_E = R_E(x)$ such that $B(x,R)$ is exceedingly good for each $R \geq R_E$, then the lower bound (5.22) holds (with constants depending in addition on $C_E$ and $C_F$) whenever $t \geq R_E$ and $t \geq d(x,y)$.*

PROOF. Fix $y \in G$ and let $D = d(x,y)$. We need to check that we can find $R$ so that we can apply Theorems 3.8, 5.3 and 5.6. We begin with (a). Let $R = t^{2/3}$, so that $R \geq D$; by (5.20) we have $R \geq R_0$. Let $B = B(x,R)$; then $N_B^{2(2+d)} \leq R$.

To apply Theorem 3.8 and obtain (5.21) we need (3.18) to hold, but this is clear since $t \geq D^{3/2} \geq D$ and

$$N_B^{2d+1} \leq N_B^{2(d+2)} \leq R = t^{2/3} \leq t = R^{3/2} \leq \frac{cR^2}{\log R}.$$

To obtain (5.22), we use Theorem 5.3. Since $R \geq D$ we have $y \in B(x,R)$. Let $R_1 = R \log R$. Since $R_1 \geq R_0$, $B_1 = B(x,R_1)$ is very good with $N_{B_1}^{3(2+d)} \leq R_1$. So

$$N_{B_1}^{2(2+d)} \leq R \log R \leq R^{3/2} = t \leq \tfrac{1}{4} R^2$$

and (5.14) is verified. Also,

$$cDN_{B_1}^{2+d} \leq c' t^{2/3} R_1^{1/3} \leq t^{2/3} R^{1/2} \leq t,$$

so (5.15) holds and we obtain (5.22).

(b) We now take $R = t$, so that $R \geq R_E$. Then $D \leq t = R$ and we apply Theorem 5.6 to $B(x,R_1)$ with $R_1 = R \log R$. Condition (5.16) is easily verified and the bound in (5.22) is immediate. $\square$

We now prove an elliptic Harnack inequality for graphs satisfying the conditions of Theorem 5.3. The approach is the same as in Section 5 of Fabes and Stroock (1986), but we need an additional argument in Theorem 5.11.

LEMMA 5.8. *Let $B = B(x_0, R \log R)$ be a very good ball with $N_B^{2(d+2)} \leq R$. Let $d(x_0, x_1) \leq \tfrac{1}{2} R \log R$, let $B_0 = B(x_1, R)$ and let $q_t^0(x,y)$ be the density of the process $Y$ killed at the exit time $\tau$ of $Y$ from $B_0$. Then*

(5.23) $\quad q_t^0(x,y) \geq c_1 t^{-d/2}, \qquad x,y \in B(x_1, 3R/4), \qquad c_2 R^2 \leq t \leq R^2.$

PROOF. Whereas this argument is quite standard [see, e.g., Lemma 5.1 of Fabes and Stroock (1986)], we only give a sketch. For $x, y \in B_0$ we have

(5.24) $$\begin{aligned} q_t^0(x,y) &\geq q_t(x,y) - E^x \mathbb{1}_{(\tau < t)} q_{t-\tau}(Y_\tau, y) \\ &\geq q_t(x,y) - \sup_{0 \leq s \leq t} \sup_{z \in \partial B_0} q_s(z,y). \end{aligned}$$

Let $\delta \in (0, \frac{1}{8})$, $d(x,y) \leq \delta R$ and $t = \delta^2 R^2$. Using the estimates in Theorems 5.3 and 3.8 in (5.24) we obtain, provided $t$ and $s = \theta t$ satisfy conditions (3.18), (5.14) and (5.15),

$$q_t^0(x,y) \geq t^{-d/2} \bigg[ c_1 \exp(-c_2) - \sup_{0 \leq \theta \leq 1} c_3 \theta^{-d/2} \exp(-c_4/(\theta \delta^2)) \bigg].$$

[If $s$ is too small to satisfy (3.18), we use Lemma 1.1.] Hence if $\delta$ is chosen small enough, we obtain

(5.25) $$q_t^0(x,y) \geq c_5 t^{-d/2}.$$

A chaining argument now gives (5.23). □

DEFINITION. Write $\overline{B}(x,R) = B(x,R) \cup \partial_e(B(x,R))$. A function $h: \overline{B}(x,R) \to \mathbb{R}$ is *harmonic* on $B(x,R)$ if

$$\mathcal{L}h(x) = 0, \qquad x \in B(x,R).$$

We write $\text{Osc}(h,A) = \sup_A h - \inf_A h$.

The following oscillation bound follows from (5.23) just as in Lemma 5.2 of Fabes and Stroock (1986).

LEMMA 5.9. *Let $B = B(x_0, R \log R)$ be a very good ball with $N_B^{2(d+2)} \leq R$. Let $d(x_0, x_1) \leq \frac{1}{2} R \log R$, $B_0 = B(x_1, R)$, $B_1 = B(x_1, \frac{1}{2}R)$ and $h$ be harmonic in $B_0$. There exists $c_1 \in (0,1)$ such that*

$$\text{Osc}(h, B_1) \leq (1 - c_1) \text{Osc}(h, B_0).$$

PROOF. By a linear transformation we can assume $\min_{B_0} h = 0$, $\max_{B_0} h = 1$ and $\int_{B_1} h \, d\mu \geq \frac{1}{2} \mu(B_1)$. Then if $x \in B_1$, by Lemma 5.8,

$$h(x) \geq \int_{B_1} q_{R^2}^0(x,y) h(y) \mu(dy) \geq \tfrac{1}{2} c_2 R^{-d} \mu(B_1) \geq c_3.$$

So $\text{Osc}(h, B_1) \leq (1 - c_3)$. □

We also need an intermediate range Harnack inequality.

LEMMA 5.10. *Let $B = B(x_0, R \log R)$ be a very good ball with $N_B^{3(d+2)} \leq R$. Let $d(x_0, x_1) \leq \frac{1}{2} R \log R$, $B_0 = B(x_1, R)$, and $h$ be nonnegative and harmonic in $B_0$. Then if $d(x_1, y) \leq R/2$ and $r = R^{1/2}$,*

$$\sup_{B(y,r)} h \leq c_1 \inf_{B(y,r)} h.$$

PROOF. Let $B_1 = B(x_1, 3R/4)$. We can assume $\inf_{B_1} h = 1$. The local Harnack inequality for graphs implies that if $x \sim y$ and $x, y \in B_1$, then $h(x) \leq C_0 h(y)$. Hence $h(x) \leq \exp(c_2 R)$, $x \in \overline{B}_0$. We extend $h$ to $G$ by taking $h \equiv 1$ outside $\overline{B}_0$. Let $\tau$ be the first exit time of $Y$ from $B_1$. Whereas $h(Y_t)$ is a martingale on $[0, \tau_0]$, for $x \in B(x_1, R/2)$,

$$\begin{aligned} h(x) &= E^x h(Y_{t \wedge \tau}) \\ &= E^x \mathbb{1}_{B_1}(Y_t) h(Y_{t \wedge \tau}) + E^x \mathbb{1}_{B_1^c}(Y_t) h(Y_\tau) \\ &= E^x \mathbb{1}_{B_1}(Y_t) h(Y_t) + E^x \mathbb{1}_{B_1}(Y_t) \mathbb{1}_{(\tau < t)}(h(Y_\tau - h(Y_t))) + E^x \mathbb{1}_{B_1^c}(Y_t) h(Y_\tau) \\ &\leq \int_{B_1} q_t(x, y) h(y) \, d\mu(y) + \exp(c_2 R) P^x(\tau < t). \end{aligned}$$

Using Proposition 3.7, the final term above is bounded by $c_3 \exp(c_2 R - c_4 R^2/t) \leq c_3$ if $t \leq Rc_2/c_4$.

Now let $d(z, x) \leq R^{1/2}$ and $s = \lambda t$ with $\lambda > 1$. Then if $\lambda \theta \leq c_2/c_4$,

$$\begin{aligned} h(z) &\geq E^z \mathbb{1}_{B_1}(Y_s) h(Y_s) - E^z \mathbb{1}_{(\tau < t)} h(Y_s) \\ &\geq \int_{B_1} q_s(z, y) h(y) \, d\mu(y) - \exp(c_2 R) P^x(\tau < s) \\ &\geq \int_{B_1} q_s(z, y) h(y) \, d\mu(y) - c_3. \end{aligned}$$

Also, if $u = R^2$, by Lemma 5.8,

$$h(z) \geq \int_{B_1} q_t^0(z, y) h(y) \, d\mu(y) \geq \int_{B_1} c_5 u^{-d/2} h(y) \, d\mu(y).$$

Hence

$$2h(z) \geq \int_{B_1} (c_5 R^{-d} + q_s(z, y)) h(y) \, d\mu(y) - c_3.$$

Using Lemma 1.1, and Theorems 3.8 and 5.3, we can choose $\lambda$ so that

$$c_5 R^{-d} + q_s(z, y) \geq c_6 q_t(x, y), \qquad y \in B_1.$$

It follows that $2h(z) \geq c_6(h(x) - c_3) - c_3$, so that [as $h(z) \geq 1$] we have $c_6 h(x) \leq (2 + c_3(1 + c_6)) h(z)$. □

THEOREM 5.11. *Let $B = B(x_0, R \log R)$ be very good, with $N_B^{4(d+2)} \leq R$. Then if $d(x_0, x_1) \leq \frac{1}{3} R \log R$, $B_0 = B(x_1, R)$, and $h: \overline{B}_0 \to \mathbb{R}$ is nonnegative and harmonic in $B_0$,*

$$\sup_{B(x_1, R/2)} h \leq c_1 \inf_{B(x_1, R/2)} h.$$

PROOF. We begin in the same way as in Fabes and Stroock (1986). Let $\alpha_3 = 1/4$. From Lemma 5.9 it follows that there exists $A$ such that

$$\text{Osc}(h, B(y, r)) \leq (2e^d)^{-1} \text{Osc}(h, B(y, Ar)) \quad \text{if } R^{\alpha_3} \leq r, Ar \leq R/8.$$
(5.26)

We normalize $h$ so that $\min_{B(x_1, R/4)} h = 1$. Let $x_2 \in B(x_1, R/4)$ satisfy $h(x_2) = 1$. By Lemma 5.8, if $y \in B(x_1, 3R/4)$ and $B(y, s)$ is good,

$$1 = h(x_2) \geq c_2 \int_{B(y,s)} R^{-d} h(y) \, d\mu(y) \geq c_3 (s/R)^d \inf_{B(y,s)} h.$$

Thus

(5.27) $$\inf_{B(y,s)} h \leq c_3^{-1} \frac{R^d}{s^d}.$$

Now let $M_r = 2c_3^{-1} e^{rd}$ and $a_r = Re^{-r}$.

Suppose that there exists $y_r \in B(x_1, R/2)$ with $h(y_r) \geq M_r$. Then, by (5.27), $\text{Osc}(h, B(y_r, a_r)) \geq c_3^{-1} e^{rd}$. So, provided

(5.28) $$a_r \geq R^{\alpha_3} \quad \text{and} \quad Aa_r \leq R/8,$$

it follows that $\text{Osc}(h, B(y_r, Aa_r)) \geq 2e^d c_3^{-1} e^{rd} = M_{r+1}$. Hence there exists $y_{r+1}$ with $d(y_r, y_{r+1}) \leq Aa_r$ with $h(y_{r+1}) \geq M_{r+1}$.

Choose $r_0$ so that $\sum_{r_0}^{\infty} Aa_r \leq R/8$. Then if $\sup_{B(x_1, R/4)} h \geq M_{r_0}$, the argument above implies that we can construct a sequence $y_j$, $r_0 \leq j \leq k$ with $y_j \in B(x_1, R/2)$ and $h(y_j) \geq M_j$. Here $k$ is the largest $r$ such that (5.28) holds. We have $M_k = 2c_3^{-1} e^{kd} = cR^d a_k^{-d} \geq c_4 R^{d(1-\alpha_3)}$. The local Harnack inequality in Lemma 5.10 implies that

(5.29) $$\inf_{B(y_k, R^{1/2})} h \geq c_5 c_4 R^{d(1-\alpha_3)},$$

while (5.27) implies that

(5.30) $$\inf_{B(y_k, R^{1/2})} h \leq c_3^{-1} R^{d/2}.$$

Since $\alpha_3 < \frac{1}{2}$ this gives a contradiction if $R \geq c_6$. So we deduce that

$$\sup_{B(x_1, R/4)} h \leq M_{r_0} \inf_{B(x_1, R/4)} h.$$

The Harnack inequality for $B(x_1, R/2) \subseteq B(x_1, R)$ now follows by an easy chaining argument. □

**6. Random walks and percolation.** In this section we tie together the results of Section 2 and Sections 3–5 and prove Theorems 1–5. We recall the notation for bond percolation: We take $\Omega = \{0,1\}^{\mathbb{E}_d}$, write $\eta_e, e \in \mathbb{E}_d$, for the coordinate maps and take $\mathbb{P}_p$ to be the probability on $\Omega$ which makes the $\eta_e$ i.i.d. Bernoulli r.v. with mean $p$. We assume that $p > p_c$ so that there exists $\Omega_0$ with $\mathbb{P}_p(\Omega_0) = 1$ such that for $\omega \in \Omega_0$ there is a unique infinite cluster $\mathcal{C}_\infty(\omega)$. As in the Introduction we take $Y$ to be the Markov process with generator $\mathcal{L}_\omega$ given by (0.1), and write $q_t^\omega(x,y)$ for the transition density of $Y$ as given by (0.2). Write $(P_\omega^x, x \in \mathcal{C}_\infty)$ for the probability law of $Y$.

PROPOSITION 6.1. *Let $p > p_c$. There exists $\Omega_1 \subseteq \Omega$ with $\mathbb{P}_p(\Omega_1) = 1$ and $S_x, x \in \mathbb{Z}^d$, such that $S_x(\omega) < \infty$ for each $\omega \in \Omega_1$, $x \in \mathcal{C}_\infty(\omega)$. There exist constants $c_i = c_i(d,p)$ such that for $x, y \in \mathcal{C}_\infty(\omega)$, $t \geq 1$ with*

$$(6.1) \qquad S_x(\omega) \vee d_\omega(x,y) \leq t,$$

*the transition density $q_t(x,y)$ of $Y$ satisfies*

$$(6.2) \quad c_1 t^{-d/2} \exp(-c_2 d_\omega(x,y)^2/t) \leq q_t^\omega(x,y) \leq c_3 t^{-d/2} \exp(-c_4 d_\omega(x,y)^2/t).$$

PROOF. Let the constants $C_0$, $C_V$, $C_P$, $C_W$, $C_E$ and $C_F$ (depending on $d$ and $p$) be as in Section 2 and, as in Section 2, let $\alpha_2^{-1} = 11(d+2)$. Let $(N_x, x \in \mathbb{Z}^d)$ be as in Lemma 2.24 and let $\Omega_1 = \{\omega : N_x(\omega) < \infty \text{ for all } x\}$. Let $\omega \in \Omega_1$ and $x \in \mathcal{C}_\infty(\omega)$. We set $S_x = R_E(x) = N_x$ and check the hypotheses of Theorem 5.7. Let $R \geq R_E$, $B = B_\omega(x, R) \subseteq Q = Q(x, R)$ and $n = 2R = s(Q)$.

Whereas $\omega \in L(Q) \subseteq D(Q, \alpha_2) \cap H(Q, \alpha_2)$, applying Theorem 2.18(c) we deduce that $B$ is very good with $N_B \leq C_F n^{\alpha_2} \leq R^{1/(10(d+2))}$.

Now let $C_E \leq r \leq N_B^{d+2}$ and let $x_1, x_2 \in B$ with $d_\omega(x_1, x_2) \geq R^{1/4}$. Whereas $D(Q, \alpha_2)$ holds, $d_\omega(x_1, x_2) \geq \frac{1}{3} n^{1/4}$. Choose $m$ so that $m/16 = \lfloor r \log r \rfloor$ and apply Theorem 2.18(b) to deduce that there exists a path $x_1' = z_0, \ldots, z_j = x_2'$ that satisfies condition (b) of Theorem 2.23. We need to verify (a)–(c) of Definition 5.4.2. Part (a) holds since $B_i = B_\omega(z_i, r \log r) = B_\omega(z_i, m/16)$ is very good and $N_{B_i} \leq m^{1/(d+4)} < r^{1/(d+2)}$. For (b) we have $j \leq c_{2.23.3}|x_1 - x_2|_\infty \leq dc_{2.23.3} d_\omega(x_1, x_2)$. Part (c) is easily verified as $d_\omega(x_i, x_i') \leq \frac{1}{3} n^{1/4} < R^{1/4}$. Thus $B$ is exceedingly good.

So we can apply Theorem 5.7 to deduce that the bound (6.2) holds for $t \geq R_E$ and $y$ such that $d(x,y) \leq t$. □

PROOF OF THEOREM 1. Given Proposition 6.1, all that remains is to replace the chemical distance $d_\omega(x,y)$ by $|x-y|_1$. Whereas $|x-y|_1 \leq d_\omega(x,y)$, the upper bound in (0.4) is immediate. For the lower bound, we take $S_x$ as in Proposition 6.1, and let $x, y \in \mathcal{C}_\infty$ and $t \geq S_x$ with $|x-y|_1 \leq t$. Choose the smallest cube $Q$ of side $n$ that contains $x$ and $y$ and with $n \geq S_x$. Since $Q$ is very good, by Proposition 2.17(d) we have $|x-y|_\infty \leq C_F(n^{\alpha_2} + 1) \vee |x-y|_\infty$.

If $|x-y|_\infty \geq 1 + n^{\alpha_2}$, then it follows that $d_\omega(x,y) \leq c|x-y|_\infty$, so the lower bound in (0.4) follows. If $|x-y|_\infty < 1 + n^{\alpha_2}$, then $d_\omega(x,y) \leq n^{\alpha_2}$. So because $t \geq S_x$, both $d_\omega(x,y)^2/t$ and $|x-y|_1^2/t$ are less than 1, and again the lower bound in (0.4) holds. □

PROOF OF THEOREM 2. Fix $x,y \in \mathbb{Z}^d$, let $D = |x-y|_1$ and fix $t \geq D$. Write $A = \{x, y \in \mathcal{C}_\infty\}$. Since $c_5 \leq \mathbb{P}_p(x, y \in \mathcal{C}_\infty) \leq 1$, it is enough to prove (0.6) for $\mathbb{E}_p(q_t^\omega(x,y)\mathbb{1}_A)$. Let $S_x$ be as in the proof of Theorem 1 and let $n = t/2$. Then

$$(6.3) \quad \mathbb{E}_p(q_t^\omega(x,y)\mathbb{1}_A) = \mathbb{E}_p(q_t^\omega(x,y)\mathbb{1}_A\mathbb{1}_{\{S_x < n\}}) + \mathbb{E}_p(q_t^\omega(x,y)\mathbb{1}_A\mathbb{1}_{\{S_x \geq n\}}).$$

By the proof of Theorem 1, if $\omega \in A$ and $S_x(\omega) \leq n$, then $q_t^\omega(x,y)$ satisfies the bounds in (0.4). So the lower bound in (0.6) is immediate since

$$(6.4) \quad \begin{aligned} \mathbb{E}_p(q_t^\omega(x,y)\mathbb{1}_A) &\geq \mathbb{E}_p(q_t^\omega(x,y)\mathbb{1}_A\mathbb{1}_{\{S_x < n\}}) \\ &\geq c_6 \exp(-c_7 D^2/t)(\mathbb{P}_p(A) - \mathbb{P}_p(A \cap \{S_x \geq n\})) \\ &\geq c_6 \exp(-c_7 D^2/t)(c_5 - c_8 \exp(-t^{\alpha\beta})). \end{aligned}$$

So if $t \geq c_9$, then the final term in (6.4) is greater than $\frac{1}{2}$ and we obtain (0.6). If $t < c_9$, then with probability at least $c_5 p^D$ there is a path of length $D$ in $\mathcal{C}_\infty$ that joins $x$ and $y$, and using (1.4), we deduce that on this event, $q_t^\omega(x,y) \geq c^D \geq c'$. So (0.6) holds in this case also.

To prove the upper bound, recall that by Lemma 1.1 we always have (when $t \geq D$) that

$$q_t^\omega(x,y) \leq c \exp(-2c_{10} d_\omega(x,y)^2/t) \leq c \exp(-2c_{10} D^2/t).$$

If $\exp(-c_{10} D^2/t) \leq t^{-d/2}$, then the upper bound in (0.6) is then immediate. If not, then by Theorem 1 we have

$$\mathbb{E}_p(q_t^\omega(x,y)\mathbb{1}_A\mathbb{1}_{\{S_x < n\}}) \leq c_{11} \exp(-c_{12} D^2/t),$$

so it remains to control the second term in (6.3). We have, provided $t \geq c_{13}$,

$$\mathbb{E}_p(q_t^\omega(x,y)\mathbb{1}_A\mathbb{1}_{\{S_x \geq n\}}) \leq \mathbb{P}_p(S_x \geq n) \leq c \exp(-ct^{\alpha_2 \beta})$$
$$\leq t^{-d} \leq t^{-d/2} \exp(-c_{10} D^2/t).$$

If $t \leq c_{13}$, then (as $D \leq t$) we obtain the upper bound in (0.6) by taking the constant $c_2$ large enough. □

PROOF OF THEOREM 3. The proof follows easily from Theorem 5.11 and Proposition 6.1. □

PROOF OF THEOREM 4. (a) This is a well-known consequence of the Harnack inequality. Let $h : \mathcal{C}_\infty \to (0, \infty)$ be a global harmonic function. Replacing $h$ by $ch$ if necessary, we can assume $h \geq 1$ everywhere and that there exists $x_0$ with $h(x_0) < 2$. Applying Theorem 3 to $B_\omega(x_0, R) \subseteq B_\omega(x_0, 2R) \subseteq \mathcal{C}_\infty$, when $R$ is large enough we deduce that $\sup_{B_\omega(x_0,R)} h \leq 2c_1$, so that $h$ is bounded.

Let $B_n = B(x_0, 2^n R)$. Then Theorem 3 implies the oscillation inequality $\inf_{B_n} h \geq c_1^{-1} \sup_{B_{n+1}} h$, so that $\mathrm{Osc}(h, B_n) \leq (1 - c_1^{-1}) \mathrm{Osc}(h, B_{n+1})$. By iterating we deduce that $\mathrm{Osc}(h, B_n) \geq (1 - c_1^{-1})^{-n} \mathrm{Osc}(h, B_0)$ and since $h$ is bounded, this implies that $\mathrm{Osc}(h, B_0) = 0$. Thus $h$ is constant on any large ball and so is constant. (b) This is also standard [Lemma 5.2 of Fabes and Stroock (1986)]. Lemma 5.8 gives lower bounds on the transition probability $q_t^{0,\omega}(x,y)$ for $Y$ killed outside a ball $B(x_0, R)$ for all sufficiently large $R$. Let $F$ be an event in the tail field and set $f(s,x) = P_\omega(F|Y_s = x)$. Then $0 \leq f \leq 1$ and $f$ satisfies

$$(6.5) \qquad f(s,x) = \int_{\mathcal{C}_\infty} q_{t-s}^\omega(x,y) f(t,y) \mu(dy), \qquad s < t.$$

Fix $x_0 \in \mathcal{C}$, $t_0 \geq 0$ and set

$$A(R) = B(x_0, \tfrac{1}{2}R) \times [t_0, t_0 + R^2].$$

Let $g(s,x)$ satisfy (6.5) with $\min_{A(2R)} g = 0$, $\max_{A(2R)} g = 1$ and

$$\int_{B(x_0, R/2)} g(t_0 + 4R^2, y) \mu(dy) \geq \tfrac{1}{2}.$$

Then if $(s,x) \in A(R)$,

$$g(s,x) \geq \int_{B(x_0,R/2)} q_{t_0+4R^2-s}^{0,\omega}(x,y) g(t_0 + 4R^2, y) \mu(dy)$$

$$\geq \int_{B(x_0,R/2)} c_2 R^{-d} g(t_0 + 4R^2, y) \mu(dy)$$

$$\geq \tfrac{1}{2} c_2 R^{-d} \mu(B(x_0, \tfrac{1}{2}R)) \geq c_3.$$

Hence there exists $\delta > 0$ such that

$$\mathrm{Osc}(f, A(R)) \leq (1 - \delta) \mathrm{Osc}(f, A(2R)),$$

and by iterating, it follows that $f$ is constant. $\square$

PROOF OF THEOREM 5. The proof is immediate on integrating the bounds (0.4) and using (1.6) to control $q_t(x,y)$ for $|x - y|_1 \geq t$. $\square$

## REFERENCES

ANTAL, P. and PISZTORA, A. (1996). On the chemical distance for supercritical Bernoulli percolation. *Ann. Probab.* **24** 1036–1048. MR1404543

BARLOW, M. T. and BASS, R. F. (1989). The construction of Brownian motion on the Sierpinski carpet. *Ann. Inst. H. Poincaré Probab. Statist.* **25** 225–257. MR1023950

Barlow, M. T. and Bass, R. F. (1992). Transition densities for Brownian motion on the Sierpinski carpet. *Probab. Theory Related Fields* **91** 307–330. MR1151799

Barlow, M. T. and Bass, R. F. (2004). Stability of parabolic Harnack inequalities. *Trans. Amer. Math. Soc.* **356** 1501–1533. MR2034316

Bass, R. F. (2002). On Aronsen's upper bounds for heat kernels. *Bull. London Math. Soc.* **34** 415–419. MR1897420

Benjamini, I. and Mossel, E. (2003). On the mixing time of simple random walk on the super critical percolation cluster. *Probab. Theory Related Fields* **125** 408–420. MR1967022

Benjamini, I., Lyons, R. and Schramm, O. (1999). Percolation perturbations in potential theory and random walks. In *Random Walks and Discrete Potential Theory* 56–84 (M. Dicardello and W. Woess, eds.). Cambridge Univ. Press. MR1802426

Carlen, E. A., Kusuoka, S. and Stroock, D. W. (1987). Upper bounds for symmetric Markov transition functions. *Ann. Inst. H. Poincaré Probab. Statist.* **23** 245–287. MR898496

Coulhon, T. and Grigor'yan, A. (2003). Pointwise estimates for transition probabilities of random walks on infinite graphs. In *Fractals* (P. Grabner and W. Woess, eds.) 119–134. Birkhäuser, Boston. MR2091701

Couronné, O. and Messikh, R. J. (2003). Surface order large deviations for 2D FK-percolation and Potts models. Preprint. MR2078538

Davies, E. B. (1993). Large deviations for heat kernels on graphs. *J. London Math. Soc. (2)* **47** 65–72. MR1200978

De Gennes, P. G. (1976). La percolation: Un concept unificateur. *La Recherche* **7** 919–927.

Delmotte, T. (1999). Parabolic Harnack inequality and estimates of Markov chains on graphs. *Rev. Mat. Iberoamericana* **15** 181–232. MR1681641

De Masi, A., Ferrari, P. A., Goldstein, S. and Wick, W. D. (1989). An invariance principle for reversible Markov processes. Applications to random motions in random environments. *J. Statist. Phys.* **55** 787–855. MR1003538

Deuschel, J.-D. and Pisztora, A. (1996). Surface order large deviations for high-density percolation. *Probab. Theory Related Fields* **104** 467–482. MR1384041

Fabes, E. B. and Stroock, D. W. (1986). A new proof of Moser's parabolic Harnack inequality via the old ideas of Nash. *Arch. Rational. Mech. Anal.* **96** 327–338. MR855753

Grigor'yan, A. A. (1992). Heat equation on a noncompact Riemannian manifold. *Math. USSR-Sb.* **72** 47–77. MR1098839

Grimmett, G. R. (1999). *Percolation*, 2nd ed. Springer, Berlin. MR1707339

Grimmett, G. R., Kesten, H. and Zhang, Y. (1993). Random walk on the infinite cluster of the percolation model. *Probab. Theory Related Fields* **96** 33–44. MR1222363

Heicklen, D. and Hoffman, C. (2000). Return probabilities of a simple random walk on percolation clusters. Preprint.

Jerison, D. (1986). The weighted Poincaré inequality for vector fields satisfying Hörmander's condition. *Duke Math. J.* **53** 503–523. MR850547

Kaimanovitch, V. A. (1990). Boundary theory and entropy of random walks in random environment. In *Probability Theory and Mathematical Statistics* 573–579. Mokslas, Vilnius.

Kesten, H. (1986a). The incipient infinite cluster in two-dimensional percolation. *Probab. Theory Related Fields* **73** 369–394. MR859839

Kesten, H. (1986b). Subdiffusive behavior of random walks on a random cluster. *Ann. Inst. H. Poincaré Probab. Statist.* **22** 425–487. MR871905

Kusuoka, S. and Zhou, X. Y. (1992). Dirichlet form on fractals: Poincaré constant and resistance. *Probab. Theory Related Fields* **93** 169–196. MR1176724

Liggett, T. M., Schonmann, R. H. and Stacey, A. M. (1997). Domination by product measures. *Ann. Probab.* **25** 71–95. MR1428500

Mathieu, P. and Remy, E. (2004). Isoperimetry and heat kernel decay on percolation clusters. *Ann. Probab.* **32** 100–128. MR2040777

Nash, J. (1958). Continuity of solutions of parabolic and elliptic equations. *Amer. J. Math.* **80** 931–954. MR100158

Penrose, M. D. and Pisztora, A. (1996). Large deviations for discrete and continuous percolation. *Adv. in Appl. Probab.* **28** 29–52. MR1372330

Pisztora, A. (1996). Surface order large deviations for Ising, Potts and percolation models. *Probab. Theory Related Fields* **104** 427–466. MR1384040

Saloff-Coste, L. (1992). A note on Poincaré, Sobolev, and Harnack inequalities. *Internat. Math. Res. Notices* **2** 27–38. MR1150597

Saloff-Coste, L. (1997). Lectures on finite Markov chains. *Lectures on Probability Theory and Statistics. Ecole d'Éte de Probabilités de Saint-Flour XXVI. Lecture Notes in Math.* **1665** 301–408. Springer, Berlin. MR1490046

Saloff-Coste, L. and Stroock, D. W. (1991). Opérateurs uniformément sous-elliptiques sur les groupes de Lie. *J. Funct. Anal.* **98** 97–121. MR1111195

Sidoravicius, V. and Sznitman, A.-S. (2003). Quenched invariance principles for walks on clusters of percolation or amoung random conductances. Preprint. MR2063376

Stroock, D. W. and Zheng, W. (1997). Markov chain approximations to symmetric diffusions. *Ann. Inst. H. Poincaré Probab. Statist.* **33** 619–649. MR1473568

Thomassen, C. (1992). Isoperimetric inequalities and transient random walks on graphs. *Ann. Probab.* **20** 1592–1600. MR1175279

Department of Mathematics
University of British Columbia
Vancouver, British Columbia
Canada V6T 1Z2
e-mail: barlow@math.ubc.ca
url: www.math.ubc.ca/~barlow